\documentclass[11pt,a4paper]{article}

\usepackage{latexsym}
\usepackage{blindtext}

\PassOptionsToPackage{force}{filehook}

\usepackage[utf8]{inputenc}
\usepackage[english]{babel}
\usepackage{graphicx}
\usepackage{graphics}
\usepackage{fancyhdr}
\usepackage{lmodern}
\usepackage{color}
\usepackage{transparent}
\usepackage{float}
\usepackage{pythontex}
\usepackage{tikz-cd}
\usepackage{listings}
\usepackage[linktocpage=true,citebordercolor=green, urlcolor=black]{hyperref}
\usepackage{tabularx}
\usepackage[titletoc,title]{appendix}
\usepackage{xcolor}
\usepackage[mathscr]{eucal}
\usepackage{calrsfs}
\usepackage{cite}
\usepackage{comment}
\usepackage{subcaption}

\usetikzlibrary{shapes.misc}
\usepackage{pgf,tikz}
\usepackage{caption}
\usetikzlibrary{decorations.markings}
\usetikzlibrary{arrows}
\usetikzlibrary{patterns}
\usepackage{tikz-3dplot}
\usetikzlibrary{calc}
\usepackage{amsfonts}
\usepackage{amsmath,amssymb,amsthm}
\usepackage{mathtools}
\usepackage{relsize}
\usepackage{titlesec}
\usepackage{arydshln}
\usepackage{multirow}

\usepackage[most]{tcolorbox}
\newcolumntype{P}[1]{>{\raggedright\arraybackslash}p{#1}}

\usepackage{stmaryrd}
\usepackage{trimclip}

\usepackage{enumitem}

\usepackage{wasysym}

\makeatletter
\DeclareRobustCommand{\shortto}{%
	\mathrel{\mathpalette\short@to\relax}%
}

\newcommand{\short@to}[2]{%
	\mkern2mu
	\clipbox{{.5\width} 0 0 0}{$\m@th#1\vphantom{+}{\shortrightarrow}$}%
}
\makeatother

\theoremstyle{definition}

\newtheorem*{proposition*}{Proposition}

\titleformat{\paragraph}
{\normalfont\normalsize\bfseries}{\theparagraph}{1em}{}
\titlespacing*{\paragraph}
{0pt}{3.25ex plus 1ex minus .2ex}{1.5ex plus .2ex}

\renewcommand{\arraystretch}{1.3}

\usepackage[margin=2.5cm]{geometry}

\usepackage{mathptmx}

\AtEndDocument{%
	\par
	\medskip
	\begin{tabular}{@{}l@{}}%
		\textsc{Cengiz Aydin}\\
		\textsc{Institut für Mathematik, Universität Heidelberg, Germany}\\
		\textit{E-mail address}: \texttt{\href{mailto:cengiz.aydin@hotmail.de}{cengiz.aydin@hotmail.de}}
\end{tabular}}

\usepackage{authblk}
\title{Exploration of vertical self-resonant bifurcations from Distant Retrograde Orbits (DROs) in the Earth--Moon Circular Restricted Three-Body Problem (CR3BP)}
\author{Cengiz Aydin}

\begin{document}
	
\setcounter{page}{1}
\pagenumbering{arabic}
	
\maketitle
	
\begin{abstract}
This paper investigates vertical self-resonant bifurcations from the Distant Retrograde Orbit (DRO) family in the Earth--Moon circular restricted three-body problem.\ Using a predictor--corrector continuation algorithm, we identify fourteen vertical-critical DROs consisting of twelve higher-order periods of multiplicity from integer multiples of five to ten and two vertical pitchfork bifurcations.\ Their continuation generates twenty-five families of spatial symmetric periodic orbits.\ Beyond the computation of new periodic orbit families, the main contribution of this work is a systematic characterization of the resulting global bifurcation topology through symplectic invariants.\ To this end, we combine Floquet multiplier transitions at local bifurcations with Conley--Zehnder indices, Krein signatures and the Euler characteristic of the local Floer homology to develop a computational symplectic framework that systematically classifies the bifurcated families and provides a topological interpretation of their numerically established connections.\ The computations reveal a global bifurcation network generated by the DRO family, including bridge families that connect the DROs with Butterfly, planar prograde and P3DRO families, together with several families consisting of long-period multi-excursion orbits that approach collision with the Moon.\ The computations also identify spatial periodic orbits satisfying Earth--Moon resonances along several bifurcated families.\
\end{abstract}
	
\begin{center}
	\begin{tabular}{ll}
		\textbf{Keywords} & $\quad$ Earth--Moon CR3BP $\cdot$ DRO $\cdot$ vertical self-resonant bifurcations $\cdot$ bridge family $\cdot$ \\
		& $\quad$ bifurcation diagram $\cdot$ Conley--Zehnder index $\cdot$ spatial resonant orbit\\
		\textbf{MSC 2020} & $\quad$ 70H12 $\cdot$ 70F07 $\cdot$ 70G45
	\end{tabular}
\end{center}
	
\tableofcontents
	
\section{Introduction}

The Circular Restricted Three-Body Problem (CR3BP) is a fundamental model in celestial mechanics that describes the motion of an infinitesimal mass under the gravitational attraction of two primary bodies moving on circular orbits about their common barycenter.\ Despite its idealized nature, the CR3BP captures many of the qualitative dynamical features observed in planetary and satellite systems and provides the mathematical framework underlying numerous studies in astrodynamics, mission design and nonlinear dynamics \cite{koon}.\ In particular, periodic orbits constitute the backbone of the phase space structure, providing natural pathways for transfers, long-term station-keeping and exploration of cislunar space.\

Among the most important planar periodic solutions are the Distant Retrograde Orbits (DROs), which form a family of retrograde periodic orbits around the smaller primary.\ Their existence and stability have been studied extensively since the pioneering works of Hénon \cite{henon_1969,henon_1974}, Broucke \cite{broucke} and Benest \cite{benest_1974,benest_1975,benest_1977}, who investigated the corresponding retrograde orbit families in the Hill problem and the CR3BP for various mass ratios.\ These studies established the DRO family as one of the principal periodic orbit families of the CR3BP and revealed its remarkable dynamical stability over a wide range of energies.\

Bifurcations of periodic orbit families provide a mechanism through which new families of periodic solutions emerge and are therefore fundamental for understanding the architecture of periodic motions.\ In particular, higher-period bifurcations generated by self-resonances from planar periodic orbits have attracted considerable attention.\ Hénon \cite{henon_1970,henon_2003} first identified planar period-three bifurcations of retrograde orbits in the Hill problem, leading to the family denoted by $g_3$.\ The analogous family in the CR3BP has subsequently been studied under the names P3DRO (for \textit{period-3} DRO) \cite{capdevilla,li} or $\text{QDRO}_{(g3)}$ (for \textit{quasi DRO}) \cite{restrepo_russell}.\ More generally, higher-period bifurcations of retrograde orbits have been investigated in several three-body systems, including the Jupiter--Europa~\cite{lam_whiffen,lara_russell_villac}, Sun--Earth~\cite{demeyer_gurfil,stramacchia}, Earth--Moon~\cite{capdevilla,conte,ming_shijie,minghu,oshima}, Mars--Phobos~\cite{oshima} and Jupiter--Ganymede~\cite{li} problems.\

For the Earth--Moon CR3BP, the basic planar symmetric periodic orbit families, including retrograde, prograde and Lyapunov orbits, together with their stability properties, were characterized by Broucke~\cite{broucke}.\ Periodic solutions emanating from the collinear libration points and their bifurcations up to period-doubling were investigated in \cite{doedel}.\ Motions in the vicinity of Near-Rectilinear Halo Orbits, generated through period-doubling and period-quadrupling bifurcations, have been studied in \cite{zimovan}, while an extensive database of planar and spatial periodic orbit families near the Moon is available in \cite{franz_russell}.\ More recently, systematic numerical techniques for detecting and continuing bifurcations of periodic orbit families have been developed in \cite{caleb}.\ Despite these significant advances and the considerable enlargement of the catalogue of known periodic orbit families, much less attention has been devoted to their global bifurcation topology, that is, to understanding how different periodic orbit families are interconnected through local bifurcations.\ In particular, although vertical self-resonant bifurcations of the DRO family have been investigated in other CR3BP configurations and in the Hill problem~\cite{aydin_cz,aydin_batkhin,kalantonis,lara_russell_villac,robin_markellos}, a systematic study of the vertical self-resonant bifurcations of the Earth--Moon DRO family, together with the global relation of the bifurcated families, has not yet been carried out.\ This constitutes the primary motivation for the present work.\

Beyond the computation of new periodic orbit families, the key objective of this paper is to understand how these families are interconnected through bifurcations and to characterize their global organization.\ Numerical continuation establishes the connections between periodic orbit families, while symplectic invariants provide additional topological and dynamical information for their classification and interpretation.\ In particular, the Conley--Zehnder index, Krein signature, and Euler characteristic of the local Floer homology help constrain the possible bifurcation configurations and verify the consistency of the resulting global bifurcation diagrams.\ This approach has recently been applied to several restricted three-body problems.\ In the Hill problem, \cite{aydin_cz} investigated retrograde and prograde periodic orbit families and established fundamental bifurcation relations between them.\ These results were substantially extended in~\cite{aydin_batkhin}, where the global bifurcation network connecting retrograde, prograde, Lyapunov, halo and vertical rectilinear collision orbit families was characterized through symplectic invariants.\ The same framework was applied to the Saturn--Enceladus and Jupiter--Europa systems in \cite{moreno_aydin}.\ Very recently, \cite{joung} combined symplectic invariants with Moser regularization to investigate bifurcations of highly inclined near halo orbits.\ In the present work, we extend this line of research to the Earth--Moon CR3BP, where the richer geometry of the DRO family gives rise to a diverse collection of vertical self-resonant bifurcations, spatial bridge families and global bifurcation connections.\

The objectives of this paper are therefore twofold.\ First, we perform a systematic numerical investigation of vertical self-resonant bifurcations of the Earth--Moon DRO family and compute the bifurcated families of spatial symmetric periodic orbits.\ Second, we interpret the resulting bifurcations through symplectic invariants in order to classify the bifurcated families and construct bifurcation diagrams that characterize the bifurcation topology.\ This combined numerical and symplectic approach provides both new families of periodic solutions and a coherent framework for understanding their relationships.\

Our computations identify fourteen vertical-critical DROs, comprising twelve higher-order periods of multiplicity from integer multiples up to ten and two vertical pitchfork bifurcations, that generate twenty-five spatial symmetric periodic orbit families.\ The resulting bifurcation diagrams reveal a global bifurcation network connecting the DRO family with the Butterfly, planar prograde and P3DRO families, together with additional families of long-period multi-excursion orbits approaching collision with the Moon.\ Spatial Earth--Moon resonant periodic orbits are identified as well.\

The paper is structured as follows.\ Section \ref{sec2} reviews the CR3BP dynamics and its symmetry properties.\ Section \ref{sec3} presents the computational methodology, including the continuation algorithm for symmetric periodic orbits, Floquet stability analysis and the symplectic framework based on Conley--Zehnder indices, Krein signatures and Euler characteristic of the local Floer homology.\ Section \ref{sec4} presents the computed vertical self-resonant bifurcations of the Earth--Moon DRO family and characterizes their global bifurcation topology.\ Section~\ref{sec5} concludes the work, discusses its limitations and outlines directions for future research.\ Appendix~\ref{appendix_s} summarizes the predictor--corrector equations for spatial symmetric periodic orbits.\ Appendix~\ref{appendix} provides numerical database of the computed periodic orbit families.\

\section{The CR3BP and Symmetry Properties}
\label{sec2}

\begin{figure}[t!]
	\centering
	\includegraphics[scale=1]{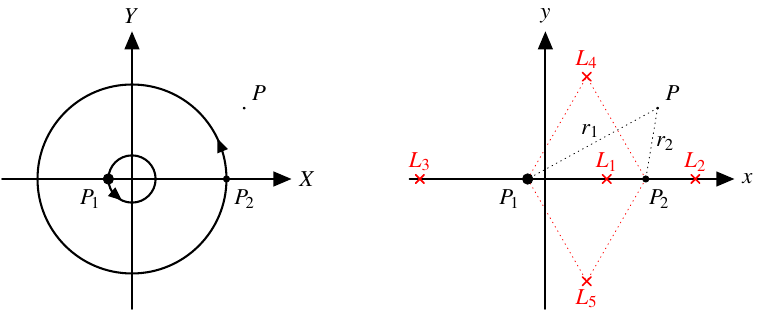}
	\caption{CR3BP in an inertial reference frame $(X,Y)$ and in a rotating reference frame $(x,y)$, with primaries $P_1$ and $P_2$, infinitesimal body $P$, libration points $L_i$ ($i=1,2,3,4,5$), and distances $r_1$ and $r_2$ from $P$ to the primaries $P_1$ and $P_2$, respectively.}
	\label{figure_1_cr3bp}
\end{figure}

In the spatial CR3BP we consider two point masses $P_1$ and $P_2$, which we respectively refer to as larger and smaller primary, moving along circles in the same plane around their common center of mass with constant angular velocity.\ A third body $P$, assumed to have negligible mass, moves under the gravitational attraction of the primaries without influencing their motion.\ Let $m_1$ and $m_2$ denote the masses of the larger and smaller primary, respectively, and define the mass parameter $\mu = m_2 / (m_1 + m_2)$.\ Throughout this paper we use the standard normalized units in which the gravitational constant is equal to one, the distance between the primaries is one, the total mass satisfies $m_1 + m_2 = 1$, and the period of the circular motion of the primaries is $2 \pi$.\ Consequently, the normalized mass of the large and smaller primary equals $1-\mu$ and $\mu$, respectively.\ Moreover, we consider a rotating frame of reference that rotates with an angular velocity of the orbital rate of the primaries.\ Choosing the origin of coordinates at the barycenter, the larger primary~$P_1$ is fixed at $(-\mu,0,0)$, while the smaller primary $P_2$ is fixed at $(1-\mu,0,0)$, as shown in Figure \ref{figure_1_cr3bp}.\ The dynamics of the infinitesimal body $P$ is described by the equations \cite[Chapter 10]{szebehely}:\
\begin{align}
	\ddot{x} &= 2 \dot{y} + x - (1 - \mu) \dfrac{x + \mu}{r_1^3} - \mu \dfrac{x - 1 + \mu}{r_2^3}, \nonumber \\
	\ddot{y} &= - 2 \dot{x} + y - \left( \dfrac{1 - \mu}{r_1^3} + \dfrac{\mu}{r_2^3} \right) y, \label{ham_equation} \\
	\ddot{z} &= - \left( \dfrac{1-\mu}{r_1^3} + \dfrac{\mu}{r_2^3} \right) z, \nonumber
\end{align}
where $r_1 = \left( (x+\mu)^2 + y^2 + z^2 \right)^{\frac{1}{2}}$ and $r_2 = \left( (x - 1 + \mu)^2 + y^2 + z^2 \right)^{\frac{1}{2}}$ indicate the distances from the infinitesimal body to the larger and smaller primary, respectively.\ Defining the canonical momenta by
\begin{align}\label{kinetic_moments}
	p_x = \dot{x} - y,\quad p_y = \dot{y} + x, \quad p_z = \dot{z},
\end{align}
the system can be written in Hamiltonian form with the corresponding Hamiltonian function
\begin{align} \label{hamiltonian}
	H(x,y,z,p_x,p_y,p_z) = \frac{1}{2} \left( p_x^2 + p_y^2 + p_z^2 \right) - \frac{1-\mu}{r_1} - \frac{\mu}{r_2} + p_x y - p_y x,
\end{align}
and its equations,
\begin{align} \label{equations}
	\left\{ \dot{k} = \frac{\partial H}{\partial p_k},\quad \dot{p}_k = - \frac{\partial H}{\partial k} \right\},\quad \text{with } k = x,y,z.
\end{align}

Since the Hamiltonian is conserved, periodic solutions form smooth one-parameter families that may be parameterized by the energy, or equivalently by the Jacobi constant $C = -2H$.\ The phase space is 6-dimensional and endowed with the standard symplectic form $\omega_0 = \sum dp_k \wedge dk$ ($k=x,y,z$).\ The flow of the Hamiltonian vector field $X_H$, which is uniquely determined by $dH(\cdot) = \omega_0( \cdot, X_H )$, is equivalent to the equations~\eqref{ham_equation} and \eqref{equations}.\ In the linear algebra expression of $\omega_0$,
$$ \omega_0 (\cdot, \cdot) = \langle J \cdot, \cdot \rangle,\quad \text{with } J = \begin{pmatrix}
	0 & I_3\\
	-I_3 & 0
\end{pmatrix}, $$
where the brackets denote the Euclidean inner product, we can represent the Hamiltonian vector field by $X_H = J \nabla H$, where the gradient, as usual, is defined with respect to the Euclidean inner product.\ With that, the equations~\eqref{ham_equation} and \eqref{equations} can be equivalently expressed by
\begin{align}\label{equation_J}
	\dot{ \bold X } = J \nabla H ( \bold X ),\quad \bold X \equiv (x,y,z,p_x,p_y,p_z).
\end{align}

The CR3BP possesses five equilibrium points, commonly referred to as \textit{libration points} or \textit{Lagrange points}.\ Three of them, denoted by $L_1$, $L_2$ and $L_3$, lie on the line joining the primaries, while $L_4$ and~$L_5$ form equilateral triangles with the primaries, as illustrated in Figure \ref{figure_1_cr3bp}.\ Explicit expressions for their coordinates may be found, for instance, in \cite[Chapter 4]{szebehely}.\

We define the following involutive maps acting on the phase space:\
\begin{align*}
	\rho_x \colon \mathbb{R}^6 \to \mathbb{R}^6,\quad &(x,y,z,p_x,p_y,p_z) \mapsto (x,-y,-z,-p_x,p_y,p_z), \\
	\rho_{xz} \colon \mathbb{R}^6 \to \mathbb{R}^6,\quad &(x,y,z,p_x,p_y,p_z) \mapsto (x,-y,z,-p_x,p_y,-p_z), \\
	\sigma \colon \mathbb{R}^6 \to \mathbb{R}^6,\quad &(x,y,z,p_x,p_y,p_z) \mapsto (x,y,-z,p_x,p_y,-p_z). 
\end{align*}
Geometrically, the map $\rho_x$ describes a half-turn about the $x$-axis, while the map $\rho_{xz}$ corresponds to reflection across the $xz$-plane.\ The involution $\sigma$ represents reflection with respect to the ecliptic plane $\{z = 0\}$.\ The momentum components transform accordingly so that each involution leaves the Hamiltonian~(\ref{hamiltonian}) invariant.\ Moreover, these transformations satisfy $\rho_x \circ \rho_{xz} = \rho_{xz} \circ \rho_x = \sigma$.\ Among these involutions, the map $\sigma$ is symplectic, while the maps $\rho_x$ and $\rho_{xz}$ are anti-symplectic.\ In the Hamiltonian setting \cite{lamb_roberts}, anti-symplectic involutions correspond to time-reversing symmetries, while symplectic involutions preserve the orientation of the Hamiltonian flow.\

For a map $\rho$ on the phase space, we denote its fixed point set by $\text{Fix}(\rho) = \{ z \in \mathbb{R}^6 \colon \rho(z) = z \}$.\ For the involutions introduced above, the corresponding fixed point sets are
$$ \text{Fix}(\rho_x) = \{ (x,0,0,0,p_y,p_z) \},\quad \text{Fix}(\rho_{xz}) = \{ (x,0,z,0,p_y,0) \},\quad \text{Fix}(\sigma) = \{ (x,y,0,p_x,p_y,0) \}. $$
In particular, the set $\text{Fix}(\sigma)$ corresponds to the planar CR3BP.\ Consequently, every planar periodic orbit is invariant under the reflection across the ecliptic plane $\{z = 0\}$.\ The fixed point sets $\text{Fix}(\rho_x)$ and $\text{Fix}(\rho_{xz})$ define symmetry sections that play a fundamental role in the computation of periodic orbits that are invariant under them.\ We distinguish four types of symmetric periodic orbits with the period $T$.\
\begin{itemize}
	\item[(I)] A periodic orbit is called \textit{simple symmetric with respect to $\rho_x$} if
	\begin{align}\label{bound_cond_1}
		y(0) = z(0) = p_x(0) = 0,\quad \text{ and } \quad y(T/2) = z(T/2) = p_x(T/2) = 0.
	\end{align}
	\item[(II)] A periodic orbit is called \textit{simple symmetric with respect to $\rho_{xz}$} if
	\begin{align}\label{bound_cond_2}
		y(0) = p_x(0) = p_z(0) = 0,\quad \text{ and } \quad y(T/2) = p_x(T/2) = p_z(T/2) = 0.
	\end{align}
	\item[(III)] A periodic orbit is called \textit{doubly symmetric with respect to $\rho_x$-$\rho_{xz}$} if
	\begin{align}\label{bound_cond_3}
		y(0) = z(0) = p_x(0) = 0,\quad \text{ and } \quad y(T/4) = p_x(T/4) = p_z(T/4) = 0.
	\end{align}
	\item[(IV)] A periodic orbit is called \textit{doubly symmetric with respect to $\rho_{xz}$-$\rho_x$} if
	\begin{align}\label{bound_cond_4}
		y(0) = p_x(0) = p_z(0) = 0,\quad \text{ and } \quad y(T/4) = z(T/4) = p_x(T/4) = 0.
	\end{align}
\end{itemize}
\noindent
Geometrically, the boundary conditions \eqref{bound_cond_1} and \eqref{bound_cond_2} imply that the orbit intersects the corresponding symmetry section perpendicularly at half period, while \eqref{bound_cond_3} and \eqref{bound_cond_4} imply that the orbit possesses two distinct symmetry sections, which are intersected by the orbit alternately at quarter period.\

\section{Computational Methods}
\label{sec3}

\subsection{Predictor--Corrector Method for Symmetric Periodic Orbits}

For the accurate computation and continuation of symmetric periodic orbits we use a classical predictor--corrector algorithm based on \cite{robin_markellos}.\ Closely related formulations can also be found in \cite{kazantzis,howell}.\ The original formulations were developed in non-canonical Lagrangian coordinates, while in the present work we apply their canonical Hamiltonian counterparts, obtained through the linear transformation \eqref{kinetic_moments} between velocities and canonical momenta.\ The predictor--corrector method requires the first-order variation of the solution with respect to the initial conditions that is described by the state-transition matrix.\ The state-transition matrix, which we denote by
\begin{align}\label{stm}
	V(t) = (v_{ij}) \in \mathbb{R}^{6 \times 6},
\end{align}
solves the variational equations,
\begin{align*}
	\dot{V}(t) = J D^2 H \left( \bold X (t) \right) V(t),\quad V(0) = \text{id},
\end{align*}
that have to be integrated simultaneously with the equations \eqref{equation_J} along a solution $\bold X \equiv (x,y,z,p_x,p_y,p_z)$, where $D^2 H $ is the Hessian of the Hamiltonian.\vspace{0.5em}\\
\textbf{Planar symmetric periodic orbits.}\ Planar motion is confined to the invariant manifold $\{ z = p_z = 0 \}$.\ We restrict the discussion to the simple case in which a planar symmetric periodic orbit possesses only two perpendicular intersections with the $x$-axis per period.\ To start integrating the equations \eqref{equation_J} from the $x$-axis, we first have to choose a decent initial guess in the form $\bold X (0) \equiv (x_0,0,0,0,p_{y0},0)$.\ The symmetric orbit will be periodic with the period $T$ if the final state vector $\bold X (t) \equiv (x_{cut},y_{cut},0,p_{x,cut},p_{y,cut},0)$ at the time $t = T/2$ fulfills the periodicity conditions:
\begin{align}
	y_{cut} = y(x_0,0,0,0,p_{y0},0;T/2) = 0, \label{per_cond_1}\\
	p_{x,cut} = p_x(x_0,0,0,0,p_{y0},0;T/2) = 0. \label{per_cond_2}
\end{align}
In the numerical implementation, the periodicity conditions \eqref{per_cond_1} and \eqref{per_cond_2} are considered to be satisfied when $|y_{cut}| < \varepsilon$ and $|p_{x,cut}| < \varepsilon$, where $\varepsilon$ is a sufficiently small prescribed numerical tolerance, which is $\varepsilon = 10^{-10}$ in this work.\ In practice, we stop the integration when we cross the $x$-axis for the first time, such that the first condition~\eqref{per_cond_1} is automatically satisfied.\ At this point, the time is defined to be $t = T/2$.\ Since the second condition~\eqref{per_cond_2} is in general only approximately satisfied, we seek appropriate corrections $\delta x_0, \delta p_{y0}, \delta t$.\ For this purpose, we consider the first-order Taylor expansion of \eqref{per_cond_1} and \eqref{per_cond_2}, and obtain
\begin{align}
	v_{21} \delta x_0 + v_{25} \delta p_{y0} + (p_{y,cut}(t) - x_{cut}(t)) \delta t = 0,\label{corr_1}\\
	v_{41} \delta x_0 + v_{45} \delta p_{y0} + \dot{p}_{x,cut}(t) \delta t = -p_{x,cut},\label{corr_2}
\end{align}
where $v_{21}$, $v_{25}$, $v_{41}$ and $v_{45}$ are corresponding entries of the state-transition matrix \eqref{stm}.\ By fixing one of the two initial conditions, we solve the equations \eqref{corr_1} and \eqref{corr_2} for the remaining one, and obtain the corrector equations:\
\begin{itemize}
	\item[a)] If we keep $x_0$ fixed, then $ \delta p_{y0} = - p_{x,cut} / u_{45}$,
	\item[b)] and holding $p_{y0}$ fixed yields $ \delta x_0 = - p_{x,cut} / u_{41}$,
\end{itemize}
where $u_{41}$ and $u_{45}$ denote the corrected variational coefficients defined by
$$ u_{41} = v_{41} - v_{21} \frac{\dot{p}_{x,cut}}{p_{y,cut} - x_{cut}},\quad u_{45} = v_{45} - v_{25} \frac{\dot{p}_{x,cut}}{p_{y,cut} - x_{cut}}. $$
With the corrected initial data, we repeat the procedure in an iterative way of correction steps until the numerical tolerance is satisfied.\ After determining a symmetric periodic orbit, one initial condition serves as the continuation parameter, while the other is predicted from \eqref{corr_1} and \eqref{corr_2}.\ Here, $p_{x,cut} = 0$, and we write $\delta x_0 \equiv \Delta x_0$ and $\delta p_{y0} \equiv \Delta p_{y0}$ to distinguish between correction and prediction steps:\
\begin{itemize}
	\item[a)] If $x_0$ is selected as the continuation parameter and varied by an arbitrarily chosen step size~$\Delta x_0$ (e.g.,\ $5 \cdot 10^{-5}$), then the corresponding variation of the initial state vector is of the form
	$$ \bold X (0) + \Delta \bold X (0) = \left(x_0 + \Delta x_0 ,0,0,0, p_{y0} - \frac{ u_{41} }{ u_{45} } \Delta x_0 ,0 \right). $$
	\item[b)] Similarly, when $p_{y0}$ is the continuation parameter and varied by an arbitrarily chosen step size~$\Delta p_{y0}$, then the corresponding variation of the initial state vector is of the form
	$$ \bold X (0) + \Delta \bold X (0) = \left(x_0 - \frac{ u_{45} }{ u_{41} } \Delta p_{y0} ,0,0,0, p_{y0} + \Delta p_{y0} ,0 \right). $$
\end{itemize}
\noindent
These predicted initial conditions shall be corrected respectively until the periodicity condition is satisfied.\ The choice of which of the parameters to fix has an important effect on the convergence of the solutions.\ If $x_0$ is used as the continuation parameter, both the prediction and the correction steps require division by $u_{45}$.\ Similarly, the continuation in $p_{y0}$ requires division by $u_{41}$.\ The continuation parameter is chosen according to the corrected variational coefficient of largest absolute value.\ If the continuation parameter reaches an extremum, the corresponding corrected variational coefficient tends to zero, indicating that the continuation should be switched to the other parameter.\ \vspace{0.5em}\\
\textbf{Spatial symmetric periodic orbits.}\ The computation and continuation of spatial symmetric periodic orbits follows the same predictor--corrector technique described above for planar symmetric periodic orbits.\ In view of the four symmetry types \eqref{bound_cond_1}--\eqref{bound_cond_4}, a spatial symmetric periodic orbit possesses three initial conditions.\ Since in practice we always stop the integration when crossing the section $y=0$, the periodicity condition $y_{cut}$ is satisfied automatically, and only the remaining two periodicity conditions must be enforced.\ We therefore select one initial condition as the continuation parameter, while the remaining two are determined by the linearized periodicity equations.\ For completeness, a summary of the corresponding initial and final states together with the predictor--corrector equations for all four symmetry types is provided in Appendix~\ref{appendix_s} (see Tables~\ref{tab:spatial_summary_1} and~\ref{tab:spatial_summary_2}).\

\subsection{Monodromy Matrix, Floquet Stability and Local Bifurcation Detection}

Once a $T$-periodic orbit has converged, its linear stability properties are determined from the monodromy matrix $M = V(T)$, where $V(T)$ is the state-transition matrix \eqref{stm} evaluated over one period.\ For symmetric periodic orbits, the anti-symplectic symmetries $\rho_x$ and $\rho_{xz}$ allow the monodromy matrix to be constructed from the state-transition matrix over one half- or one quarter-period.\ Accordingly~\cite{bray_goudas}:
\begin{align*}
	\text{(I) simple symmetric with respect to $\rho_x$:}\quad & M = L_X V^{-1}(T/2) L_X V(T/2),\\
	\text{(II) simple symmetric with respect to $\rho_{xz}$:}\quad & M = L_{XZ} V^{-1}(T/2) L_{XZ} V(T/2), \\
	\text{(III) doubly symmetric with respect to $\rho_x$-$\rho_{xz}$:}\quad & M = \left[ L_{X} V^{-1}(T/4) L_{XZ} V(T/4) \right]^2,\\
	\text{(IV) doubly symmetric with respect to $\rho_{xz}$-$\rho_x$:}\quad & M = \left[ L_{XZ} V^{-1}(T/4) L_{X} V(T/4) \right]^2,
\end{align*}
where $L_{X}$ and $L_{XZ}$ denote the matrix representations of $\rho_x$ and $\rho_{xz}$, respectively.\ This substantially reduces the computational effort.\ Since the monodromy matrix is a $6 \times 6$ symplectic matrix, the eigenvalues occur in reciprocal pairs and are of the form
$$\{1,1,\lambda_1,\lambda_1^{-1},\lambda_2,\lambda_2^{-1}\},$$
where the double unit eigenvalue follows from the autonomous Hamiltonian structure of the CR3BP, reflecting both time invariance and energy conservation.\ The non-trivial eigenvalues are referred to as the \textit{Floquet multipliers}.\ For each pair of Floquet multipliers the \textit{stability index} is defined as
$$ s_k = \frac{1}{2} \left( \lambda_k + \lambda_k^{-1} \right),\quad k=1,2. $$
These can be computed via the following formulas \cite{bray_goudas,broucke_stab}:
$$ s_1 = - \frac{1}{4} \left( \alpha + \sqrt{\beta} \right),\quad s_2 = - \frac{1}{4} \left( \alpha - \sqrt{\beta} \right), $$
where $\alpha = 2 - \text{tr}(M)$ and $\beta  = 2 \text{tr}(M^2) - \alpha^2 + 4$.\ Each stability index is classified as:
\begin{itemize}[noitemsep]
	\item \textit{Elliptic:}\ $-1 < s_k = \cos \theta < 1$.\ It is equivalent to $\lambda_k = e^{i \theta}$.
	\item \textit{Positive hyperbolic:}\ $s_k > 1$.\ It is equivalent to $\lambda_k \in \mathbb{R}_{>0} \setminus \{1\}$.
	\item \textit{Negative hyperbolic:}\ $s_k < -1$.\ It is equivalent to $\lambda_k \in \mathbb{R}_{<0} \setminus \{-1\}$.
	\item \textit{Complex instability:}\ Case of quadruples $\{ \lambda, \lambda^{-1}, \overline{\lambda}, \overline{\lambda}^{-1} \}$.
\end{itemize}
During numerical continuation, the stability indices are evaluated for all computed periodic orbits.\ Periodic orbits are considered as (linearly) stable when both stability indices satisfy $|s_1| < 1$ and $|s_2| < 1$.\ By virtue of the $\sigma$-symmetry, the two stability indices of planar periodic orbits naturally separate into horizontal and vertical components.\ The vertical stability index $s_v$ plays a central role in the computation of bifurcated spatial periodic orbit families.\

Candidate bifurcation points are detected whenever the Floquet multipliers approach the critical values $\pm 1$ or satisfy the self-resonance condition (see below).\ Throughout this work, we refer to the Floquet multipliers associated with the bifurcation as the \textit{critical Floquet multipliers}.\ The precise bifurcation type is then determined from the corresponding stability index, critical Floquet multiplier configuration and local continuation behavior.\ Whenever a bifurcation point is approached during continuation, the continuation step size is successively reduced to determine accurately its location.\ The bifurcation types relevant to this work are summarized below (and illustrated in Figure~\ref{fig:bifurcations}).\
\begin{itemize}
	\item[1)] \textbf{Pitchfork (tangent) bifurcation \cite{campbell,bosanac}.}\ A pitchfork (tangent) bifurcation is characterized by the crossing of one of the stability indices through 1, corresponding to a unit Floquet multiplier, while the energy (or equivalently the Jacobi constant) varies monotonically along the continued family.\ This crossing is accompanied by a change in the stability of the periodic orbits and the branching of two additional periodic orbit families that possess the same stability as members of the original family prior to the bifurcation.\
	\item[2)] \textbf{Fold bifurcation \cite{campbell,bosanac,sadovskii}.}\ A fold (or saddle--node, or birth--death) bifurcation is also identified by one of the stability indices crossing $s_k = 1$, corresponding to a unit Floquet multiplier.\ Unlike a pitchfork bifurcation, the energy (or equivalently the Jacobi constant) attains a local maximum or minimum at the bifurcation point, where two families of periodic orbits merge.\ As the bifurcation is approached, the critical Floquet multipliers of one family converge to unity along the unit circle, while those of the other approach unity along the positive real axis, reflecting the different stability characteristics of the two families.\
	\item[3)] \textbf{Period-multiplying (or self-resonant) bifurcations \cite{campbell,bosanac,sadovskii}.}\ A period-multiplying bifurcation occurs when a pair of elliptic Floquet multipliers satisfies the self-resonance condition
	$$\lambda_k = e^{\pm 2 \pi i p / q},$$
	where $p,q \in \mathbb{Z}$ with $p < q$ and $p / q$ irreducible fraction.\ Equivalently, the stability index satisfies
	\begin{align}\label{vsr}
		s_k = \cos \left( 2 \pi \frac{p}{q} \right).
	\end{align}
	The corresponding $q$-fold covering has monodromy matrix $M^q$, where $M$ is the monodromy matrix of the single-turn orbit.\ Consequently, each Floquet multiplier is raised to the $q$-th power, i.e.,\ $\lambda_k^q$.\ In particular, the critical resonant Floquet multipliers become unit Floquet multipliers for the $q$-fold covering, allowing the bifurcation of higher-period orbits.\ The case $q = 2$ corresponds to the period-doubling bifurcation $(\lambda_k = s_k = -1)$.\ Higher-order resonances $(q>2)$ generate self-resonant bifurcations, from which new periodic orbit families may emanate with period $q$ times the period of the underlying family.\ Throughout this work, we denote a self-resonance bifurcation point by its self-resonance ratio $p$:$q$ (or $(q-p)$:$q$), where $q$ is referred to as the \textit{multiplicity of the self-resonance}.\ Furthermore, particular attention is given to the \textbf{Vertical Self-Resonant (VSR) bifurcations}~\cite{robin_markellos}, where the self-resonance occurs within the vertical Floquet multipliers of planar periodic orbits.\ These vertical bifurcations generate families of spatial symmetric periodic orbits, with simple or doubly symmetric orbits depending on the parity of the multiplicity $q$ \cite{robin_markellos}:\ For odd $q$, the bifurcated families consist of simple symmetric orbits, with one family symmetric with respect to $\rho_x$ and the other with respect to $\rho_{xz}$.\ By the $\sigma$-symmetry, each family possesses a symmetric counterpart and therefore exists as a pair of equivalent families.\ For even $q$, both bifurcated families consist of doubly symmetric orbits.\ In all cases, the two bifurcated families have different stability properties.\ Along one family, the critical Floquet multipliers remain on the unit circle, while along the other family, the critical Floquet multipliers move onto the positive real axis.\
	\item[4)] \textbf{Touch-and-go bifurcations \cite{sadovskii}.}\ A special case of the period-multiplying bifurcation may occur for $q = 3,4$ that is known as a touch-and-go bifurcation.\ In this case, two period-tripling (or period-quadrupling) periodic orbit families bifurcate from the original family while the critical Floquet multipliers remain positive hyperbolic on both families.\ The two bifurcated families continue in opposite directions with respect to the energy (or equivalently the Jacobi constant), giving the characteristic ``touch-and-go'' continuation behavior.\ These bifurcations arise in the planar \cite{henon_1970,henon_2003} or spatial \cite{walawska,joung} dynamics.\ In the cases of the VSR bifurcations and the vertical rectilinear collision orbits in the Hill problem~\cite{aydin_batkhin}, touch-and-go bifurcations do not occur.\
	\item[5)] \textbf{Secondary Hopf and modified secondary Hopf bifurcations \cite{campbell}}.\ A secondary Hopf bifurcation occurs when two Floquet multipliers collide on the unit circle at a point other than $\pm 1$ and subsequently leave the unit circle as a complex quadruple, resulting in a change of stability.\ Such a bifurcation may occur when the colliding Floquet multipliers have opposite Krein signatures (see below).\ In contrast, a modified secondary Hopf bifurcation is determined by the collision of two Floquet multipliers on the real axis, away from $\pm 1$, followed by their departure from the real axis as a complex quadruple.\ However, no change in stability occurs.\ Since Hopf bifurcations are not the focus of the present study, their detailed analysis is beyond the scope of this work.\
\end{itemize}

\subsection{Conley--Zehnder Index and Krein Signature}

The Conley--Zehnder (CZ) index~\cite{conley_zehnder} is a symplectic invariant associated with a non-degenerate periodic orbit, i.e., a periodic orbit whose Floquet multipliers do not contain the unit eigenvalue.\ Geometrically, it measures the winding of the linearized Hamiltonian flow relative to a fixed symplectic reference frame and remains constant under continuation as long as the orbit does not pass through a bifurcation.\ Therefore, the CZ index changes only when the orbit becomes degenerate.\ In the formal definition, the CZ index is defined through Maslov index theory \cite{robbin_salamon,salamon_zehnder} as the intersection number of a path of symplectic matrices with the Maslov cycle (that consists of all symplectic matrices having the unit eigenvalue), whose direct computation requires constructing a symplectic trivialization and evaluating the associated crossing forms~\cite{moreno_aydin}.\ Instead, we compute CZ indices by a continuation-based approach following~\cite{aydin_cz,aydin_batkhin}:\ Starting from periodic orbit families whose CZ indices are known analytically, the CZ indices are computed along the numerical continuation by monitoring the Floquet multipliers together with Krein signatures.\ The required formulas and constructions for the CZ index and the Krein signature are summarized below.\

Let $\bold X$ denote a $T$-periodic planar orbit.\ Due to the $\sigma$-symmetry, the CZ index of $\bold X$ splits additively into
$$\mu_{CZ}(\bold X) = \mu_{CZ}^p (\bold X) + \mu_{CZ}^s (\bold X),$$
where $\mu_{CZ}^p$ and $\mu_{CZ}^s$ denote the planar and spatial CZ indices, respectively.\ We denote the $q$-fold covering of $\bold X$ by $\bold X^q$ whose period equals $qT$, and assume that $\bold X^q$ is non-degenerate for all $q \geq 1$.\ Then, the CZ indices of $\bold X^q$ are expressed in terms of the Floquet multiplier type of the underlying orbit $\bold X$ via the iteration formulas provided by Hofer--Wysocki--Zehnder~\cite{hofer_w_z}.\ For elliptic Floquet multiplier pairs, the CZ index is determined by the accumulated rotation angles.\ Formally, denoting the accumulated planar and spatial rotation angles of $\bold X$ by $\varphi_p$ and $\varphi_s$, respectively, the CZ indices of $\bold X^q$ are given by
\begin{align}\label{indices}
	\mu_{CZ}^p (\mathbf X^q) = 1 + 2 \cdot \lfloor q \cdot \varphi_p /(2\pi) \rfloor, \quad \mu_{CZ}^s (\mathbf X^q) = 1 + 2 \cdot \lfloor q \cdot \varphi_s /(2\pi) \rfloor.
\end{align}
The floor function counts the winding numbers (the number of complete revolutions) performed by the pair of Floquet multipliers on the unit circle during the $q$-fold covering.\ In particular, elliptic Floquet multipliers always contribute an odd CZ index which increases by two whenever the accumulated rotation angle exceeds another multiple of $2 \pi$.\ For hyperbolic Floquet multipliers, the CZ index depends only on the hyperbolic type and the covering number, with even and odd indices corresponding to positive and negative hyperbolic Floquet multipliers, respectively:\
\begin{align}\label{indices_2}
	\mu_{CZ}^{p/s}(\mathbf X^q) = q n,\quad n \in \begin{cases}
		2 \mathbb{Z} & \text{ for the pos. hyperbolic case},\\
		2 \mathbb{Z} + 1 & \text{ for the neg. hyperbolic case.}
	\end{cases}
\end{align}

\begin{figure}[t!]
	\centering
	\includegraphics[scale=1]{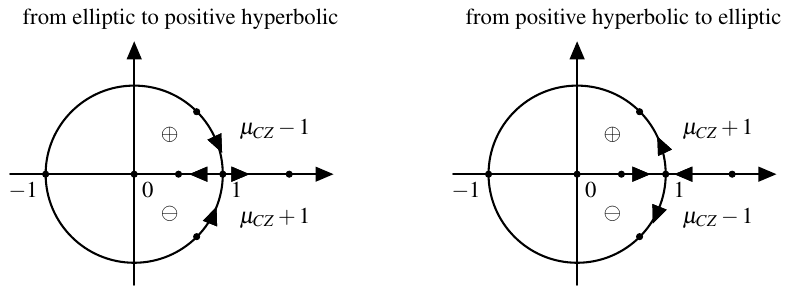}
	\caption{The index jump (reproduced from \cite{aydin_cz} with minor modifications).\ Left:\ When eigenvalue 1 is crossed from above (or below), the index goes down (or up) by 1.\ Right:\ When Floquet multipliers after crossing eigenvalue 1 goes up (or down), the index goes up (or down) by 1.\ The direction of crossing is determined by the $\pm$ Krein signature.}
	\label{figure_2}
\end{figure}

When a pair of Floquet multipliers changes from elliptic to positive hyperbolic, or vice versa, the corresponding CZ index changes by $\pm 1$, depending on the direction in which the Floquet multipliers cross the unit eigenvalue, as illustrated in Figure~\ref{figure_2}.\ If the Floquet multipliers only touch the unit eigenvalue without changing their stability type, the index changes by $\pm 2$, where the sign is determined by the direction of the rotation angles.\ To determine the direction of rotation of elliptic Floquet multipliers, we take into account the \textbf{Krein signature}, a symplectic invariant associated with each pair of elliptic eigenvalues~\cite[Appendix~29]{arnold_avez}.\ Let~$M$ be a symplectic matrix, and let $\lambda = e^{i \theta}$ and $\overline{\lambda} = e^{-i \theta}$ be a pair of elliptic eigenvalues of $M$ with $\lambda \neq \overline{\lambda}$.\ If $v$ is an eigenvector corresponding to $\lambda$, the real 2-dimensional plane $P := \langle \text{Re}(v), \text{Im}(v) \rangle $ is equipped with the quadratic form
\begin{align*}
	Q(\xi,\xi) := \omega_0 \left( \xi, M \xi \right),\quad \xi \in P.
\end{align*}
The quadratic form $Q$ is non-degenerate and either positive or negative definite.\ The \textbf{Krein signature} of the elliptic pair $\{\lambda,\overline{\lambda}\}$ is defined as the sign of the signature of the quadratic form $Q$, which is independent of the choice of the eigenvector.\ The Krein signature uniquely determines the orientation of the rotation on the unit circle.\ A positive Krein signature corresponds to the rotation angle $\theta \in (0,\pi)$, while a negative Krein signature corresponds to the complementary angle $- \theta$ in the range $(\pi,2\pi)$.\ Therefore, the Krein signature specifies the direction of the CZ index jump at bifurcation points, as illustrated in Figure~\ref{figure_2}.\ Furthermore, for VSR bifurcations the Krein signature distinguishes between the self-resonance ratios $p$:$q$ and $(q-p)$:$q$.\ We further notice that a fundamental result of Krein states that the collision of two pairs of elliptic Floquet multipliers with opposite Krein signatures can provoke complex instability, while the collision with the same Krein signatures does not cause a move off the unit circle.\ This criterion underlies the occurrence of secondary Hopf bifurcations discussed previously.\ It is also worthwhile to notice that this combination of CZ indices, Floquet multipliers and Krein signatures have been used to compute numerically the lunar periods known from astronomical observations \cite{aydin_babylonian}.\

\subsection{Conley--Zehnder Indices at Local Bifurcations}

Although Floquet multipliers identify the occurrence and type of local bifurcations, they do not provide a topological classification of the emerging periodic orbit families or explain how these families fit into the global bifurcation structure.\ From a topological perspective, the CZ index defines the integer grading of the local Floer homology whose underlying chain complexes are generated by periodic orbits.\ Since the local Floer homology is invariant under bifurcations~\cite{ginzburg}, its Euler characteristic, the alternating sum of the ranks of the homology groups, is invariant as well.\ Therefore, the Euler characteristic provides a topological relation between the CZ indices of all periodic orbit families involved in a local bifurcation, such as approached in~\cite{aydin_cz,aydin_batkhin,moreno_aydin}.\ These properties make the CZ index a natural topological label for periodic orbit families.\

When studying higher-period periodic orbits, it is necessary to distinguish between \textit{good} and \textit{bad orbits} according to \cite{eliashberg_givental_hofer}.\ The $q$-fold covering $\bold X^q$ of a planar orbit $\bold X$ is called \textit{good} if
\begin{align} \label{index_parity}
	\mu_{CZ}^p(\mathbf X ^q) \equiv \mu_{CZ}^p(\mathbf X) \text{ mod } 2, \quad \mu_{CZ}^s(\mathbf X^q) \equiv \mu_{CZ}^s(\mathbf X) \text{ mod } 2,
\end{align}
or both equations in \eqref{index_parity} are not simultaneously satisfied.\ Otherwise, the $q$-fold covering $\mathbf X^q$ is a \textit{bad orbit}.\ Only good periodic orbits contribute to the local Floer homology and its Euler characteristic.\ While all single-turn periodic orbits are good, bad orbits occur as $q$-fold covering whose exactly one pair of Floquet multipliers is of negative hyperbolic type, where~$q$ is even.\

For spatial periodic orbits, the iteration formulas~\eqref{indices} and~\eqref{indices_2}, and the equations~\eqref{index_parity}, are no longer directly applicable.\ Nevertheless, whenever a spatial periodic orbit family bifurcates from a planar periodic orbit family, or whenever its CZ index is known analytically, the CZ indices can be propagated along the continuation by monitoring each Floquet multiplier pair together with the Krein signatures.\

By combining the Floquet multiplier transitions at local bifurcation with Krein signatures, iteration formulas~\eqref{indices} and~\eqref{indices_2}, and the good--bad orbit classification, we provide explicit rules for assigning CZ indices of the continued and bifurcated families.\ These rules are summarized in Figure~\ref{fig:bifurcations}.\ For each local bifurcation, Figure~\ref{fig:bifurcations} specifies the Floquet multiplier transitions, the Krein signature of the critical elliptic Floquet multiplier pair, and the resulting CZ indices of all participating periodic orbit families.\ The corresponding local bifurcation configurations form the basis for constructing the global bifurcation diagrams presented in Section~\ref{sec4}.\ The assigned CZ indices verify the invariance of the local Euler characteristic:\
\begin{itemize}
	\item As a representative example, consider the first pitchfork bifurcation shown in Figure~\ref{fig:bif_pitchfork}.\ The incoming family has CZ index $k$.\ After the bifurcation, the continued family has CZ index $k-1$ due to the positive Krein signature.\ The two bifurcated families have the same stability as the orbits of the incoming family prior to the bifurcation, and thus have the same CZ index $k$.\ Counting the periodic orbits with CZ indices before and after the bifurcation point yields the equation
	$$ (-1)^k = (-1)^{k-1} + 2 \cdot (-1)^k, $$
	which verifies the invariance of the local Euler characteristic.\
	\item For the first period-doubling bifurcation in Figure~\ref{fig:bif_pd}, the continued period-doubled orbits are bad and therefore do not contribute.\ Therefore, the Euler characteristic before and after the bifurcation remains equal to $(-1)^k$.\
	\item Similarly, for the first VSR bifurcation shown in Figure~\ref{fig:bif_vsr}, the Euler characteristic relation reads
	$$(-1)^k = (-1)^{k-2} + (-1)^{k-1} + (-1)^k,$$
	where one representative of each bifurcated family is counted.\ For odd covering multiplicities~$q$, each spatial family has a $\sigma$-symmetric counterpart with the same CZ index.\ Counting both counterparts doubles the contributions of the two spatial families, which still cancel because their CZ indices differ by one.\ Thus, the Euler characteristic remains invariant.
\end{itemize}

The remaining configurations shown in Figure~\ref{fig:bifurcations} are verified analogously.\

\begin{figure}[t!]
	\centering
	\begin{subfigure}{\textwidth}
		\centering
		\includegraphics[width=\textwidth]{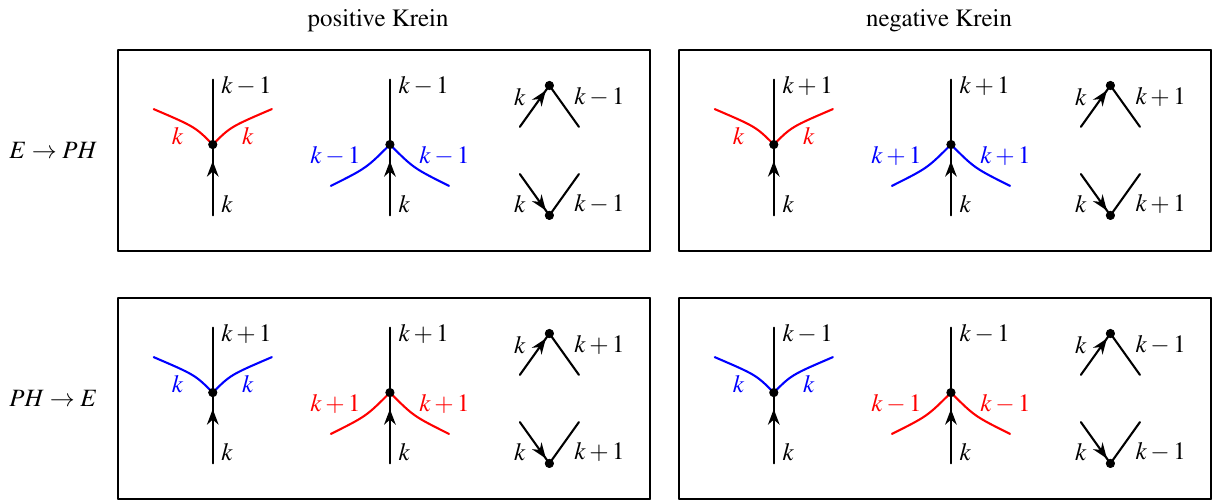}
		\caption{Pitchfork and birth--death bifurcations.\ The four boxes classify the local transitions according to the Floquet multiplier transition ($E \leftrightarrow$ PH) and the Krein signature of the elliptic Floquet multiplier pair.}
		\label{fig:bif_pitchfork}
	\end{subfigure}\vspace{1em}\\
	
	\begin{subfigure}{\textwidth}
		\centering
		\includegraphics[scale=0.8]{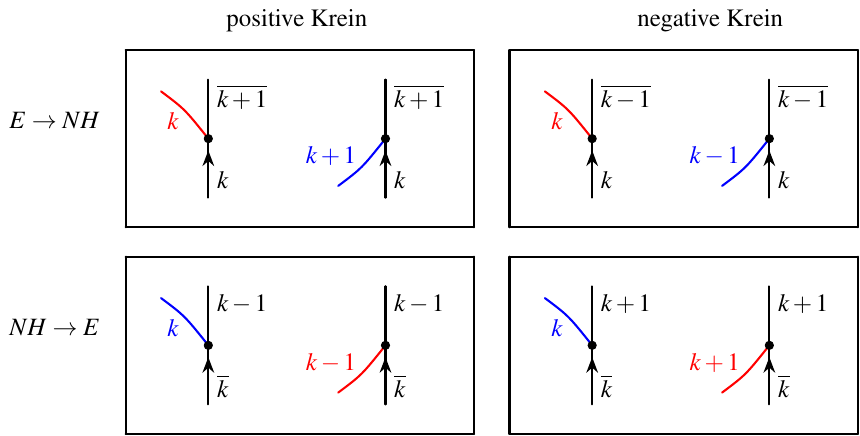}
		\caption{Period-doubling bifurcations.\ The four boxes classify the local transitions $(E \leftrightarrow NH)$ according to the Krein signature of the elliptic Floquet multiplier pair of the underlying single-turn orbit.\ It is assumed that no other Floquet multiplier pair is negative hyperbolic, so that the continued period-doubled orbits are bad; their Conley--Zehnder indices are denoted by an overline.}
		\label{fig:bif_pd}
	\end{subfigure}\vspace{1em}\\
	
	\begin{subfigure}{\textwidth}
		\centering
		\includegraphics[width=\textwidth]{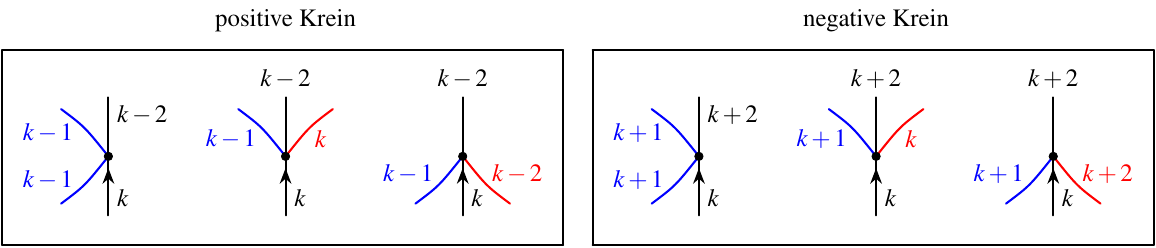}
		\caption{Touch-and-go and VSR bifurcations for positive (left box) and negative (right box) Krein signatures.\ In each box, the left diagram shows a touch-and-go bifurcation, while the middle and right show elliptic and positive hyperbolic VSR families, respectively.\ For odd $q$, only one representative of each $\sigma$-symmetric pair is shown.}
		\label{fig:bif_vsr}
	\end{subfigure}
	\caption{Conley--Zehnder indices at local bifurcations.\ $E$, $PH$ and $NH$ denote elliptic, positive and negative hyperbolic, respectively.\ Arrows indicate the continuation direction with respect to the energy (or equivalently the Jacobi constant).\ The labels denote the Conley--Zehnder indices of the continued and bifurcated families, where $k$ is the index of the incoming family.\ Red and blue families represent bifurcated families whose critical Floquet multipliers are elliptic and positive hyperbolic, respectively.}
	\label{fig:bifurcations}
\end{figure}

\section{Numerical Results}
\label{sec4}

In this section we present the families of spatial symmetric periodic orbits obtained from the continuation of Vertical Self-Resonant (VSR) bifurcations of the Distant Retrograde Orbit (DRO) family in the Earth--Moon CR3BP.\ The computations were carried out using the computation and continuation procedure together with the symplectic classification described in Section \ref{sec3}.\ The equations of motion and the associated variational equations were numerically integrated using the explicit twelfth-order Runge--Kutta method from \cite{feagin} that we have implemented in Python, accelerated with the Numba package.\ Appendix~\ref{appendix} contains the database of all computed periodic solutions that satisfy the periodicity conditions with the numerical tolerance of $\varepsilon = 10^{-10}$.\ All computations were performed using the Earth--Moon mass ratio $\mu = 1/82.27 \approx 0.0121550991$ and the normalized units introduced in Section \ref{sec2}.\

The main outcome is a collection of bifurcation diagrams (shown in Figures~\ref{figure_bif_1}, \ref{figure_bif_2}, \ref{plot_6_b}, and~\ref{figure_bif_4}) describing how the spatial families generated from the DROs are interconnected.\ Starting from each detected VSR bifurcation, all emerging families are initialized and continued until they terminate at another known orbit family or approach collision with one of the primaries.\ In the bifurcation diagrams, that are formed with the techniques shown in Figure \ref{fig:bifurcations}, vertices represent local bifurcation points and edges correspond to periodic orbit families labeled by their CZ indices.\ The diagrams are drawn schematically upward in the direction of decreasing Jacobi constant $C$, with vertical axes not to scale.\ Planar families are shown in black, while bifurcated families are distinguished by color (red and blue, as in Figure~\ref{fig:bifurcations}).\ Thick and thin edges denote doubly and simply symmetric orbit families, respectively, while dashed edges indicate the corresponding $\sigma$-symmetric counterparts.\ Overlined indices identify bad orbits.\ Birth--death bifurcations are labeled by ``b--d'', and the symbol ``$\times_M$'' denotes approaching collision with the Moon.\ Along each bifurcated family, we also identify $M$:$N$ Earth--Moon resonant orbits, where $M$ corresponds to the number of revolutions on the orbit and $N$ to the number of revolutions of the Earth and the Moon.\

\subsection{Distant Retrograde Orbit (DRO) Family}

In the vicinity of the Moon, the CR3BP is well approximated by the Kepler problem, whose planar circular solutions are divided into prograde and retrograde motions, rotating with and against the rotating reference frame, respectively.\ The DRO family originates from small circular retrograde Keplerian approximation centered at the Moon and constitutes a continuous family of planar periodic orbits that are symmetric with respect to the $x$-axis \cite{henon_1969,henon_1974,broucke}.\ Figure~\ref{figure_plot_1} shows the evolution of the DRO family together with the vertical stability index $s_v$ curve.\ The DRO family is continued until the orbits approach collision with the Earth, reaching a minimum normalized distance of $0.002252$.\ Their data are collected in Table~\ref{data_1}.\ As the Jacobi constant decreases, the family passes through three geometrical regimes according to the regions explored by the orbits:\
\begin{itemize}
	\item[(i)] In the first geometric regime, the DROs remain close to the Moon, illustrated by the first six green orbits in Figure~\ref{plot_1_c}.\ Along the vertical stability index curve in Figure~\ref{plot_1_a}, the vertical stability index~$s_v$ decreases until it reaches its minimum.\ We observe that the minimal multiplicity where VSR bifurcations occur is $q=5$.\ We consider the first six VSR orbits in Figure~\ref{plot_1_a} with self-resonance ratios of the form $1$:$q$ or $(q-1)$:$q$, where $q=5,...,10$.\ Since the elliptic Floquet multipliers all have negative Krein signatures (see Table~\ref{data_1}), the self-resonance ratios are identified as $(q-1)$:$q$.\
	\item[(ii)] Beyond the minimum of the vertical stability index curve, we observe from Figure~\ref{plot_1_c} that the DROs expand to move around the Moon together with the collinear libration points $L_1$ and $L_2$, forming the second geometric regime (see the blue DROs in Figure~\ref{plot_1_c}).\ Along this part of the family, the same sequence of multiplicities $q=5,…,10$ is encountered with the same type of self-resonance ratios $(q-1)$:$q$.\
	\item[(iii)] Figure~\ref{plot_1_b} shows the continuation of the $s_v$ curve beyond the range from Figure~\ref{plot_1_a} that leads to the third geometric regime, formed by the two largest red DROs in Figure~\ref{plot_1_c}.\ The orbits enclose the Moon together with the collinear and triangular libration points $L_1$, $L_2$, $L_4$ and $L_5$.\ Each of these two DROs corresponds to vertical pitchfork bifurcations.\
\end{itemize}

\begin{figure}[t!]
	\centering
	
	\begin{subfigure}{0.95\textwidth}
		\centering
		\includegraphics[width=\textwidth]{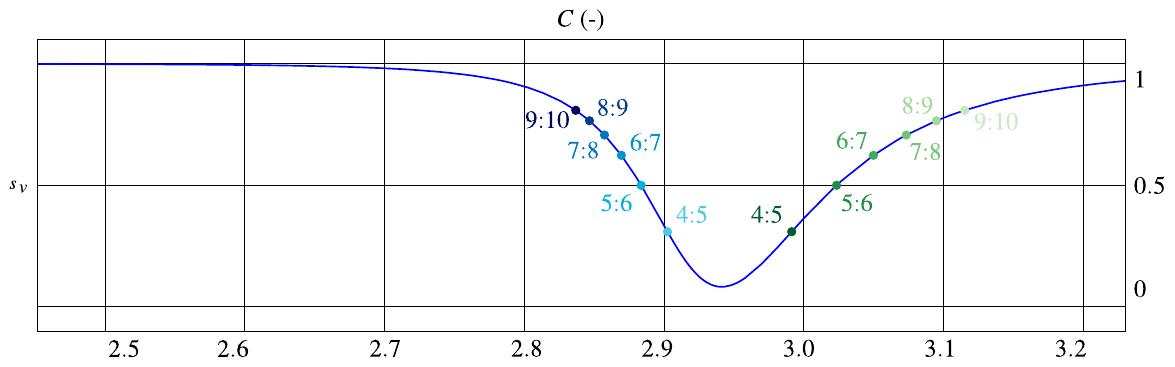}
		\caption{Vertical stability index $s_v$ as a function of the Jacobi constant of the DRO family.}
		\label{plot_1_a}
	\end{subfigure}
	\vspace{0.5cm}
	\begin{subfigure}{0.48\textwidth}
		\centering
		\includegraphics[width=\textwidth]{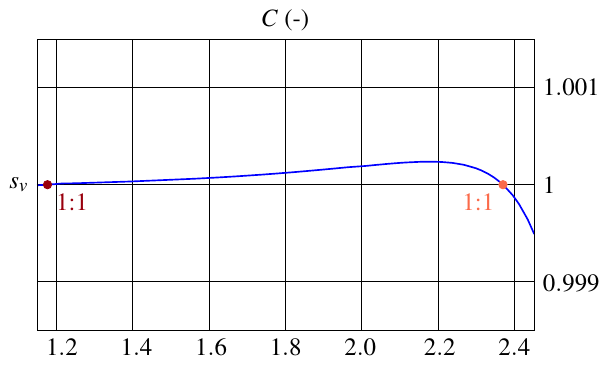}
		\caption{Continuation of the vertical stability index $s_v$ beyond the range shown in (a), zoomed to reveal the two vertical pitchfork bifurcations labeled by 1:1.}
		\label{plot_1_b}
	\end{subfigure}
	\hfill
	\begin{subfigure}{0.48\textwidth}
		\centering
		\includegraphics[width=\textwidth]{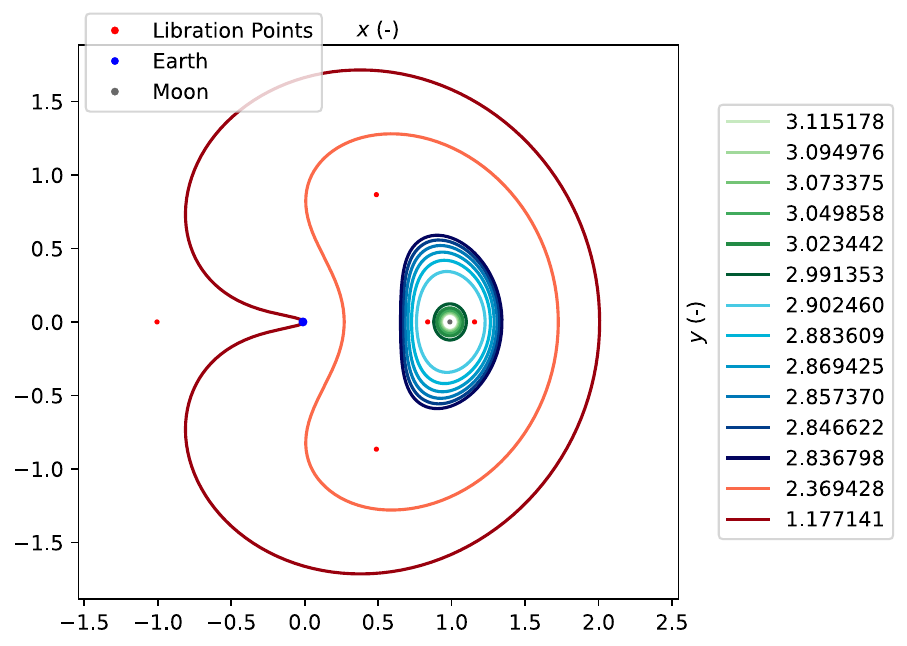}
		\caption{Representative DROs.}
		\label{plot_1_c}
	\end{subfigure}
	\caption{Vertical self-resonant DROs.\ The marked points in (a) and (b) indicate the detected vertical self-resonances and are labeled by their self-resonance ratios $p$:$q$.\ The colors used for the orbit plots in (c) coincide with those marking the corresponding vertical self-resonant points on the vertical stability index curve.}
	\label{figure_plot_1}
\end{figure}

The determination of the CZ indices, that are listed in Table~\ref{data_1}, are based on CZ indices that are known from analytical considerations when the planar orbits are generated from the Kepler problem \cite{aydin_babylonian}:
\begin{align}\label{indices_basic}
	\mu_{CZ} = \begin{cases}
		6 = \mu_{CZ}^p + \mu_{CZ}^s = 3 + 3 & \text{for the prograde periodic orbit family},\\
		2 = \mu_{CZ}^p + \mu_{CZ}^s = 1 + 1 & \text{for the retrograde periodic orbit family}.
	\end{cases}
\end{align}
In view of \eqref{indices} we see that both the planar and vertical winding numbers of the prograde orbits exceed those of the retrograde orbits by two.\ This observation plays an important role in the further analysis in this section.\ However, starting with the CZ indices \eqref{indices_basic}, monitoring the Floquet multipliers with Krein signatures as discussed in Section \ref{sec3}, yields the CZ indices of the DRO family in Table~\ref{data_1}.\

\subsection{P3DRO Family}

\begin{figure}[t!]
	\centering
	\begin{subfigure}{0.24\textwidth}
		\centering
		\includegraphics[scale=0.45]{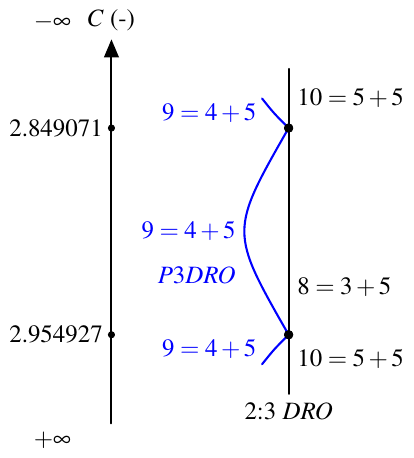}
		\caption{}
		\label{plot_1_1_a}
	\end{subfigure}
	\hfill
	\begin{subfigure}{0.24\textwidth}
		\centering
		\includegraphics[width=\textwidth]{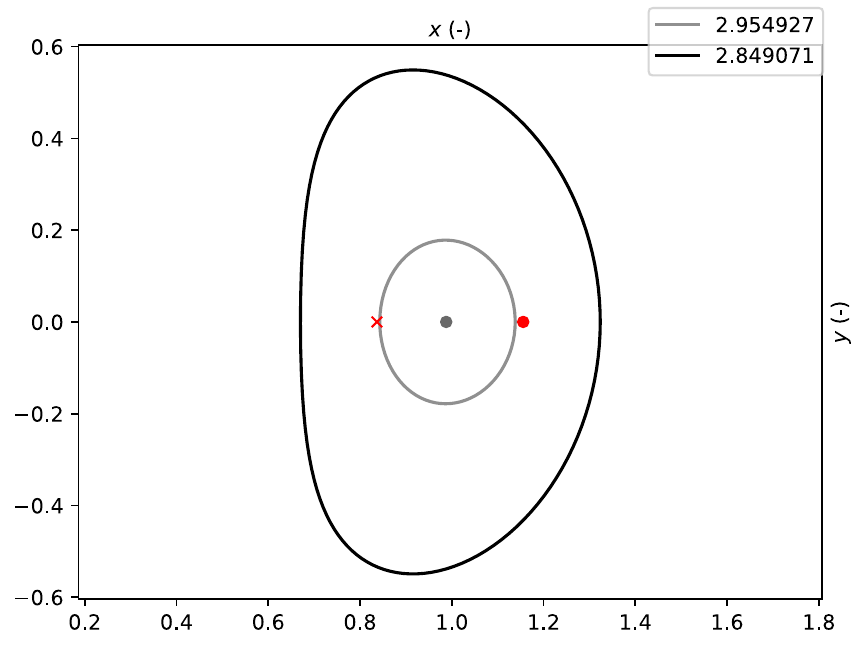}
		\caption{}
		\label{plot_1_1_b}
	\end{subfigure}
	\hfill
	\begin{subfigure}{0.24\textwidth}
		\centering
		\includegraphics[width=\textwidth]{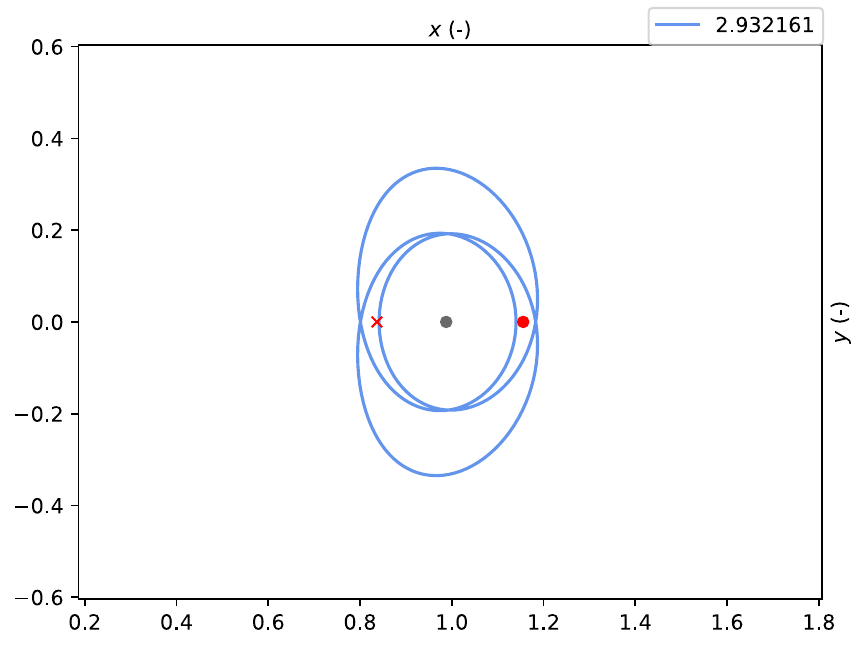}
		\caption{}
	\end{subfigure}
	\hfill
	\begin{subfigure}{0.24\textwidth}
		\centering
		\includegraphics[width=\textwidth]{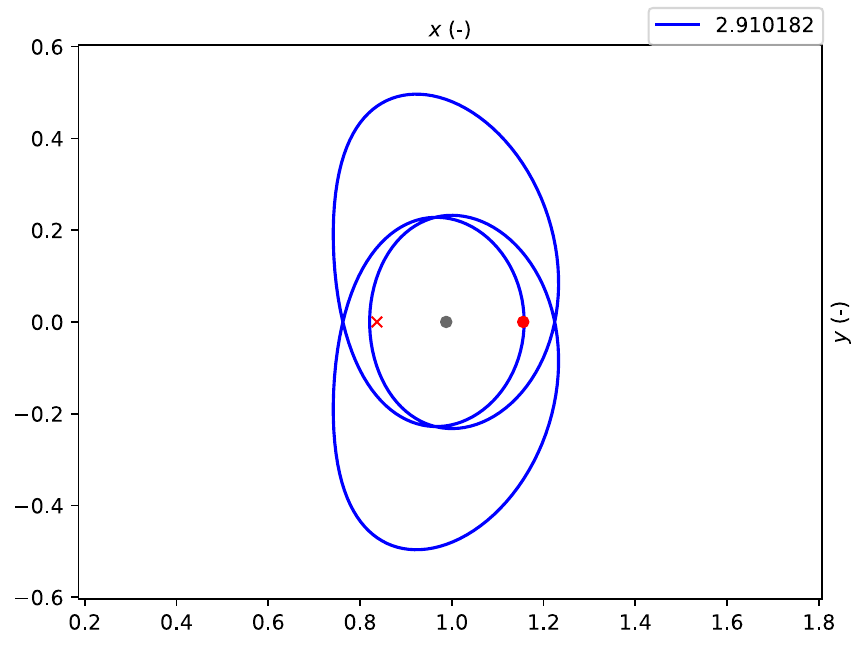}
		\caption{}
		\label{plot_1_1_d}
	\end{subfigure}
	\caption{Generation of the P3DRO family through two touch-and-go bifurcations.\ (a) Bifurcation diagram illustrating the two touch-and-go bifurcations associated with the planar period-tripling bifurcations of the DRO family.\ Both bifurcations belong to the same P3DRO family.\ (b) shows 2:3 DROs at the corresponding period-tripling bifurcation points.\ (c) and (d) show representative P3DROs continued from the lower to the upper touch-and-go bifurcations.\ Within the orbit plots:\ \textcolor{gray}{$\bullet$}, \textcolor{red}{$\times$} and \textcolor{red}{$\bullet$} denotes the Moon and the libration points $L_1$ and $L_2$, respectively.}
\end{figure}

The horizontal stability of the DRO family exhibits a behavior similar to that of the vertical stability index.\ As the Jacobi constant decreases, the horizontal stability index first decreases and then increases, while the corresponding planar Floquet multipliers remain elliptic with negative Krein signature throughout the considered range (see Table~\ref{data_1}).\ The lowest multiplicity where planar self-resonant bifurcations occur is $q = 3$, first at $C \approx 2.954927$ and second at $C \approx 2.849071$.\ These two period-tripling bifurcations correspond to touch-and-go bifurcations, that is illustrated in Figure~\ref{plot_1_1_a} with corresponding orbit plots in Figures~\ref{plot_1_1_b}--\ref{plot_1_1_d}.\ Such three-revolution orbits are denoted as P3DRO (period-3 DRO).\ At the touch-and-go bifurcations, two families bifurcate into opposite directions of the Jacobi constant, and along both bifurcated families, the critical planar Floquet multipliers become positive hyperbolic.\ Following the family toward decreasing Jacobi constant from the first touch-and-go bifurcation, the P3DROs reconnect with the DRO family at the second touch-and-go bifurcation, as shown in Figure~\ref{plot_1_1_a}.\ Therefore, the bifurcated families can be described by a single family that consists of P3DROs.\ This phenomenon was first described by Hénon \cite{henon_1970,henon_2003} in the Hill problem, and was also discussed in the Earth--Moon CR3BP \cite{capdevilla}.\ We also notice that the computed critical Jacobi constants agree closely with the values reported in the latter reference.\ In this work, we assign the CZ indices to the P3DRO family whose data are collected in Table~\ref{data_1}.\ At the first touch-and-go bifurcation, the planar CZ index of the period-3 DRO jumps from 5 to~3, while at the second touch-and-go bifurcation, the planar CZ index jumps from 3 back to 5.\ The spatial CZ index remains constant 5 in both cases.\ Due to the touch-and-go characteristics, the P3DROs close to the bifurcation points have planar CZ index 4 with positive hyperbolic Floquet multipliers, while the CZ index is 5.\ In the range between the two touch-and-go bifurcations, these CZ indices remain constant while the spatial Floquet multipliers have elliptic and negative hyperbolic phases indicating the occurrence of vertical period-doubling bifurcations (see Table~\ref{data_1}).\ We also notice that these CZ indices are consistent with those computed for the entire family within the Hill problem (denoted by $f_3$) \cite{aydin_batkhin}.\

\subsection{$L_2$ Planar Lyapunov, Halo and Butterfly Orbit Family}

A second class of basic periodic orbit families consists of the planar and vertical Lyapunov families that are generated around the collinear libration points $L_1$, $L_2$ and $L_3$.\ These families originate from the planar and vertical frequencies $\omega_p$ and $\omega_v$, respectively, associated with the linear dynamics of type $\text{saddle} \times \text{center} \times \text{center}$.\ From the planar frequency $\omega_p$, one obtains small retrograde elliptic periodic solutions centered at each of the collinear libration points, that are symmetric with respect to the $x$-axis (see \cite{broucke}, \cite[Chapter 5]{szebehely}).\ It was shown in \cite{aydin_batkhin} that, in the vicinity of the collinear libration points, the frequencies satisfy $\omega_v < \omega_p < 2 \omega_v$, and that the CZ indices of the planar Lyapunov orbits are given by
\begin{align}\label{indices_basic_2}
		\mu_{CZ} = \mu_{CZ}^p + \mu_{CZ}^s = 2 + 1.
\end{align}

\begin{figure}[t!]
	\centering
	
	\begin{subfigure}{0.32\textwidth}
		\centering
		\includegraphics[width=\textwidth]{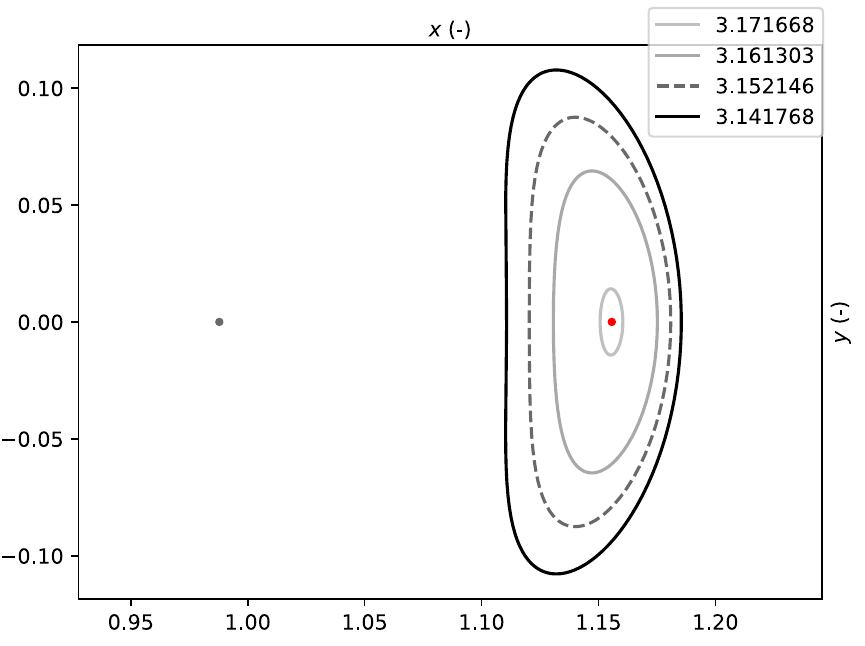}
		\caption{}
		\label{plot_2_a}
	\end{subfigure}
	\hfill
	\begin{subfigure}{0.32\textwidth}
		\centering
		\includegraphics[width=\textwidth]{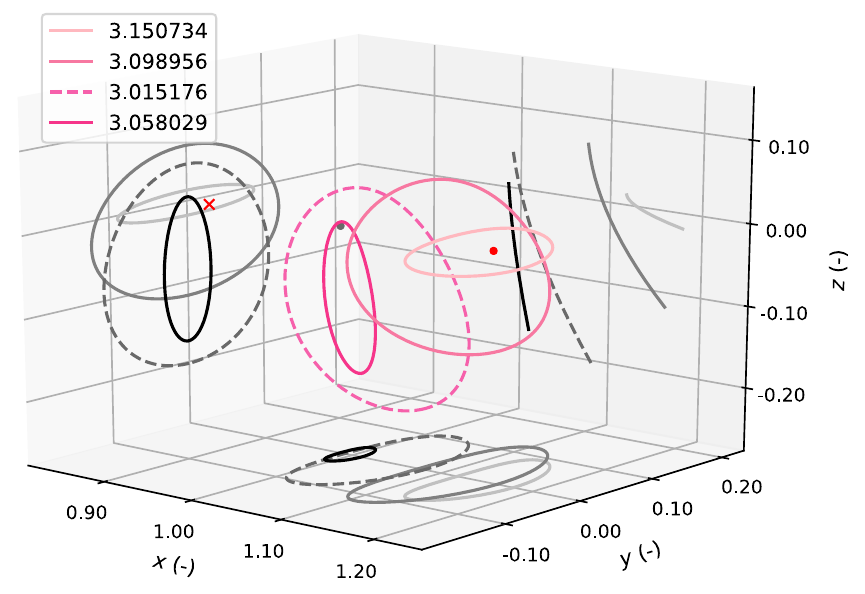}
		\caption{}
		\label{plot_2_b}
	\end{subfigure}
	\hfill
	\begin{subfigure}{0.32\textwidth}
		\centering
		\includegraphics[width=\textwidth]{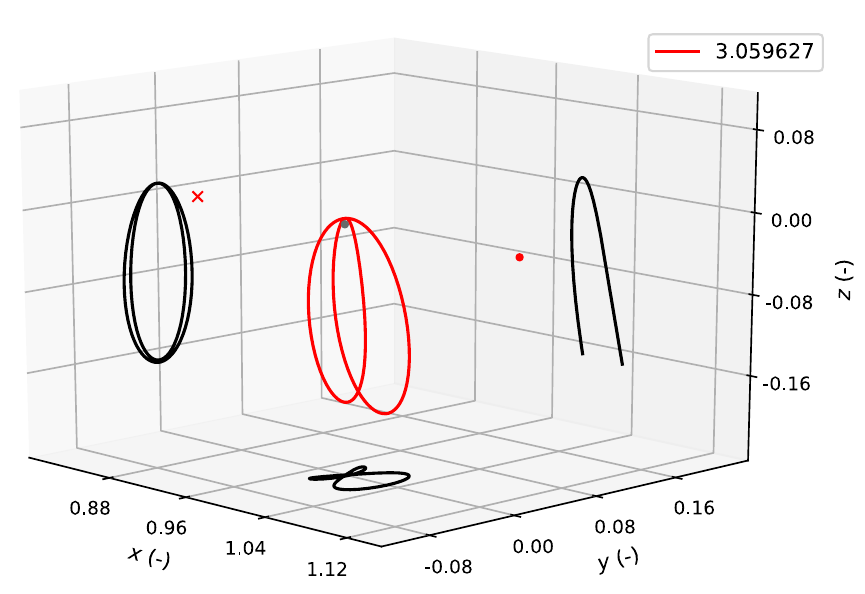}
		\caption{}
		\label{plot_2_c}
	\end{subfigure}
	\vspace{0.5cm}
	
	\begin{subfigure}{0.32\textwidth}
		\centering
		\includegraphics[width=\textwidth]{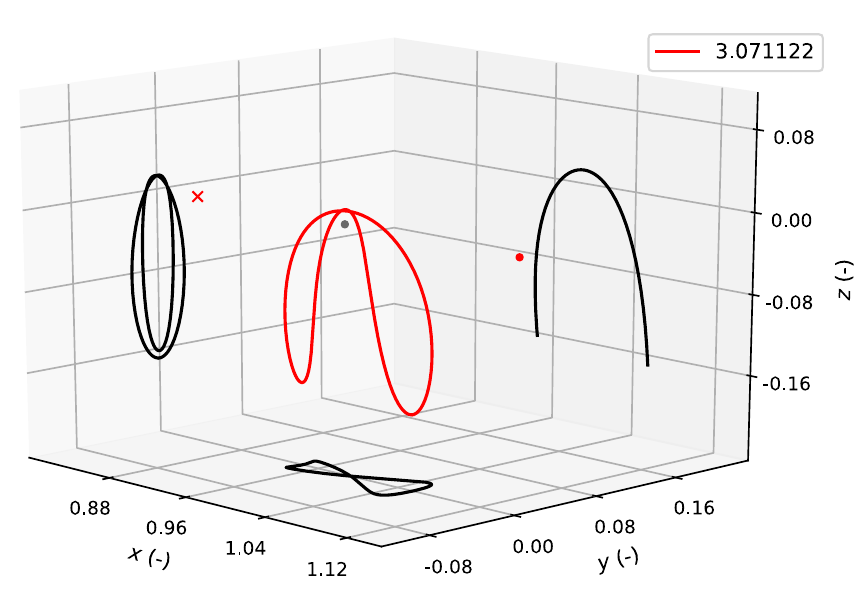}
		\caption{}
	\end{subfigure}
	\hfill
	\begin{subfigure}{0.32\textwidth}
		\centering
		\includegraphics[width=\textwidth]{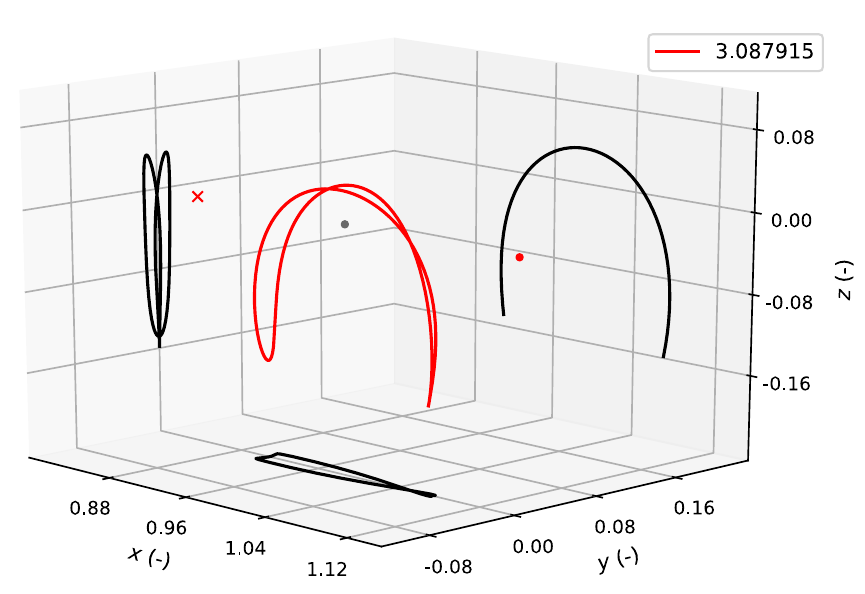}
		\caption{}
	\end{subfigure}
	\hfill
	\begin{subfigure}{0.32\textwidth}
		\centering
		\includegraphics[width=\textwidth]{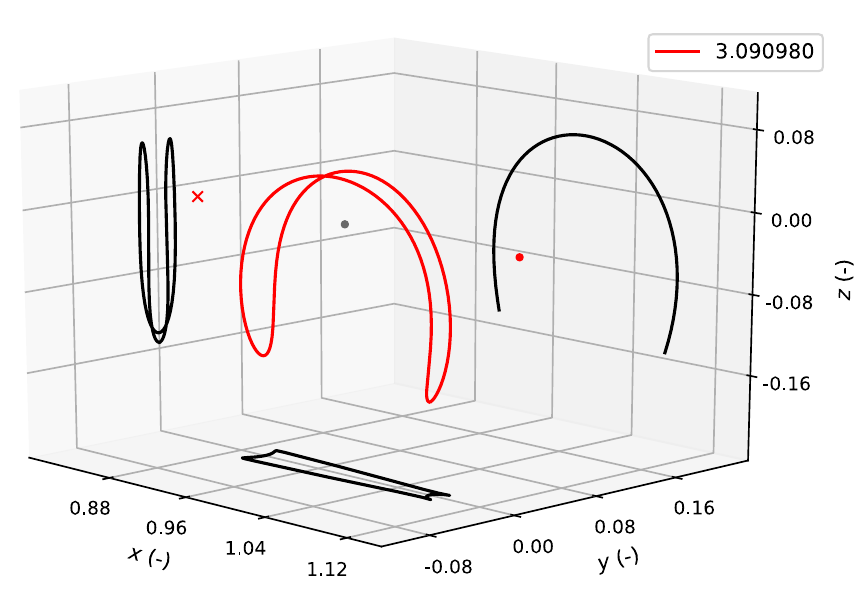}
		\caption{}
		\label{plot_2_f}
	\end{subfigure}
	\caption{(a) Planar $L_2$ Lyapunov orbits; from the gray dashed the $L_2$ halo family bifurcates.\ (b) Southern $L_2$ halo orbits; pink dashed corresponds to birth--death bifurcation point.\ (c)--(f) Evolution of southern Butterfly orbits.\ Within the plots:\ \textcolor{gray}{$\bullet$}, \textcolor{red}{$\times$} and \textcolor{red}{$\bullet$} denotes the Moon and the libration points $L_1$ and $L_2$, respectively.}
\end{figure}

In this paper, we consider the $L_2$ planar Lyapunov orbits, some of which are shown Figure~\ref{plot_2_a}.\ Their data can be found in Table~\ref{data_1}.\ Together with the CZ indices \eqref{indices_basic_2}, the $L_2$ planar Lyapunov orbits start being planar positive hyperbolic and spatial elliptic with negative Krein signature.\ After short continuation, at $C \approx 3.152146$ the orbits become spatial positive hyperbolic, where the spatial CZ index jumps from 1 to 2 due to the negative Krein signature.\ This vertical pitchfork bifurcation gives rise to the well-known $L_2$ halo family \cite{farquhar,breakwell,howell,gomez}, whose orbits are symmetric with respect to $\rho_{xz}$.\ Owing to the $\sigma$-symmetry, the $L_2$ halo family consists of two symmetry-related families, commonly referred to as the \textit{northern} and \textit{southern halo families}, depending on whether the orbits extend above or below the ecliptic plane $\{z = 0\}$.\ Some southern $L_2$ halo orbits are shown in Figure~\ref{plot_2_b} and their data are collected in Table~\ref{data_2}.\ The $L_2$ halo orbits start with CZ index 3.\ Along the continuation, the Jacobi constant reaches a minimum value of $C \approx 3.015176$ where a birth--death bifurcation occurs.\ Beyond this bifurcation point, the CZ index changes to 4 and the orbits become near-rectilinear, commonly known as \textit{Near-Rectilinear Halo Orbits} (NRHOs).\ In the NRHO regime, a period-doubling bifurcation occurs at $C \approx 3.058029$, generating the well-known Butterfly family \cite{howell_campbell,grebow}, whose orbits are symmetric with respect to $\rho_{xz}$.\ Since the bifurcation takes place on the northern and southern NRHO family, the Butterfly family likewise consists of two $\sigma$-symmetry-related branches, referred to as the \textit{northern} and \textit{southern Butterfly families}.\ Figures~\ref{plot_2_c}--\ref{plot_2_f} show some southern Butterfly orbits, and their data can be seen in Table~\ref{data_2}.\

\subsection{VSR Bifurcations from near-Moon DROs}

We now investigate the six VSR bifurcations of the near-Moon DRO family in the first geometric regime, corresponding to the self-resonances $(q-1)$:$q$ with $q=5,6,7,8,9,10$.\ Recall from Section~\ref{sec3} that each VSR bifurcation generates two spatial symmetric periodic orbits families whose symmetry depends on the parity of $q$:\ odd-period VSR bifurcations produce $\rho_{xz}$- and $\rho_x$-symmetric families (together with their $\sigma$-symmetric counterparts), while even-period VSR bifurcations produce doubly symmetric families.\ Throughout this work, in view of Figure~\ref{fig:bifurcations}, these bifurcated families are referred to as the \textit{red} and \textit{blue family} according to the stability of their critical Floquet multipliers.\ The red families possess elliptic critical Floquet multipliers, while the blue families possess positive hyperbolic critical Floquet multipliers.\ All computed bifurcated families are continued in increasing direction of the Jacobi constant until they naturally terminate at another known periodic orbit family, forming bridge families in the sense of \cite{llibre}.\

Besides the basic DRO family, the CR3BP possesses several planar prograde periodic orbit families around the smaller primary.\ Following the notation from \cite{restrepo_russell}, these are denoted by $LPO_1$, $DPO$ and $LPO_2$, where $LPO_1$ (Low Prograde Orbit) represents the family of small prograde periodic orbits around the smaller primary (that is the family with CZ indices \eqref{indices_basic}), while $LPO_2$ forms a second planar prograde family connected to the $DPO$ (Distant Prograde Orbit) family through a birth--death (fold) bifurcation.\ The CZ indices of these prograde families for small mass ratios were studied in \cite{moreno_aydin}.\ As will be shown below, these prograde families provide the terminal families of almost all VSR bridge families.\

\begin{figure}[t!]
	\centering
	\includegraphics[scale=0.8]{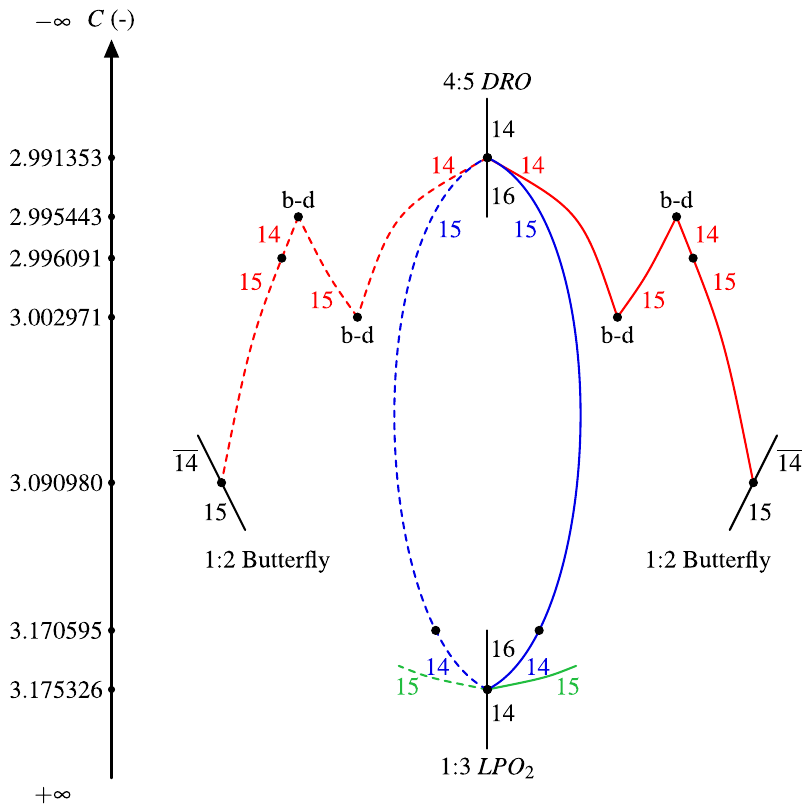}
	\caption{Bifurcation graph showing the bridge families between 4:5 DRO, 1:2 Butterfly and 1:3 $LPO_2$.}
	\label{figure_bif_1}
\end{figure}

\begin{figure}[t!]
	\centering
	
	\begin{subfigure}{\textwidth}
		\centering
		\includegraphics[width=0.32\textwidth]{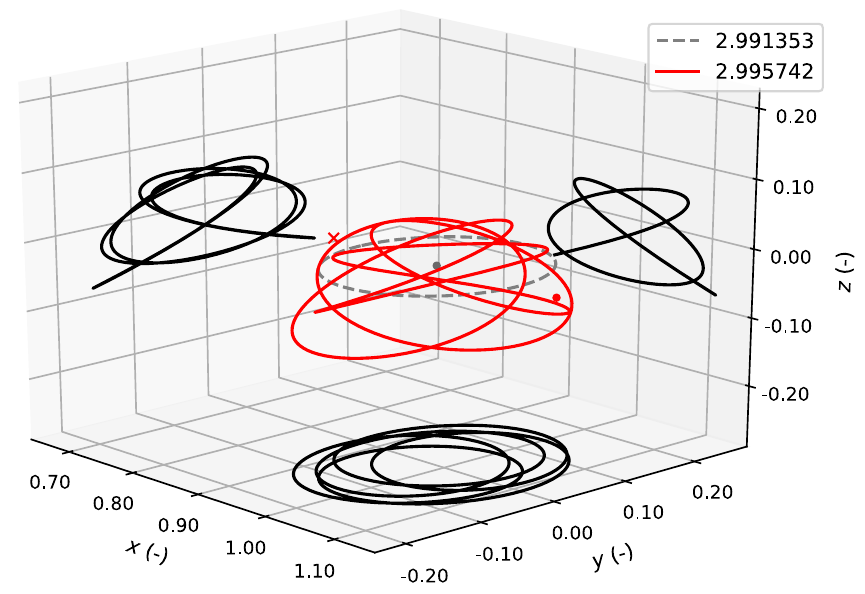}
		\hfill
		\includegraphics[width=0.32\textwidth]{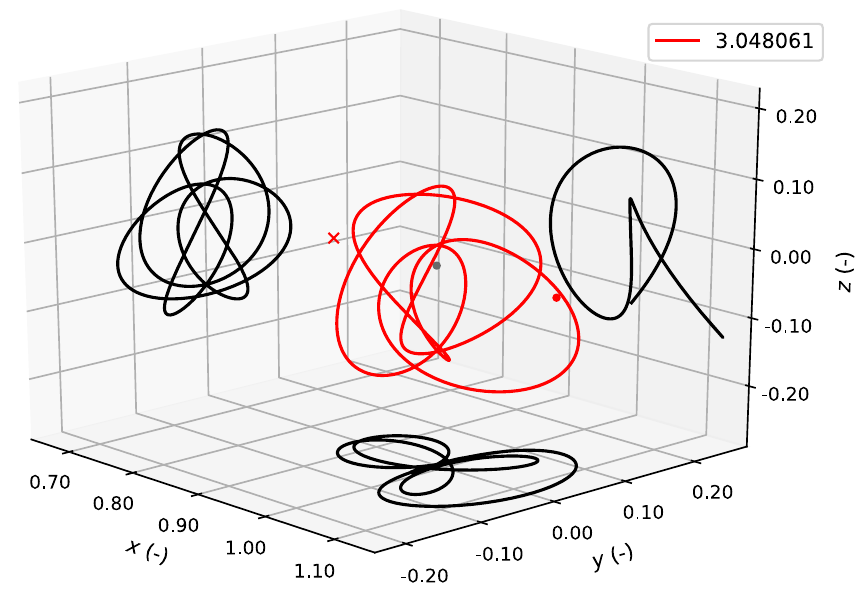}
		\hfill
		\includegraphics[width=0.32\textwidth]{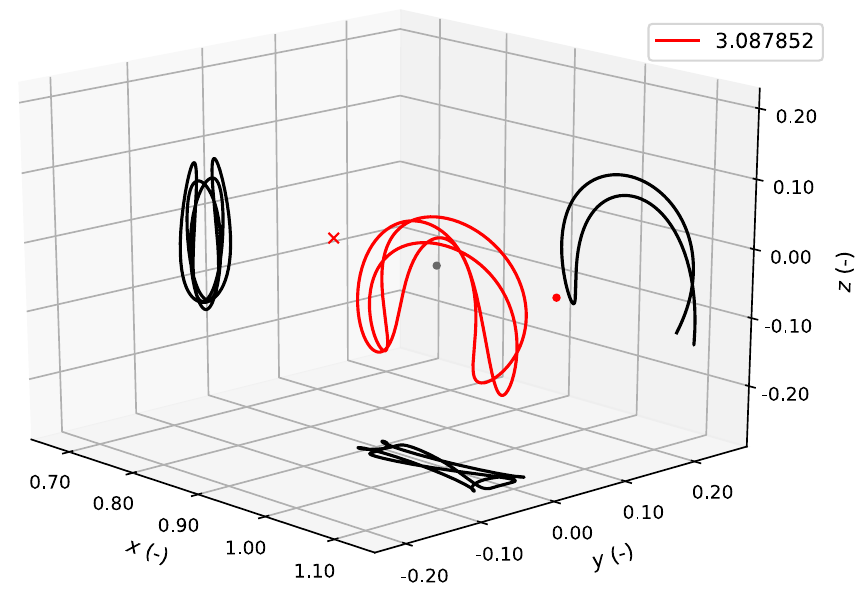}
		\caption{Representative $\rho_{xz}$-symmetric orbits of the red family, continued from left to right.\ They bifurcate from the 4:5 DRO (gray dashed) and terminate at the period-doubling bifurcation of the Butterfly family.}
	\end{subfigure}
	\vspace{0.5em}
	
	\begin{subfigure}{\textwidth}
		\centering
		\includegraphics[width=0.32\textwidth]{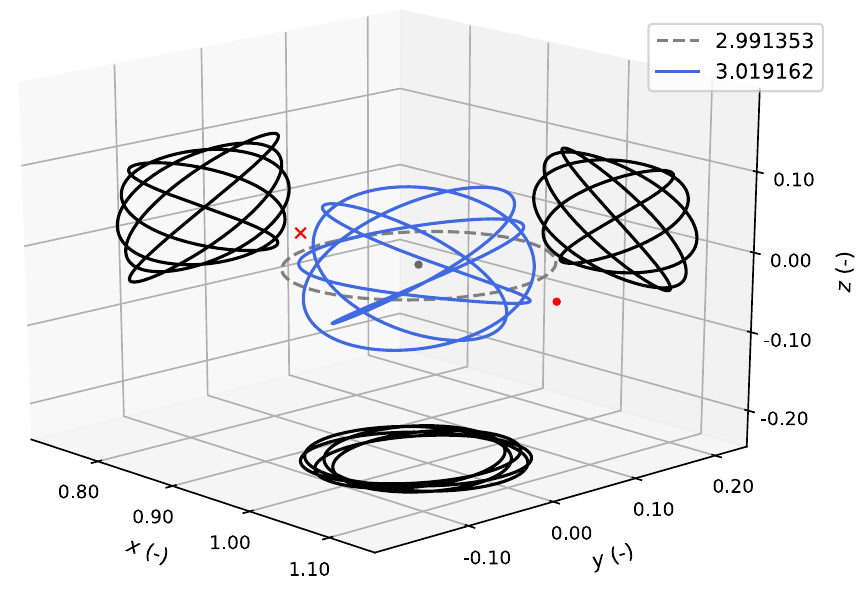}
		\hfill
		\includegraphics[width=0.32\textwidth]{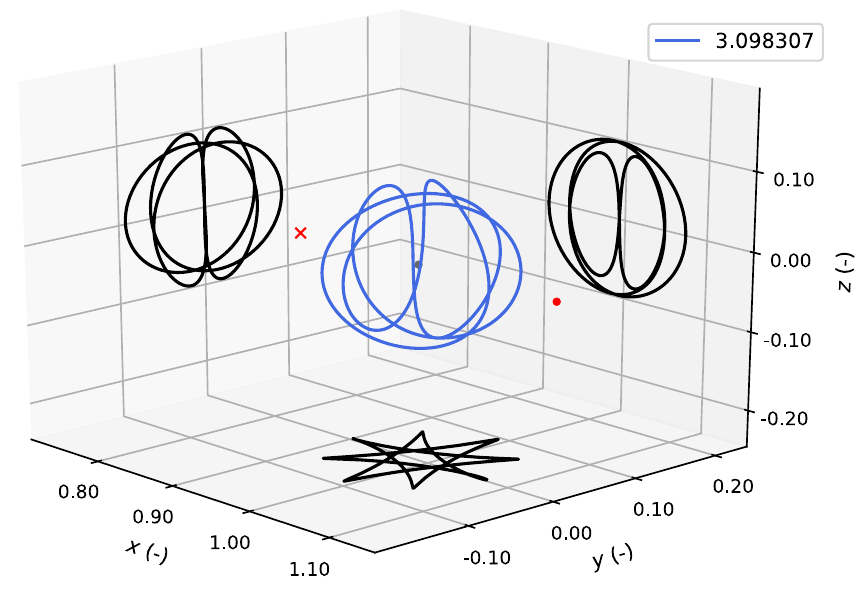}
		\hfill
		\includegraphics[width=0.32\textwidth]{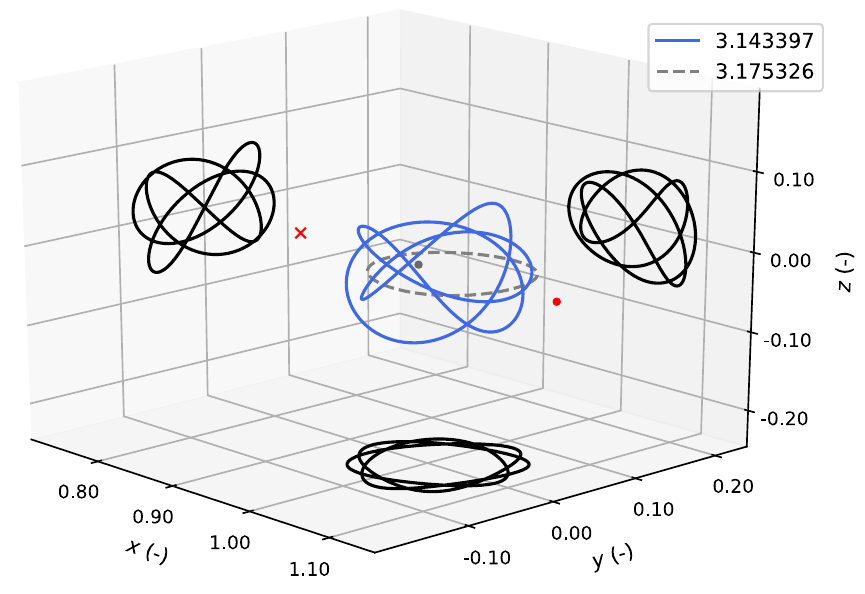}
		\caption{Representative $\rho_x$-symmetric orbits of the blue family, continued from left to right.\ They bifurcate from the 4:5 DRO (gray dashed) and terminate at the 1:3 $LPO_2$ (gray dashed).\ Middle shows a 1:1 Earth--Moon resonant orbit.}
	\end{subfigure}
	\vspace{0.5em}
	
	\begin{subfigure}{\textwidth}
		\centering
		\includegraphics[width=0.45\textwidth]{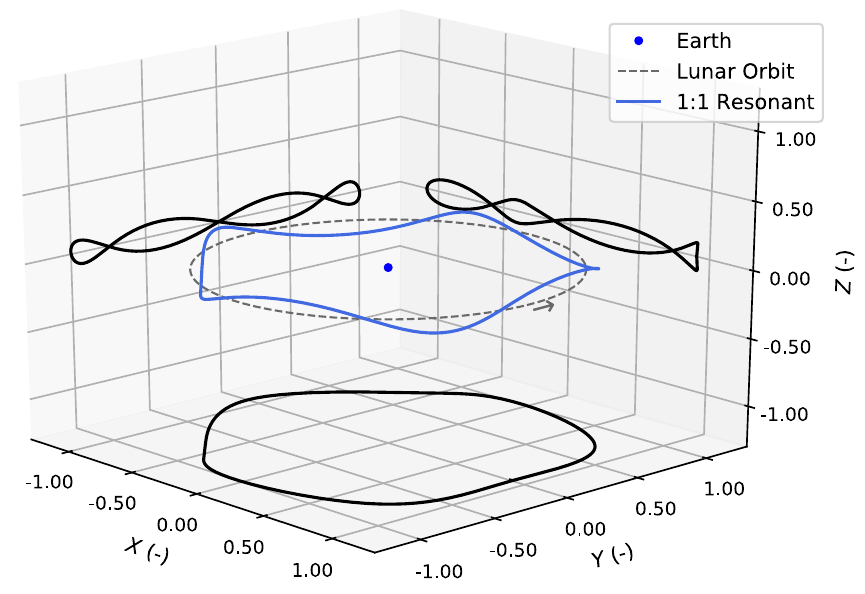}
		\includegraphics[width=0.45\textwidth]{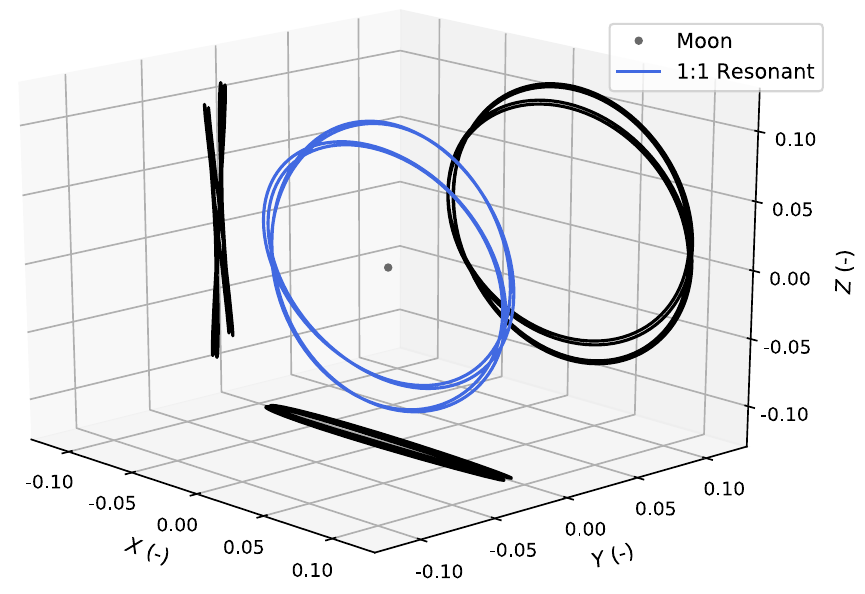}
		\caption{The 1:1 Earth--Moon resonant orbit from (b) is shown in the Earth-centered inertial frame (left) and in the Moon-centered inertial frame (right).}
		\label{plot_3_c}
	\end{subfigure}

	\caption{Vertical self-resonant bifurcations from the 4:5 DRO (near-Moon).\ Within the plots (a) and (b):\ \textcolor{gray}{$\bullet$}, \textcolor{red}{$\times$} and \textcolor{red}{$\bullet$} denotes the Moon and the libration points $L_1$ and $L_2$, respectively.}
	\label{figure_plot_3}
\end{figure}

\begin{figure}[t!]
	\centering
	\includegraphics[width=1\linewidth]{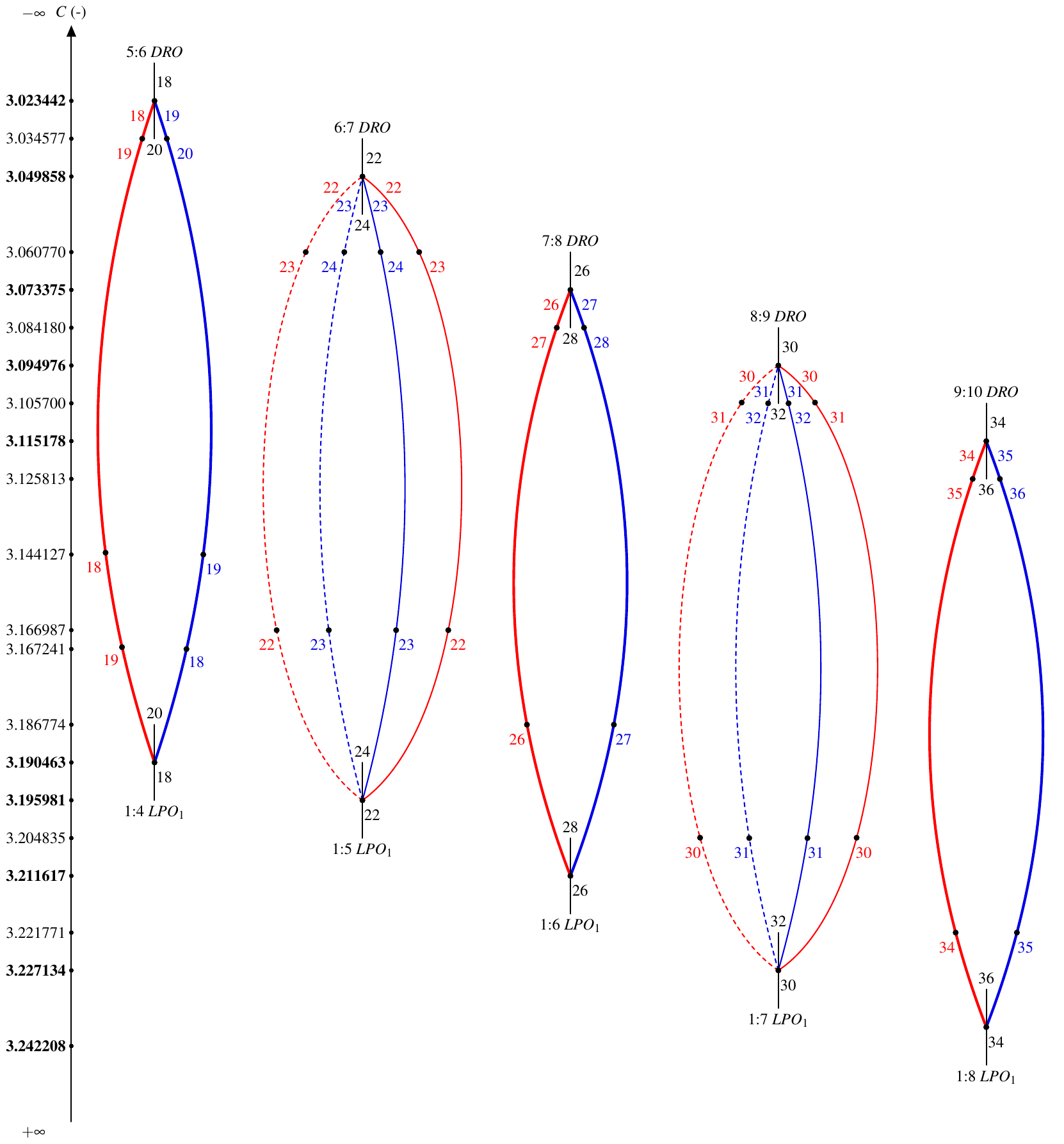}
	\caption{Bifurcation diagram that illustrates the bridge families between the $(q-1)$:$q$ DRO and the $1\text{:}(q-2)$ $LPO_1$, where $q=6,7,8,9,10$.\ Bold Jacobi constants on the left correspond to the VSR orbits of the DRO and $LPO_1$ families.}
	\label{figure_bif_2}
\end{figure}

\begin{figure}[t!]
	\centering
	
	\begin{subfigure}{\textwidth}
		\centering
		\includegraphics[width=0.32\textwidth]{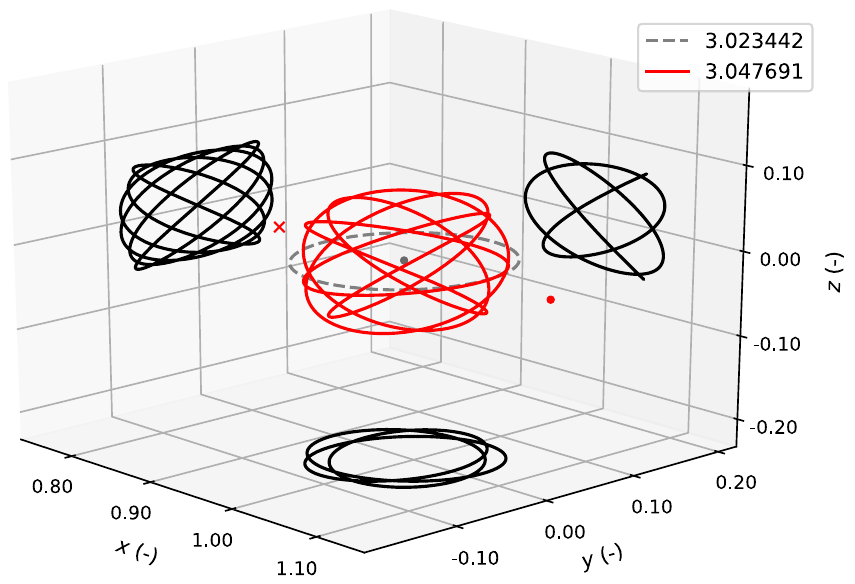}
		\hfill
		\includegraphics[width=0.32\textwidth]{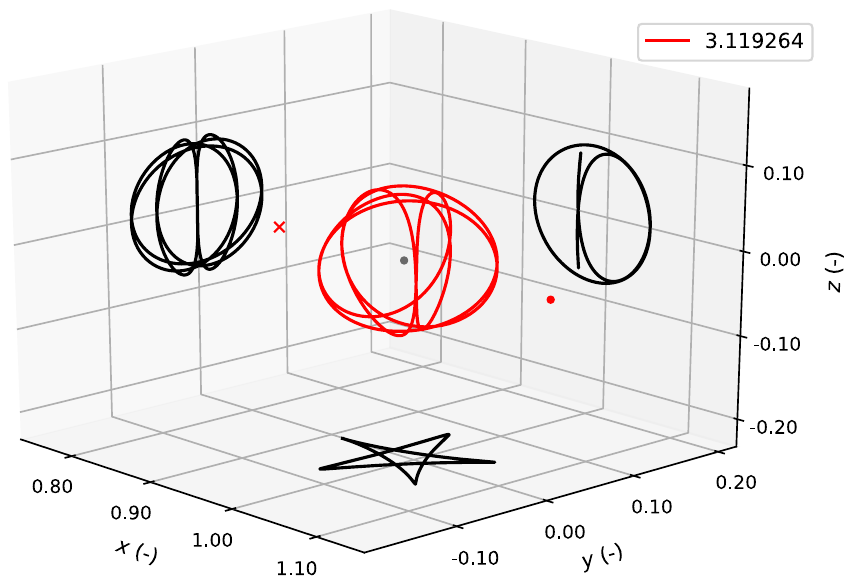}
		\hfill
		\includegraphics[width=0.32\textwidth]{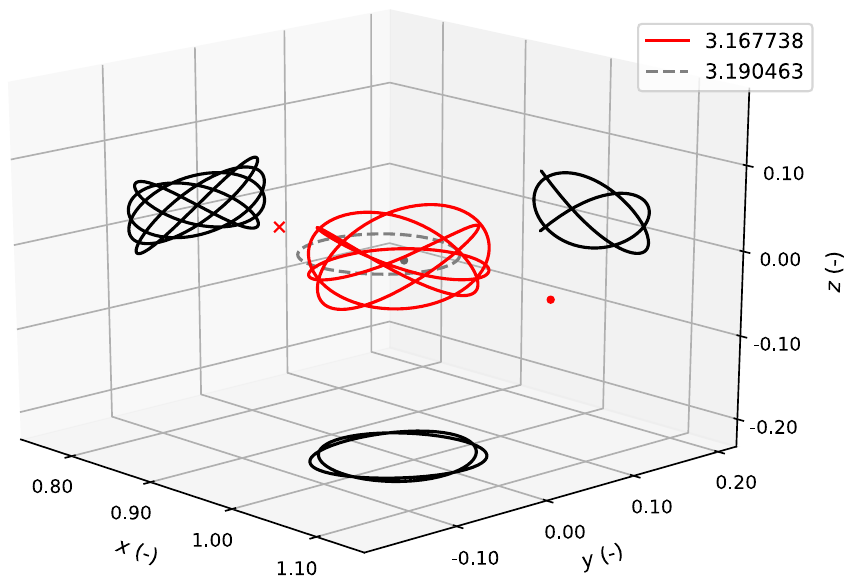}
		\caption{}
	\end{subfigure}
	\vspace{0.5em}
	
	\begin{subfigure}{\textwidth}
		\centering
		\includegraphics[width=0.32\textwidth]{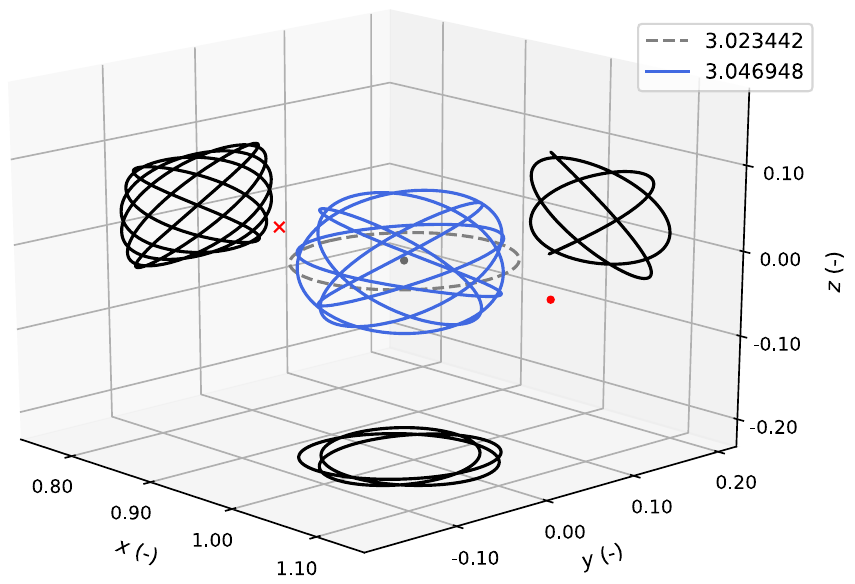}
		\hfill
		\includegraphics[width=0.32\textwidth]{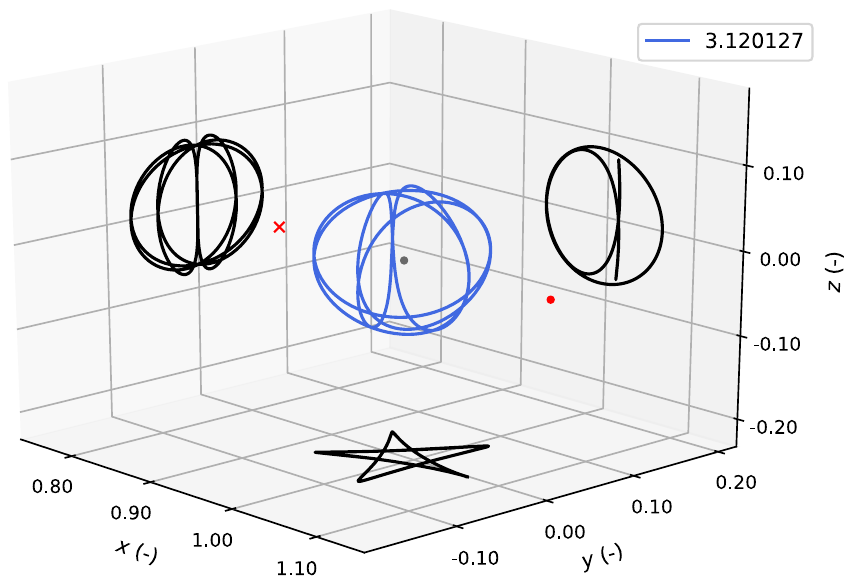}
		\hfill
		\includegraphics[width=0.32\textwidth]{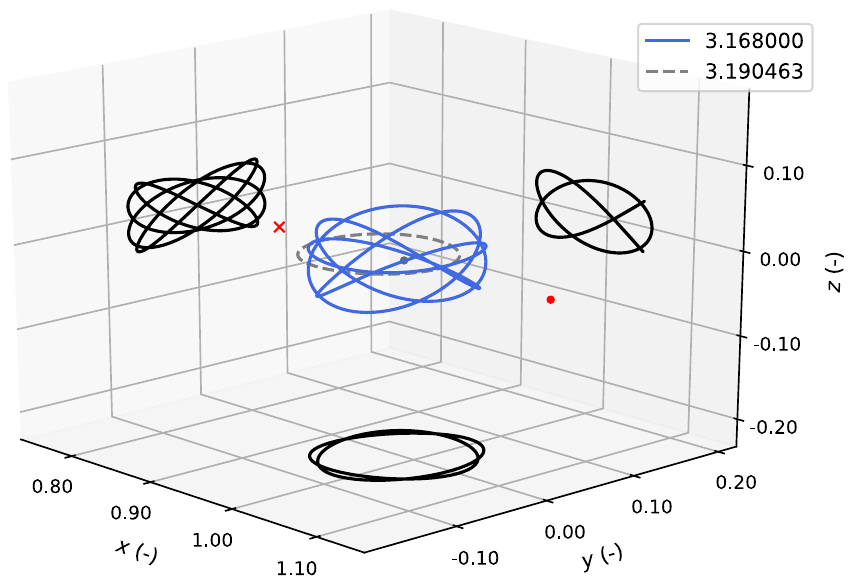}
		\caption{}
	\end{subfigure}
	\caption{Vertical self-resonant bifurcations from the 5:6 DRO (near-Moon).\ (a) and (b) show doubly symmetric orbits of the red and blue family, respectively, continued from left to right.\ They bifurcate from the 5:6 DRO (gray dashed) and terminate at the 1:4 $LPO_1$ (gray dashed).\ Within the plots:\ \textcolor{gray}{$\bullet$}, \textcolor{red}{$\times$} and \textcolor{red}{$\bullet$} denotes the Moon and the libration points $L_1$ and $L_2$, respectively.}
	\label{figure_plot_4}
\end{figure}

\begin{figure}[t!]
	\centering
	
	\begin{subfigure}{0.48\textwidth}
		\centering
		\begin{tabular}{cc}
			\includegraphics[width=.47\linewidth]{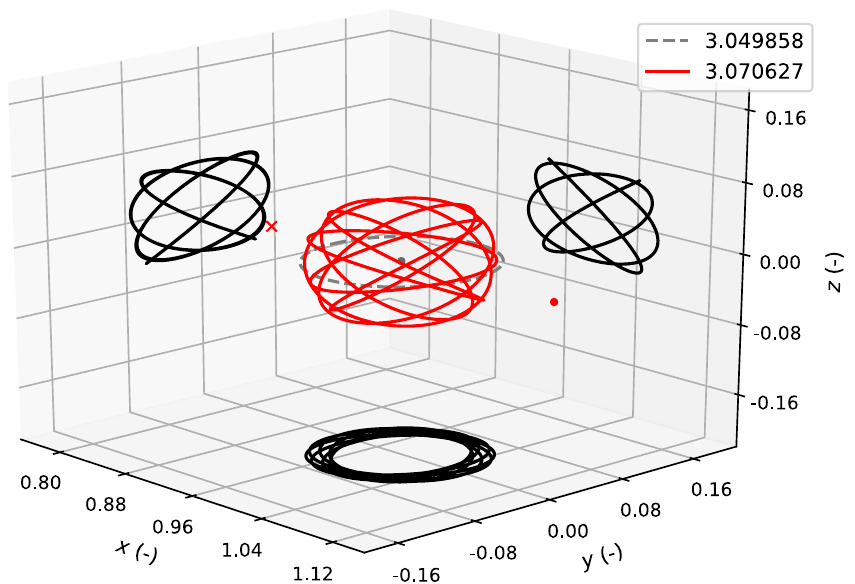} &
			\includegraphics[width=.47\linewidth]{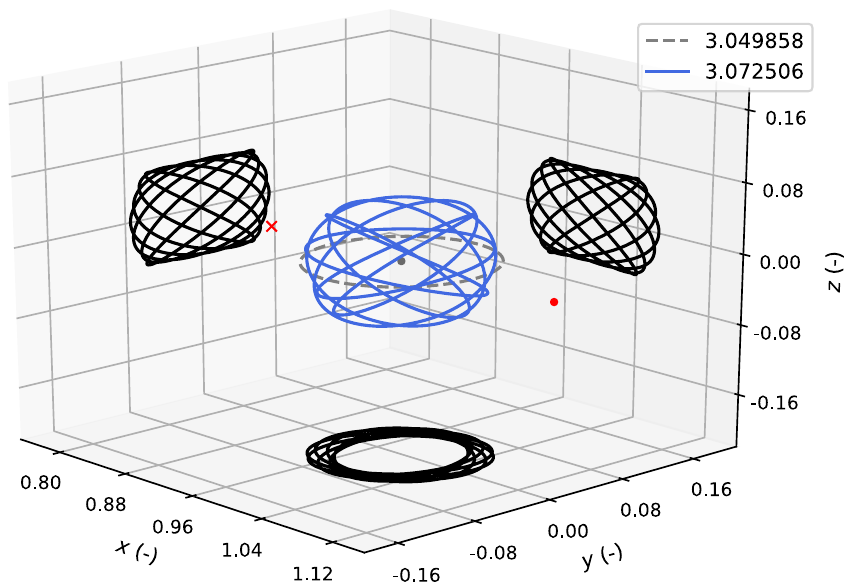} \\
			\includegraphics[width=.47\linewidth]{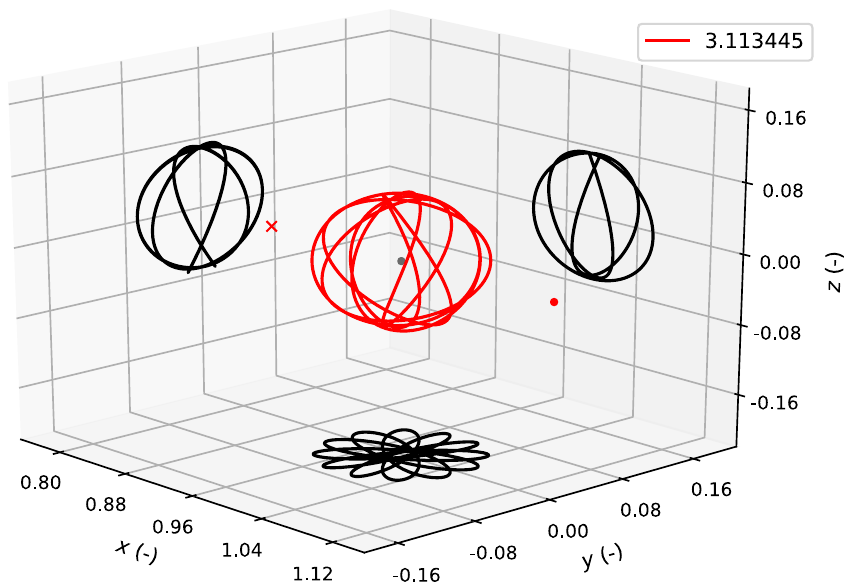} &
			\includegraphics[width=.47\linewidth]{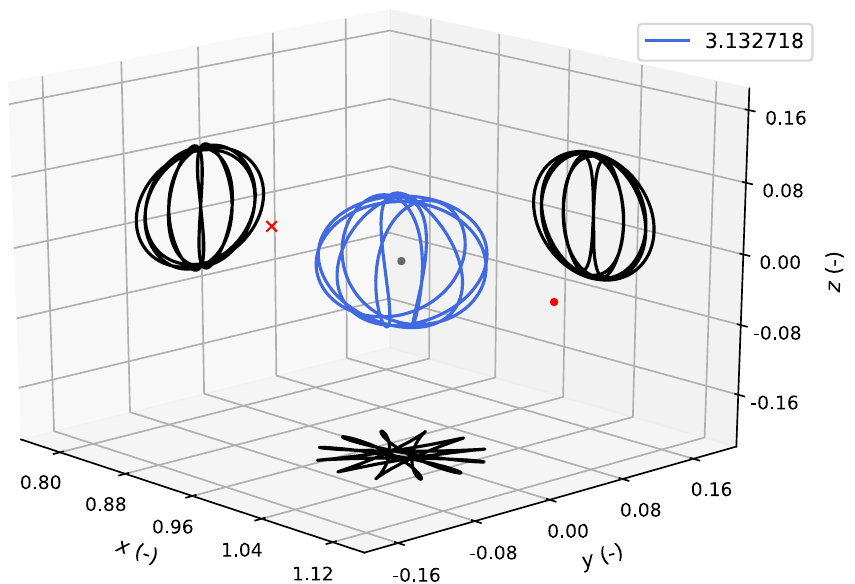} \\
			\includegraphics[width=.47\linewidth]{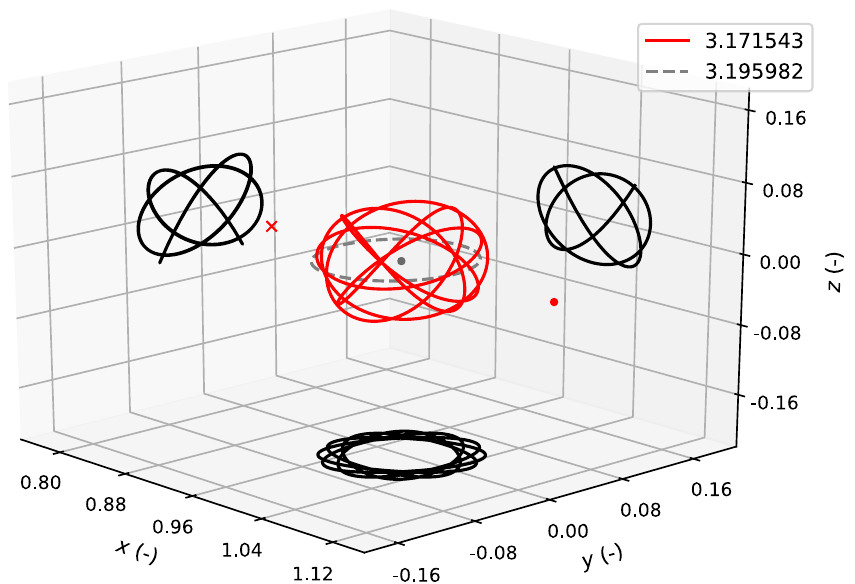} &
			\includegraphics[width=.47\linewidth]{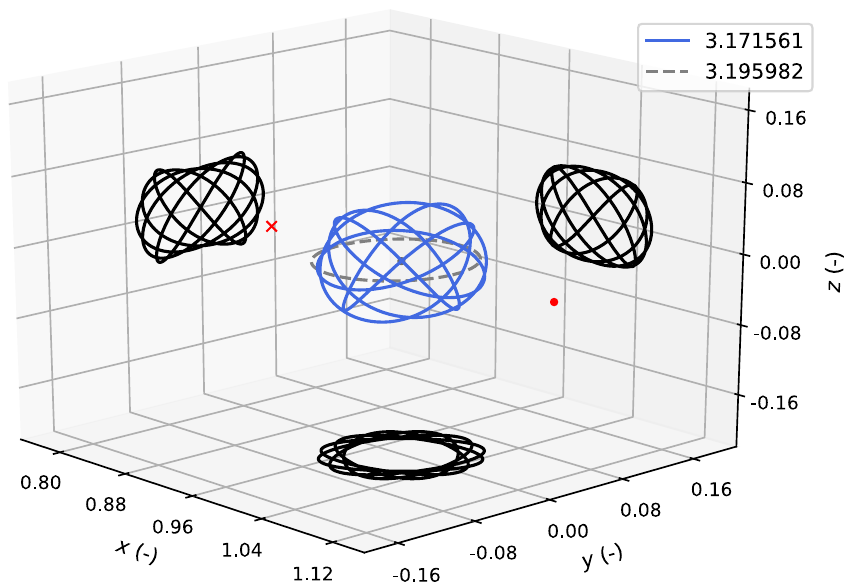}
		\end{tabular}
		\caption{6:7 DRO $\leftrightarrow$ 1:5 $LPO_1$ bridge families.}
	\end{subfigure}
	\hfill
	\begin{subfigure}{0.48\textwidth}
		\centering
		\begin{tabular}{cc}
			\includegraphics[width=.47\linewidth]{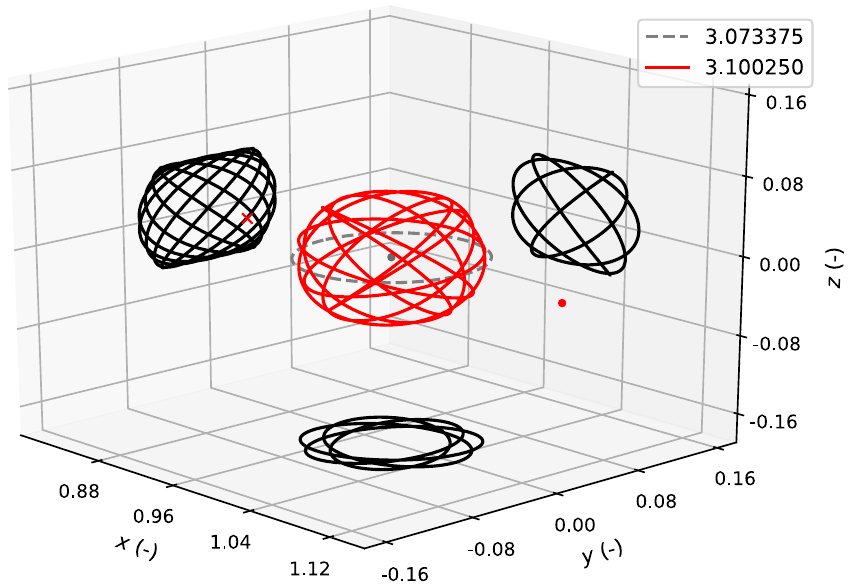} &
			\includegraphics[width=.47\linewidth]{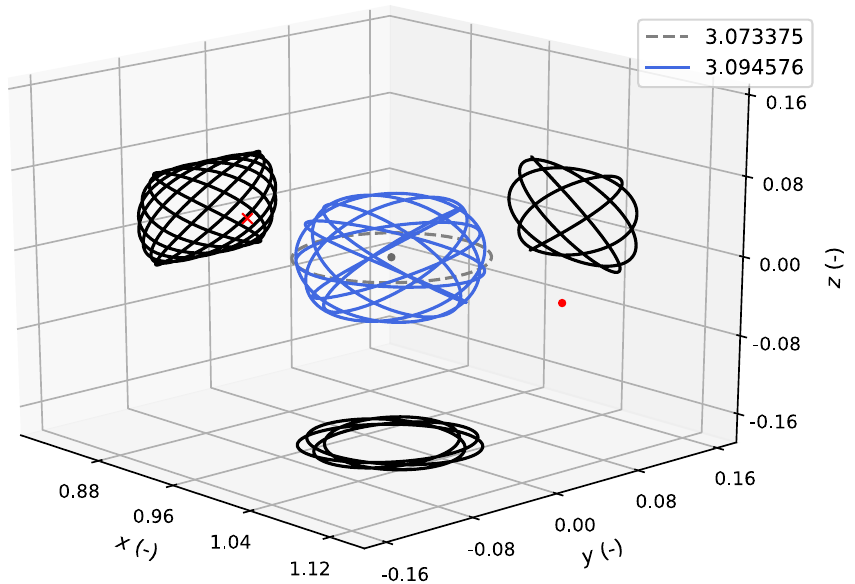} \\
			\includegraphics[width=.47\linewidth]{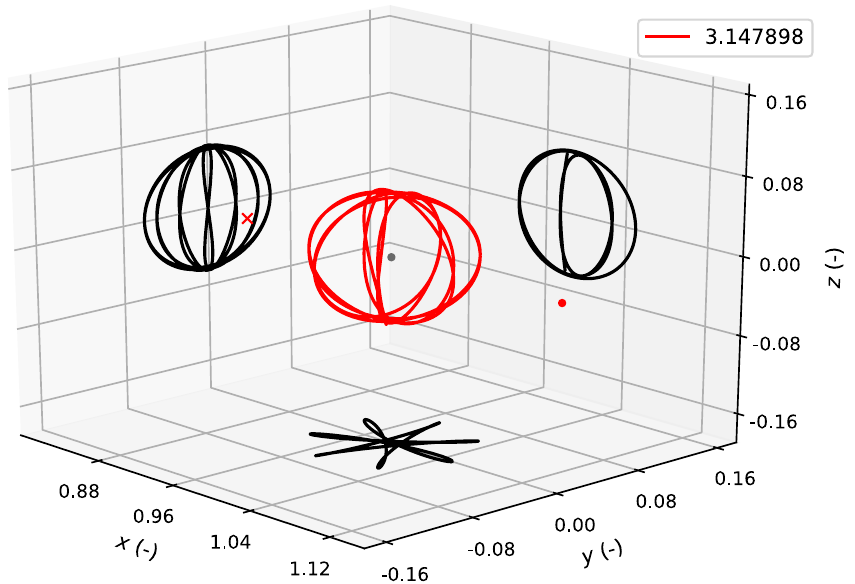} &
			\includegraphics[width=.47\linewidth]{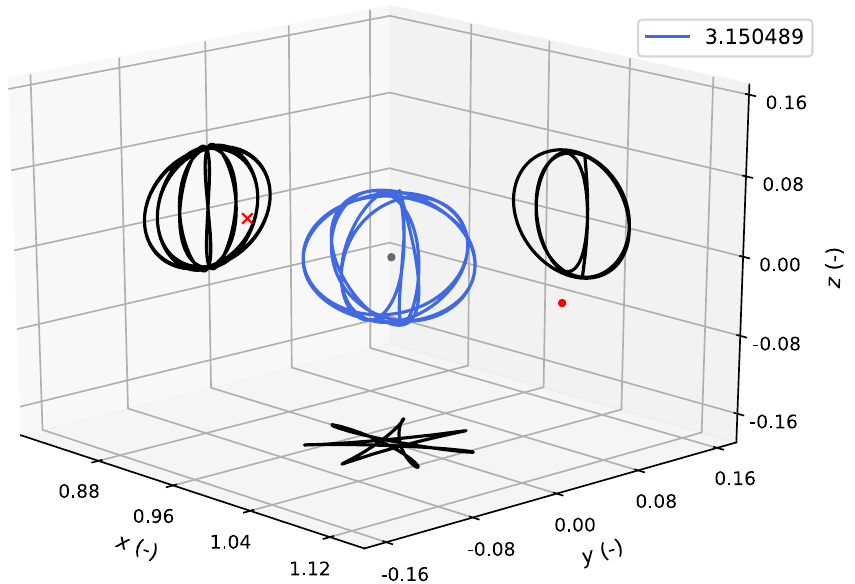} \\
			\includegraphics[width=.47\linewidth]{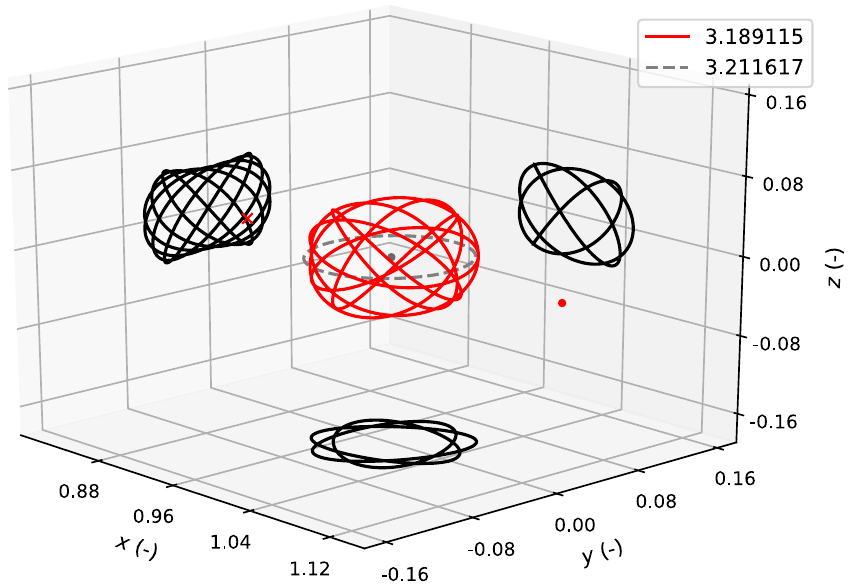} &
			\includegraphics[width=.47\linewidth]{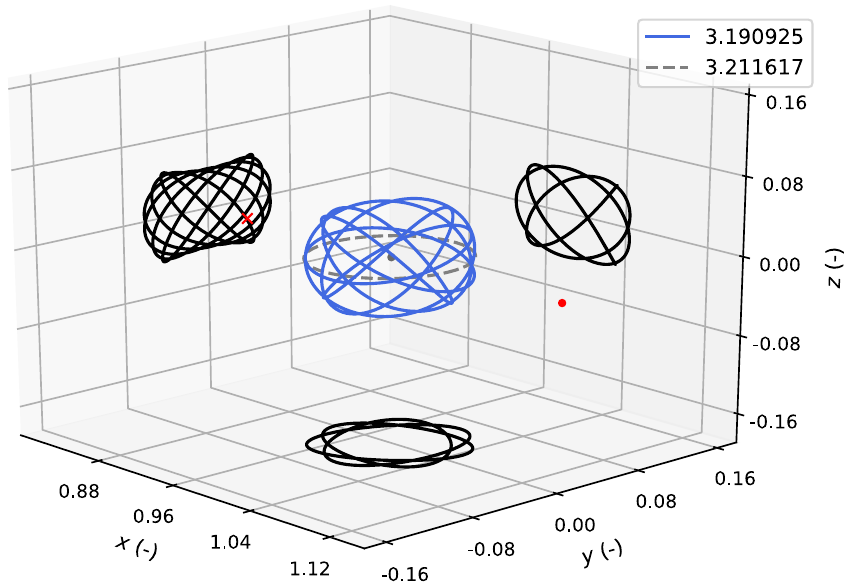}
		\end{tabular}
		\caption{7:8 DRO $\leftrightarrow$ 1:6 $LPO_1$ bridge families.}
	\end{subfigure}
	
	\vspace{0.8em}
	
	\begin{subfigure}{0.48\textwidth}
		\centering
		\begin{tabular}{cc}
			\includegraphics[width=.47\linewidth]{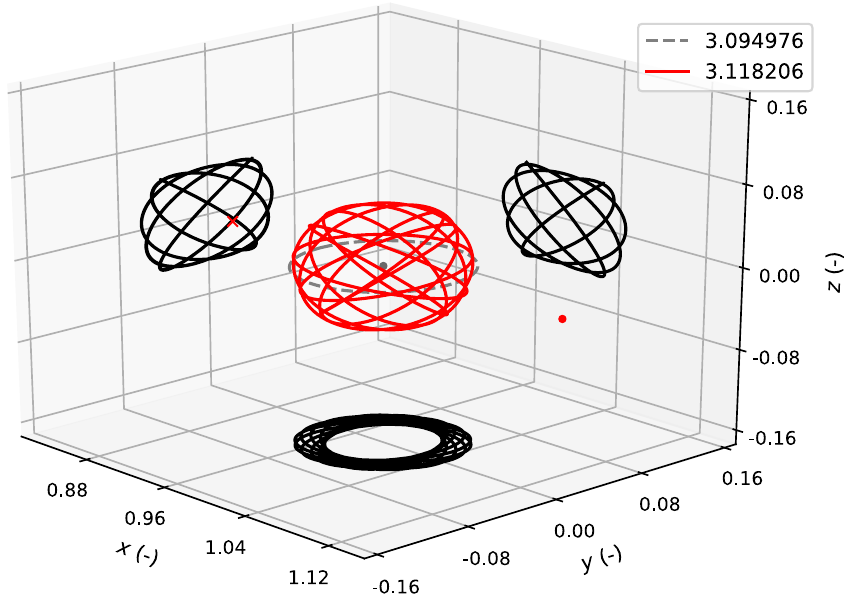} &
			\includegraphics[width=.47\linewidth]{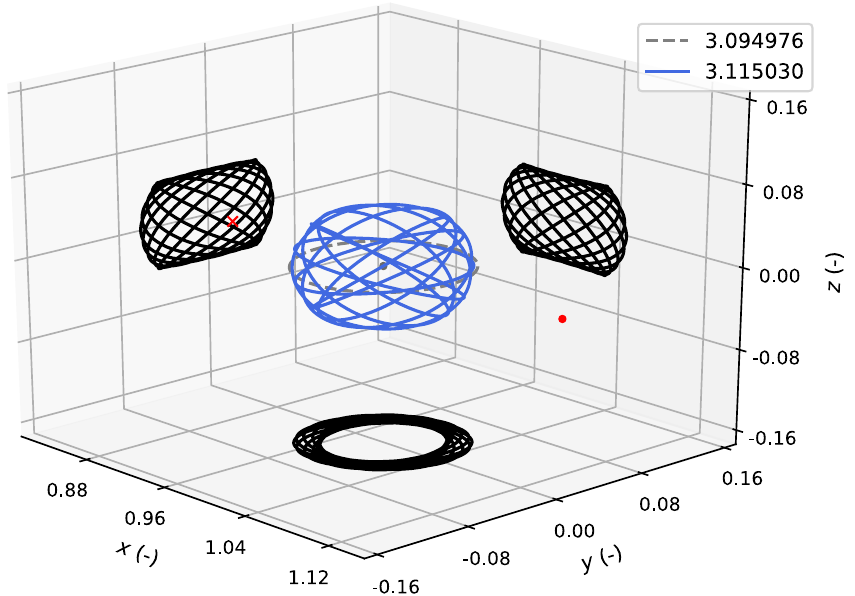} \\
			\includegraphics[width=.47\linewidth]{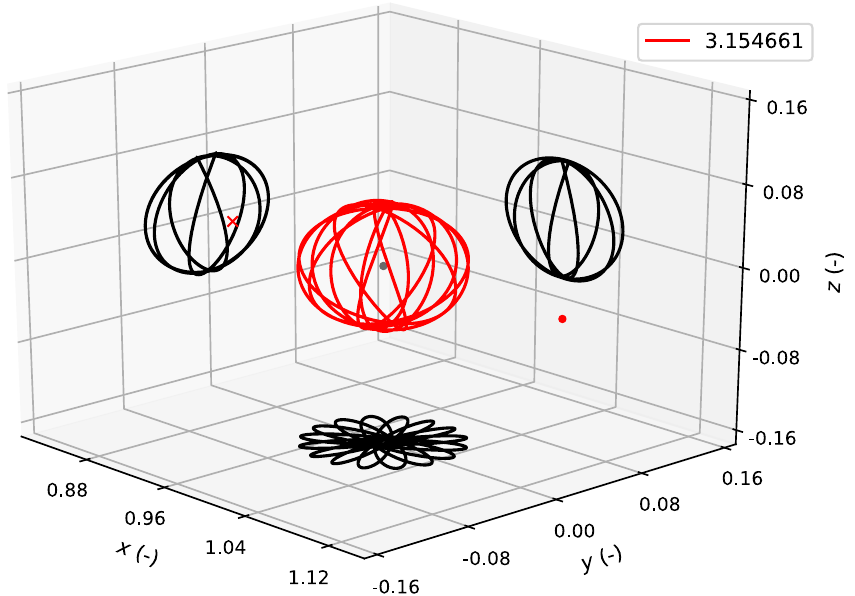} &
			\includegraphics[width=.47\linewidth]{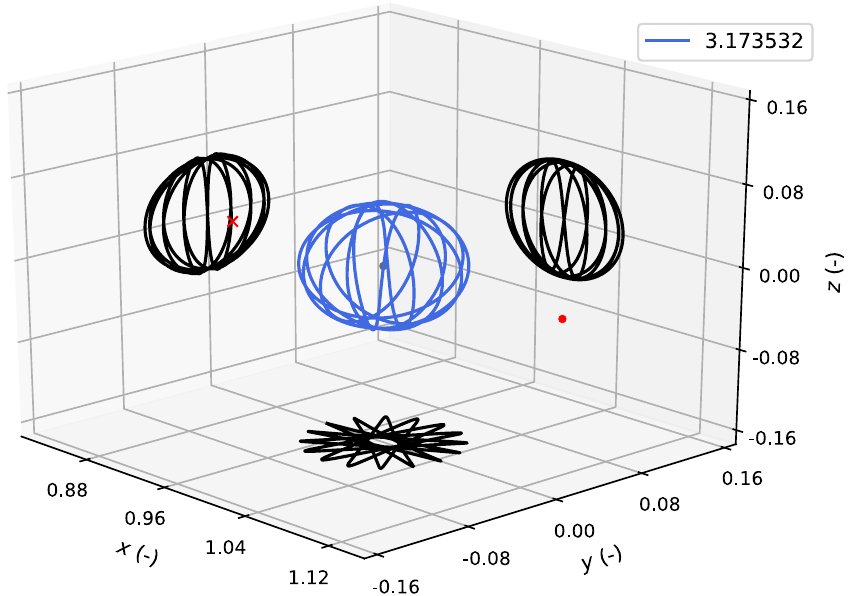} \\
			\includegraphics[width=.47\linewidth]{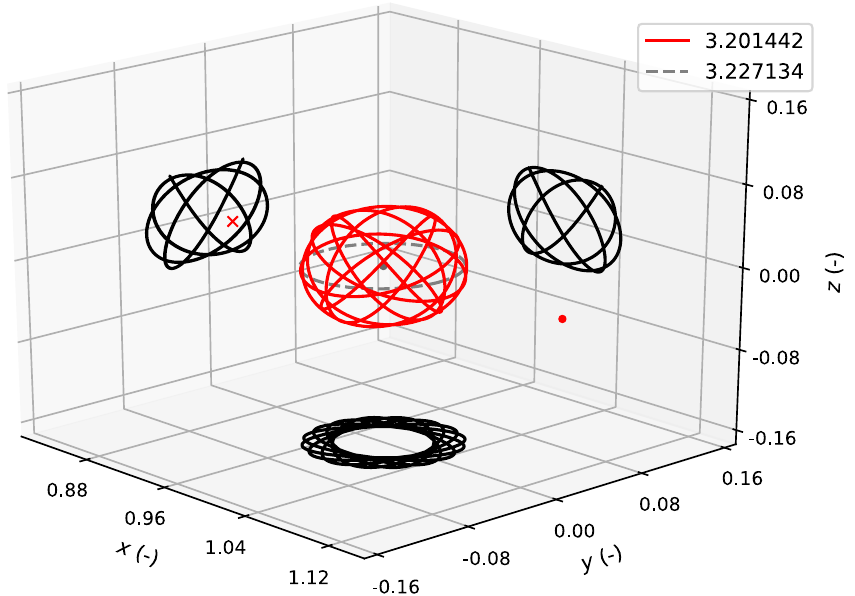} &
			\includegraphics[width=.47\linewidth]{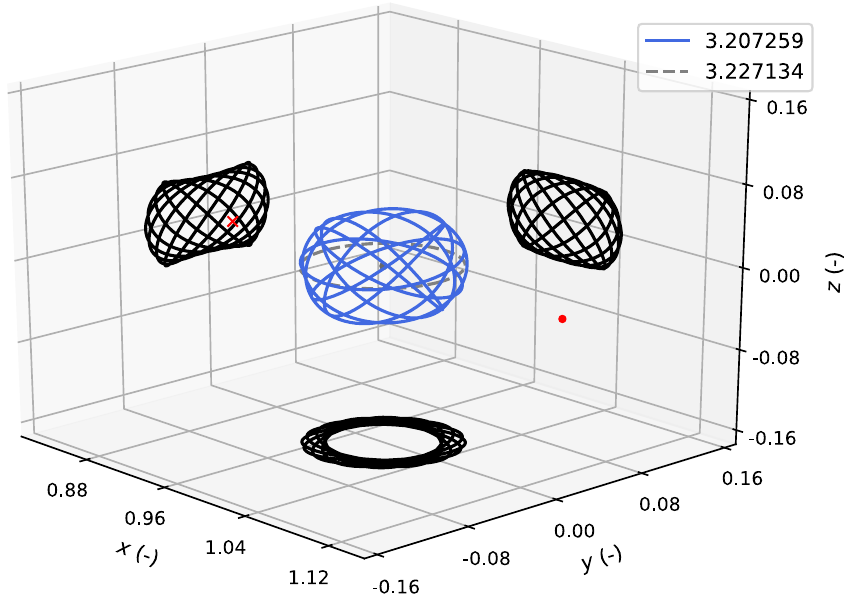}
		\end{tabular}
		\caption{8:9 DRO $\leftrightarrow$ 1:7 $LPO_1$ bridge families.}
	\end{subfigure}
	\hfill
	\begin{subfigure}{0.48\textwidth}
		\centering
		\begin{tabular}{cc}
			\includegraphics[width=.47\linewidth]{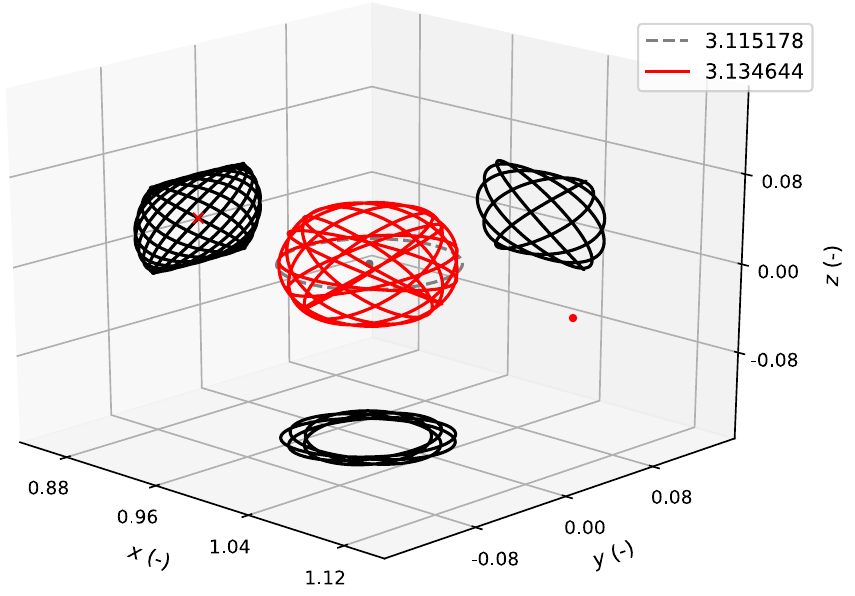} &
			\includegraphics[width=.47\linewidth]{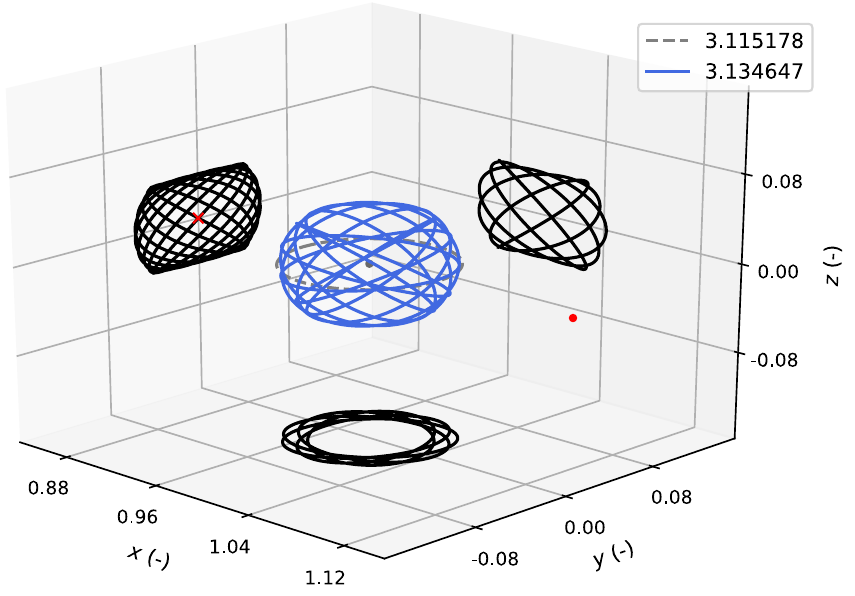} \\
			\includegraphics[width=.47\linewidth]{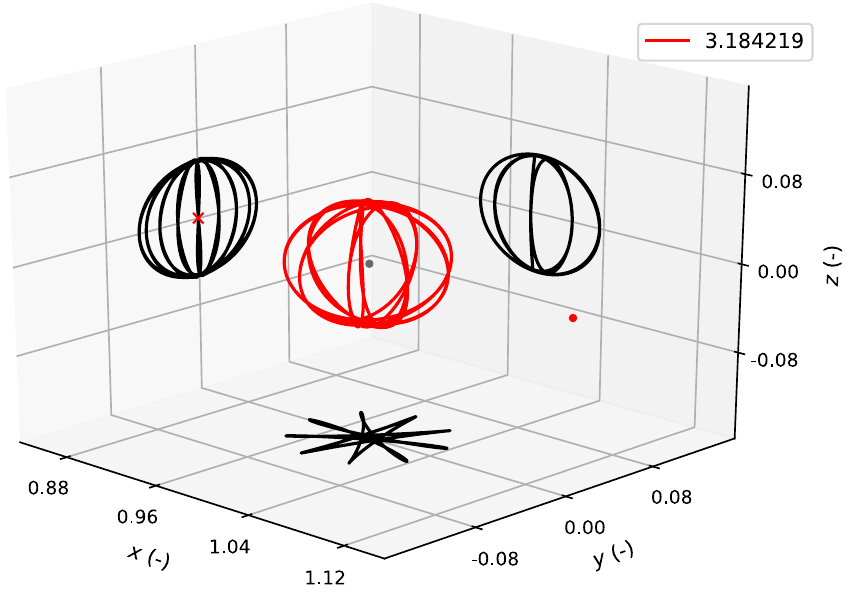} &
			\includegraphics[width=.47\linewidth]{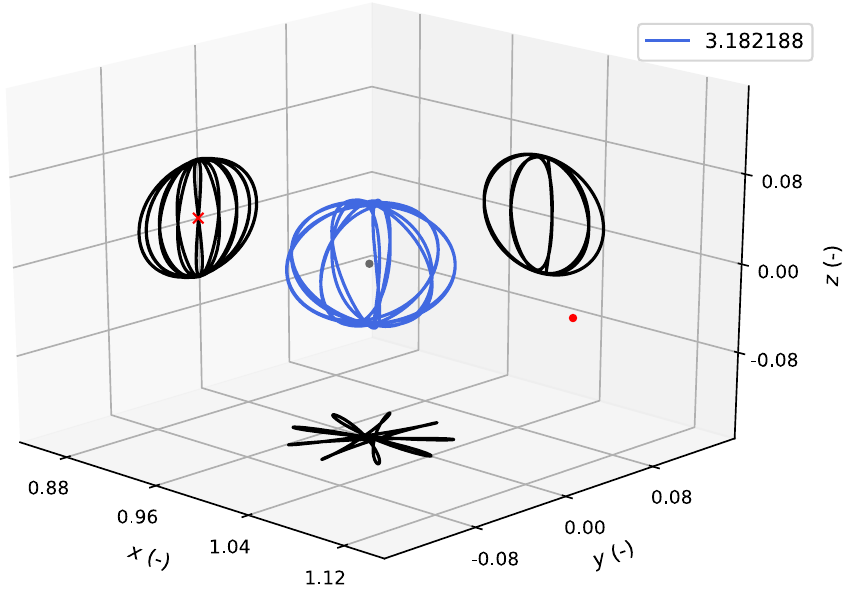} \\
			\includegraphics[width=.47\linewidth]{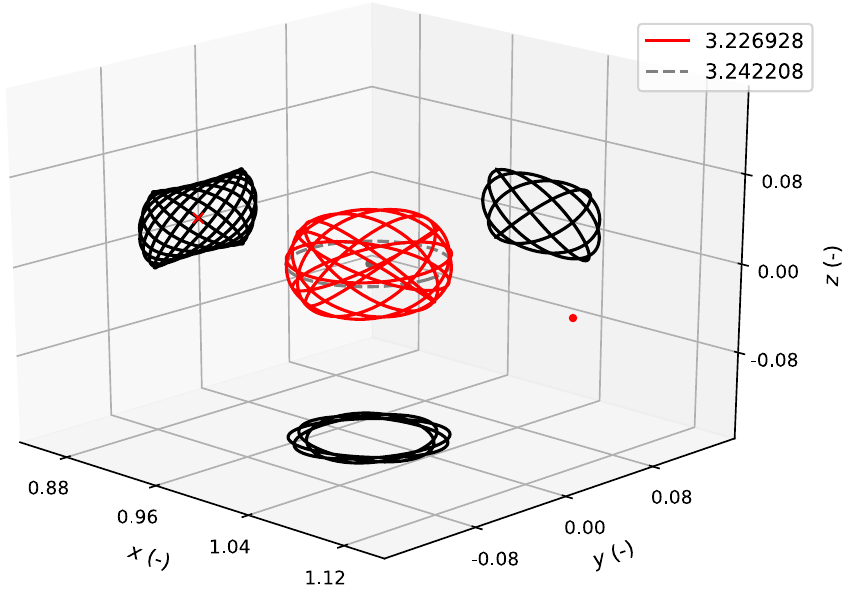} &
			\includegraphics[width=.47\linewidth]{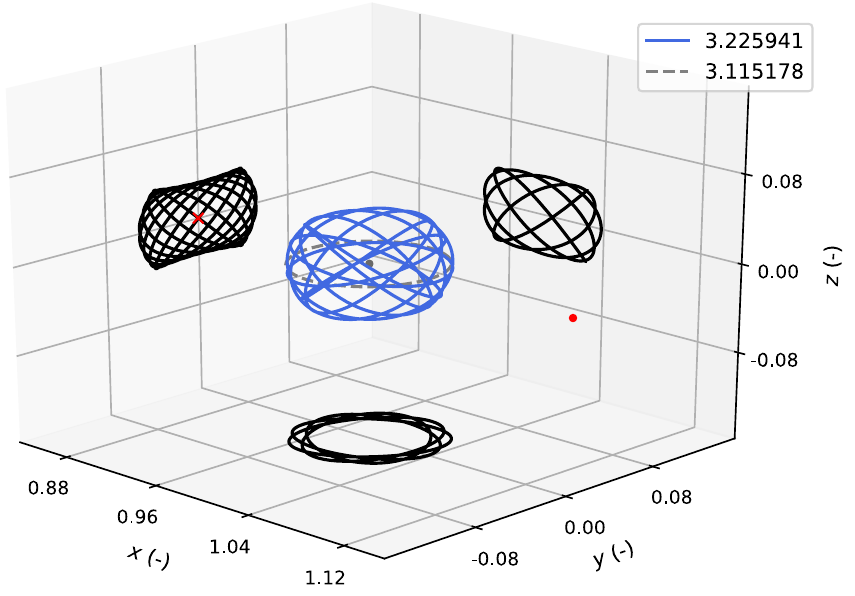}
		\end{tabular}
		\caption{9:10 DRO $\leftrightarrow$ 1:8 $LPO_1$ bridge families.}
	\end{subfigure}
	
	\caption{Representative orbits of the bifurcated families from the higher-order VSR bifurcations of the near-Moon DRO family $(q=7,8,9,10)$.\ In each case, the left column shows the red family and the right column the blue family.\ The families are continued from top to bottom, starting at the corresponding DRO (gray dashed) and terminating at the associated $LPO_1$ (gray dashed).\ For the odd-period resonances, (a) and (c), the red and blue families consist of $\rho_{xz}$- and $\rho_x$-symmetric orbits, respectively, while for the even-period resonances, (b) and (d), both families consist of doubly symmetric orbits.\ Within the plots:\ \textcolor{gray}{$\bullet$}, \textcolor{red}{$\times$} and \textcolor{red}{$\bullet$} denotes the Moon and the libration points $L_1$ and $L_2$, respectively.}
	\label{figure_plot_5}
\end{figure}

The global bifurcation structure obtained from the numerical continuation is summarized in Figures~\ref{figure_bif_1} and~\ref{figure_bif_2}, while representative orbits are shown in Figures~\ref{figure_plot_3}, \ref{figure_plot_4} and~\ref{figure_plot_5}.\ These bifurcation diagrams constitute one of the main results of the present work.\ Rather than producing isolated spatial families, every VSR bifurcation generates a bridge family, revealing how the DRO family is connected to other fundamental periodic orbit families of the Earth--Moon CR3BP.\ The detailed numerical characteristics of all continued families are collected in the Appendix B (see Tables~\ref{data_3}--\ref{data_8}).\ The CZ indices assigned to the numerically continued families provide a topological classification of the resulting bifurcation network and verify the invariance of the Euler characteristic of the local Floer homology at each bifurcation point.\

The first bridge families, associated with the 4:5 DRO, differ qualitatively from the other higher-order multiplicities and is shown in the bifurcation diagram in Figure~\ref{figure_bif_1}.\ At the 4:5 DRO, continuation with respect to $\Delta z$ generates the family of $\rho_{xz}$-symmetric periodic orbits, while continuation with respect to $\Delta p_z$ produces the family of $\rho_x$-symmetric periodic orbits.\ The former possesses elliptic critical Floquet multipliers, while the latter possesses positive hyperbolic critical Floquet multipliers (see Table~\ref{data_3}).\ Accordingly, the red and blue families consist of $\rho_{xz}$- and $\rho_x$-symmetric periodic orbits, respectively, together with $\sigma$-symmetric counterparts.\ Upon continuation, the red family terminates at the period-doubling bifurcation of the Butterfly family, while the blue family terminates at the period-tripling bifurcation of the $LPO_2$, as illustrated in Figure~\ref{figure_bif_1}.\ Therefore, these bridges provide a bifurcation network connection between three basic periodic orbit families:\ the DRO, the Butterfly and the planar prograde $LPO_2$ family.\ Representative members of both families are shown in Figure~\ref{figure_plot_3}.\ We notice that similar results were documented previously within the Sun--Jupiter CR3BP \cite{robin_markellos}.\ We further notice that there is a second family that bifurcates from the 1:3 $LPO_2$ (green family in Figure~\ref{figure_bif_1}).\ Since this work focuses on the VSR bifurcations originating from the DRO family, this secondary family is not continued further.\

The remaining five VSR bifurcations $(q=6,7,8,9,10)$ exhibit a remarkably regular global bifurcation structure, that is summarized by the bifurcation diagram in Figure~\ref{figure_bif_2}.\ Every VSR bifurcation produces bridge families connecting the DRO family with the $LPO_1$ family.\ More precisely, the $q$-fold covered DRO family is connected to the $(q-2)$-fold covered $LPO_1$ family, giving rise to the bifurcation sequence with respect to the self-resonance ratios:
$$ 5\text{:}6 \leftrightarrow 1\text{:}4,\quad 6\text{:}7 \leftrightarrow 1\text{:}5,\quad 7\text{:}8 \leftrightarrow 1\text{:}6,\quad 8\text{:}9 \leftrightarrow 1\text{:}7,\quad 9\text{:}10 \leftrightarrow 1\text{:}8. $$
This systematic pattern, established through numerical continuation, is underpinned by the corresponding CZ indices.\ As already noticed after~\eqref{indices_basic}, the planar and vertical winding numbers of the $LPO_1$ family exceed those of the DRO family by two.\ This means that in the considered regimes, the CZ indices of the $q$-fold covered DRO family coincide with those of the $(q-2)$-fold covered $LPO_1$, as can be observed in Figure~\ref{figure_bif_2}.\ Therefore, this bridge pattern reveals a systematic geometric structure that is also reflected in the corresponding symplectic invariants.\ Although similar bridge families were previously reported for the Hill problem~\cite{aydin_cz,aydin_batkhin}, Figure~\ref{figure_bif_2} demonstrates that similar global bifurcation network persists in the Earth--Moon CR3BP.\

Beyond the systematic connectivity between the DRO and $LPO_1$ families, the bifurcation diagram in Figure~\ref{figure_bif_2} reveals a remarkable bifurcation topology.\ The even-period bridge families possess similar bifurcation results.\ Likewise, the odd-period bridge families exhibit the same topological structure.\ This suggests that the parity of the multiplicity $q$ within the self-resonance determines the qualitative organization of the bridge families in the vicinity of the Moon.\ We further observe from the orbit plots in Figures~\ref{figure_plot_4} and~\ref{figure_plot_5} associated to the even-period bridge families, that the red and blue families appear to be related by an approximate reflection with respect to the $yz$-plane through the Moon.\ Although this symmetry is not a symmetry of the CR3BP, the corresponding orbit geometries are notably similar after such a reflection, indicating that the two families may be regarded as geometric counterparts of one another.\ By contrast, the orbits in Figure~\ref{figure_plot_5} related to the odd-period bridge families remain distinguished by their different $\rho_{xz}$- and $\rho_x$-symmetry types, producing geometrically different orbit families despite having the same bifurcation topology.\

Furthermore, Figure~\ref{figure_plot_3} shows a spatial 1:1 Earth--Moon resonant orbit.\ We observe from Figure~\ref{plot_3_c} that this resonant orbit is in the Earth-centered inertial frame of waveform with phases above and below the ecliptic near the lunar orbit, while in the Moon-centered inertial frame the orbit's geometry is of ring-shaped around the Moon.\ This resonant orbit has CZ index~15 with elliptic and positive hyperbolic Floquet multipliers of the form $\{ e^{\pm i \theta}, \lambda, \lambda^{-1} \}$, where $\theta \approx 5.998$ (negative Krein) and $\lambda \approx 29.695$.\ Similar spatial 1:1 Earth--Moon resonant orbits, however, are encountered along every computed bridge families connecting the DRO family to the planar prograde families, demonstrating that these bifurcation connections consistently pass through dynamically important regions of the phase space, that is relevant for the design of transfer trajectories \cite{vaquero}.\ Moreover, (linearly) stable spatial periodic orbits are concentrated on the red families, whose orbits are stable near the DRO family for all bridge families and, for the 7:8, 8:9, and 9:10 VSR bifurcations, also near the $LPO_1$ family.\ The only stable orbits within the blue families occur for the 4:5 and 5:6 VSR bifurcations, near the $LPO_2$ and $LPO_1$ families, respectively.\

\subsection{VSR Bifurcations from DROs around the Moon, $L_1$ and $L_2$}

\begin{figure}[t!]
	\centering
	\begin{subfigure}[]{0.42\textwidth}
		\centering
		\includegraphics[width=0.72\linewidth]{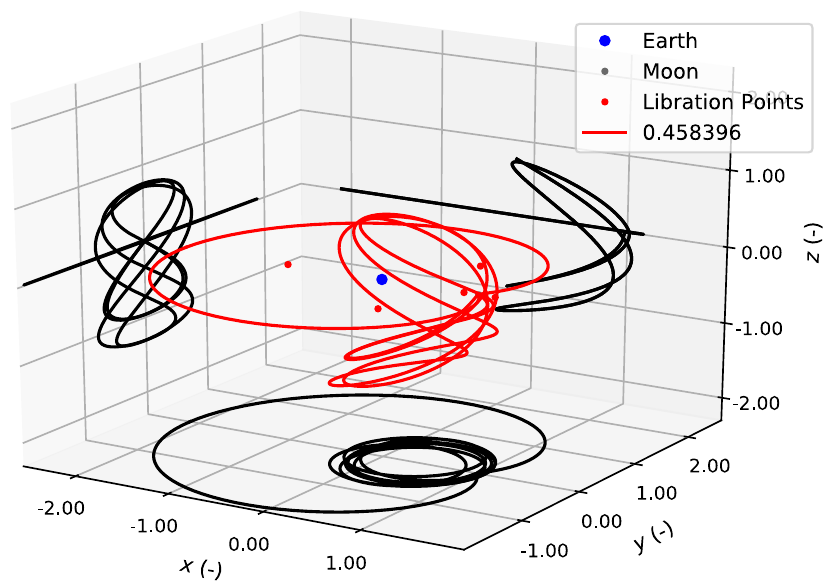}	
		\vfill
		\includegraphics[width=0.72\linewidth]{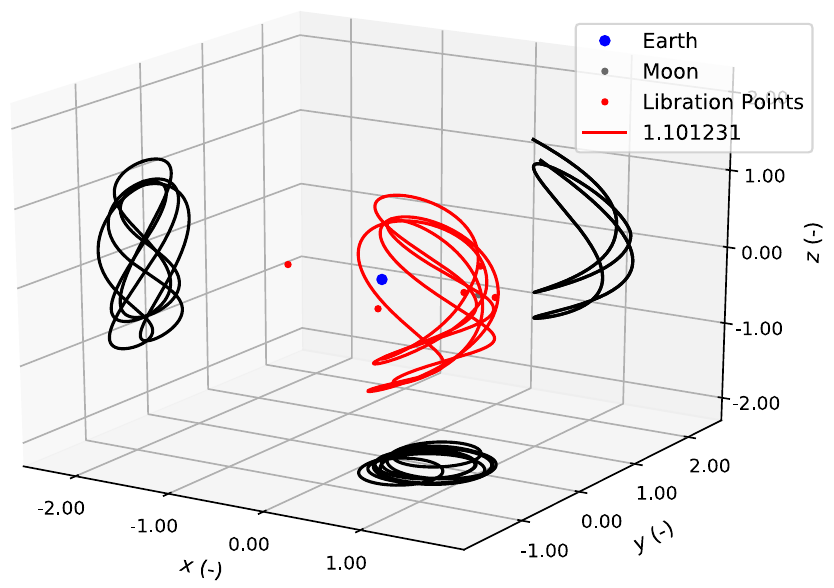}
		\vfill
		\includegraphics[width=0.72\linewidth]{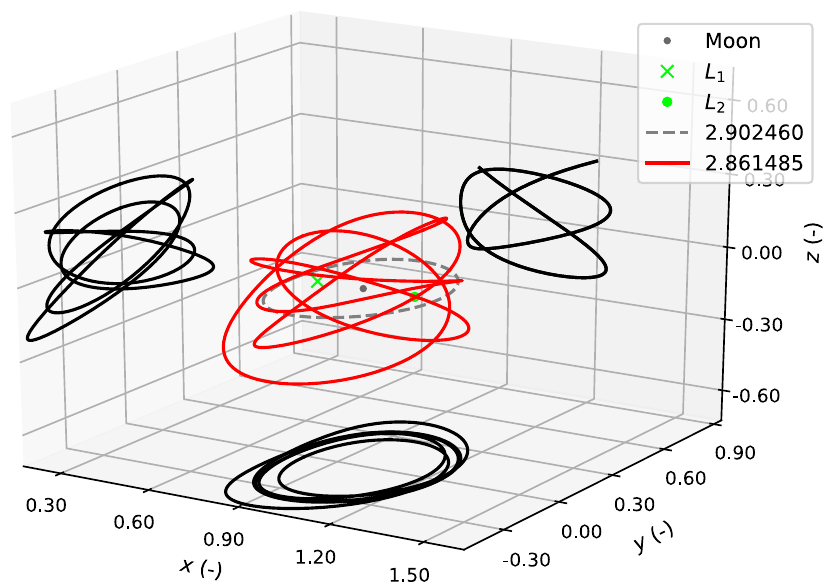}
		\caption{Representative $\rho_{xz}$-symmetric orbits of the red family, continued from bottom to top.}
		\label{plot_6_a}
	\end{subfigure}
	\hfill
	\begin{subfigure}[]{0.57\textwidth}
		\centering
		\includegraphics[width=\linewidth]{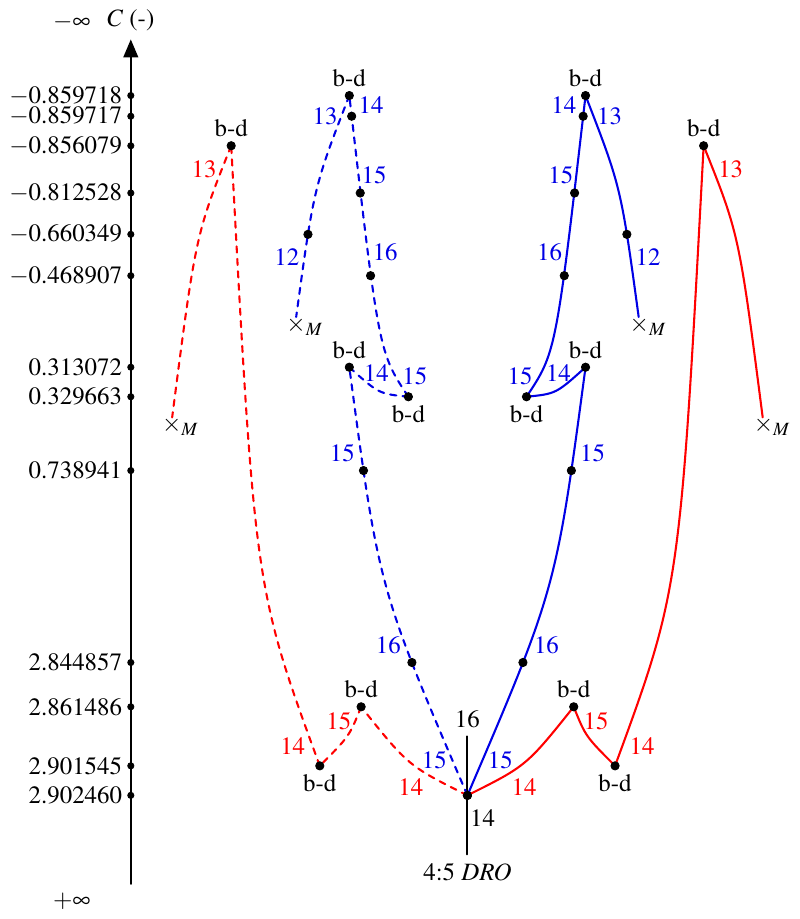}
		\caption{Bifurcation diagram associated to the 4:5 DRO that moves around the Moon, $L_1$ and $L_2$.}
		\label{plot_6_b}
	\end{subfigure}
	\vspace{0.5em}
	
	\begin{subfigure}[t]{\textwidth}
		\centering
		\includegraphics[width=0.32\linewidth]{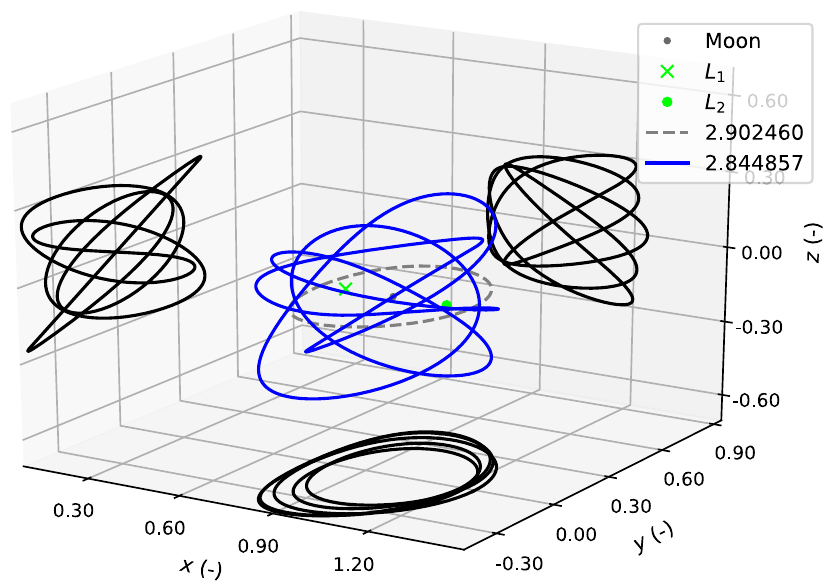}
		\hfill
		\includegraphics[width=0.32\linewidth]{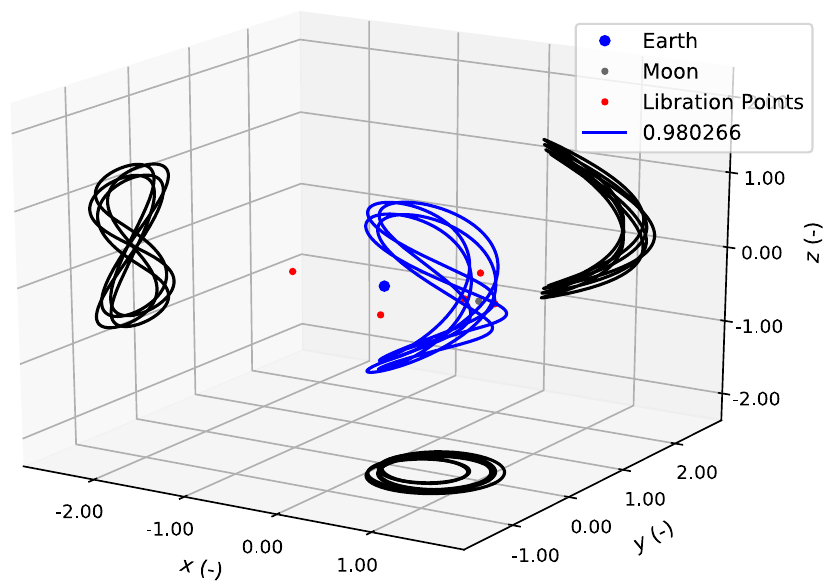}
		\hfill
		\includegraphics[width=0.32\linewidth]{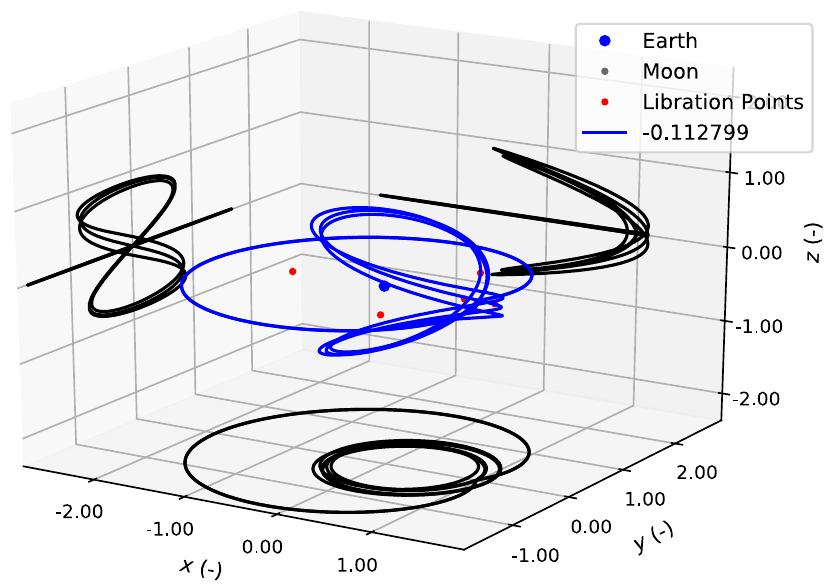}
		\caption{Representative $\rho_x$-symmetric orbits of the blue family, continued from left to right.}
		\label{plot_6_c}
	\end{subfigure}
	\vspace{0.5em}
	
	\begin{subfigure}[t]{0.49\textwidth}
		\centering
		\includegraphics[width=0.49\linewidth]{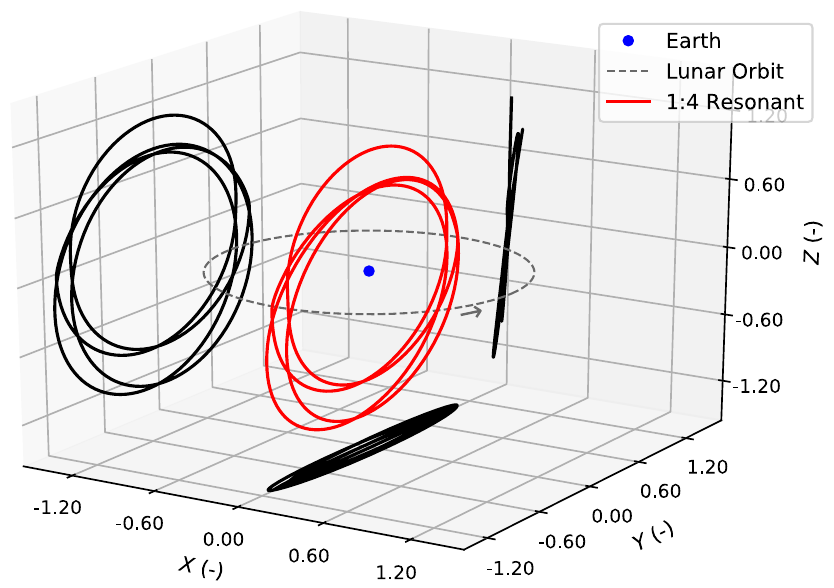}
		\hfill
		\includegraphics[width=0.49\linewidth]{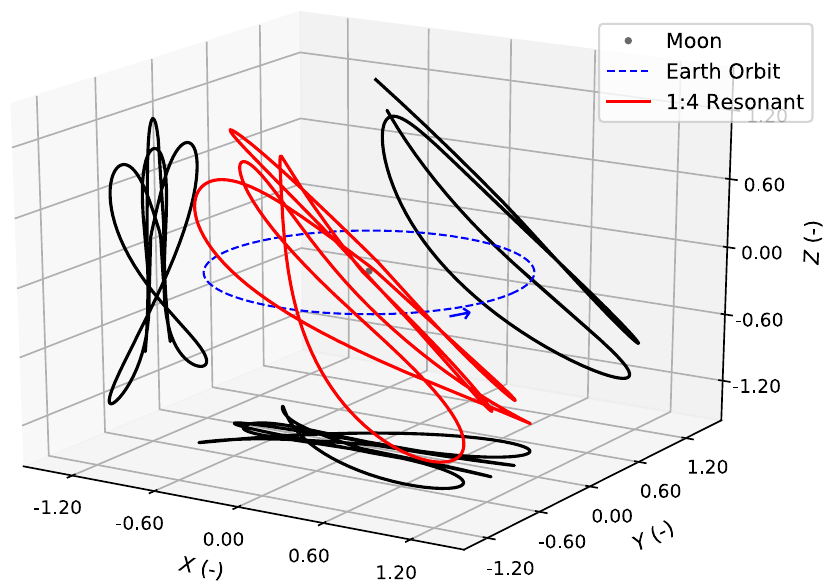}
		\caption{Spatial 1:4 Earth--Moon resonant periodic orbit of the red family, shown in the Earth-centered inertial frame (left) and in the Moon-centered inertial frame (right).\ The corresponding orbit in the rotating frame is shown in the middle of (a).}
		\label{plot_6_d}
	\end{subfigure}
	\hfill
	\begin{subfigure}[t]{0.49\textwidth}
		\centering
		\includegraphics[width=0.49\linewidth]{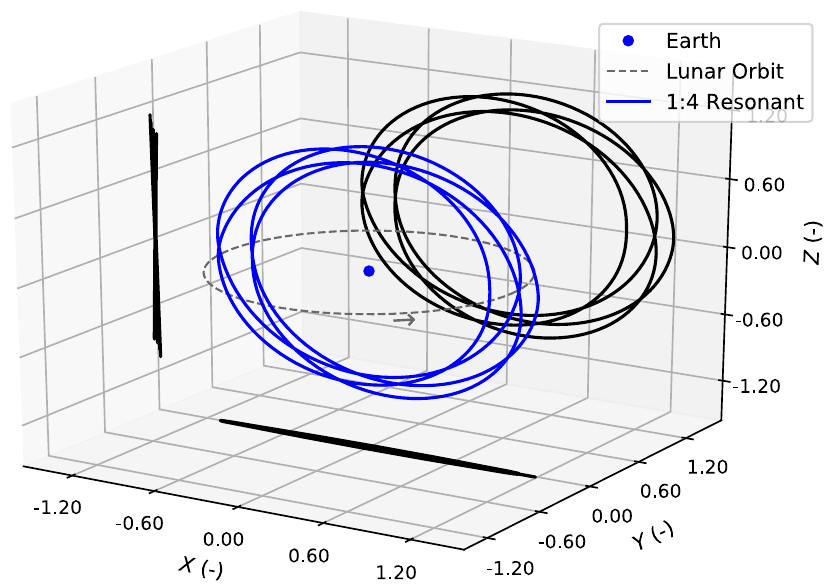}
		\hfill
		\includegraphics[width=0.49\linewidth]{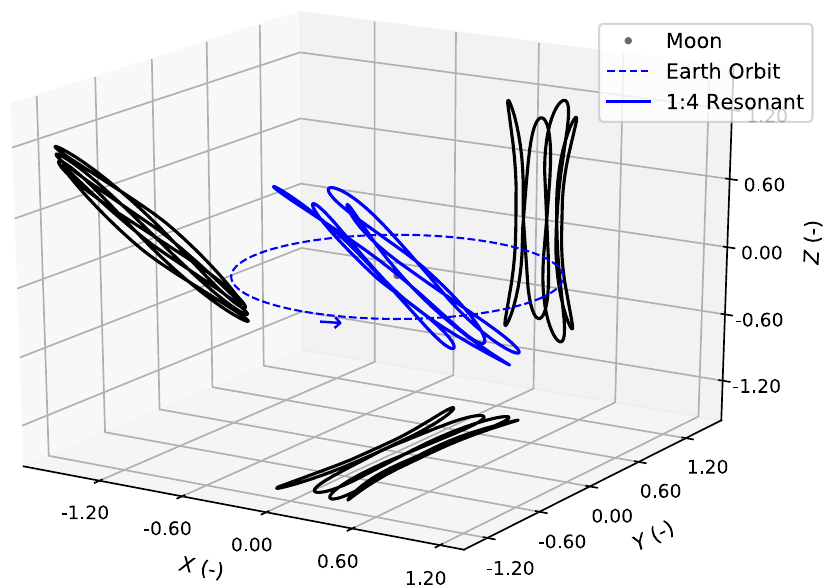}
		\caption{Spatial 1:4 Earth--Moon resonant periodic orbit of the blue family, shown in the Earth-centered inertial frame (left) and in the Moon-centered inertial frame (right).\ The corresponding orbit in the rotating frame is shown in the middle of (c).}
		\label{plot_6_e}
	\end{subfigure}
	
	\caption{VSR bifurcations from the 4:5 DRO that moves around the Moon, $L_1$ and $L_2$.}
	\label{figure_plot_6}
\end{figure}

\begin{figure}[t!]
	\centering
	\includegraphics[width=1\linewidth]{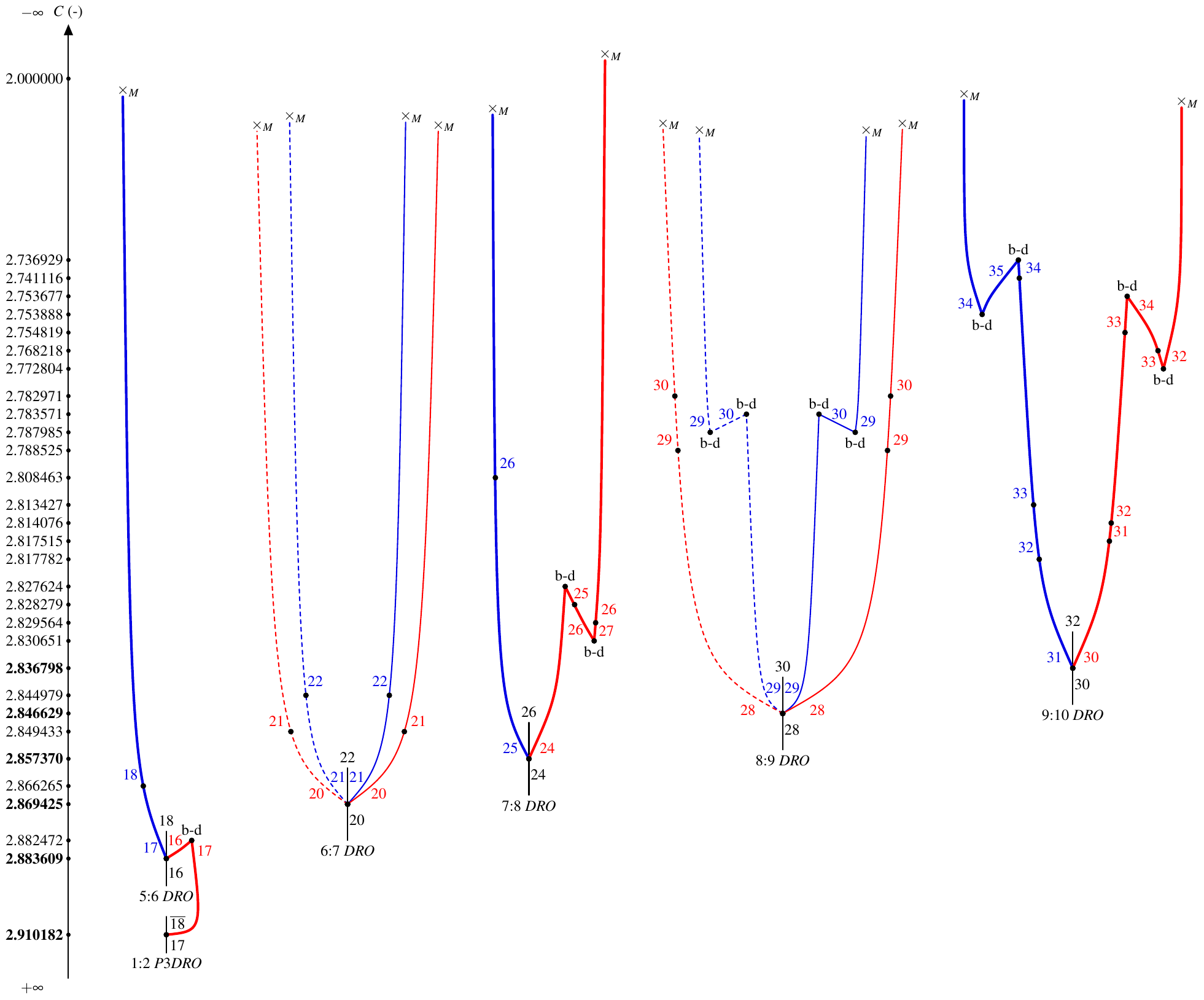}
	\caption{Bifurcation diagram associated to the VSR bifurcations from $(q-1)\text{:}q$ DRO that move around the Moon, $L_1$ and $L_2$, where $q=6,7,8,9,10$.\ Bold Jacobi constants on the left correspond to the VSR orbits.\ Only the red family that branches off from 5:6 DRO forms a bridge to 1:2 P3DRO.\ Orbits of all the other bifurcated families become of long-period multi-excursion type.}
	\label{figure_bif_4}
\end{figure}

We now investigate the next six VSR bifurcations associated with the DRO family with the same multiplicities $(q=5,6,7,8,9,10)$ in the second geometric regime, where the DROs move not only around the Moon but also around the collinear libration points $L_1$ and $L_2$.\ In contrast to the near-Moon DROs, these VSR bifurcations are continued in decreasing direction of the Jacobi constant and exhibit a fundamentally different global organization (see bifurcation diagrams in Figures~\ref{plot_6_b} and~\ref{figure_bif_4}).\ While the near-Moon VSR bifurcations generate an interconnected network of bridge families linking the DRO family with the Butterfly and planar prograde families, only one bridge family is found in the present second geometric regime.\ The remaining bifurcated families evolve into a common class of highly unstable orbits that combine extended near-planar motion in the ecliptic with multiple spatial excursions, during which the orbits repeatedly alternate between large and small distances from the Moon.\ The projections of these trajectories onto the coordinate planes exhibit characteristic insect-like patterns, resembling Bumble Bee, Hoverfly, or Dragonfly shapes.\ We refer to this class of spatial periodic orbits as \textit{long-period multi-excursion orbits}.\

The continuation of the VSR bifurcations from the 4:5 DRO is summarized in Figure~\ref{plot_6_b}, and the numerical data are collected in Table~\ref{data_9}.\ Both of the bifurcated red and blue families undergo a sequence of three birth--death bifurcations before evolving toward orbits of long-period multi-excursion type.\ The corresponding bifurcation diagram is illustrated in Figure~\ref{plot_6_b} with representative orbits that are shown in Figures~\ref{plot_6_a} and~\ref{plot_6_c}.\ Only the red family provides stable solutions that are near the 4:5 DRO.\ The final computed long-period multi-excursion orbits within the red $(C \approx 0.458396)$ and blue $(C \approx -0.112799)$ families have minimum lunar distances of 0.006423 and 0.007918 (normalized units), respectively.\ Furthermore, along both families, we find 1:4 Earth--Moon resonant orbits shown in Figures~\ref{plot_6_d} and~\ref{plot_6_e}, that have large negative and positive hyperbolic Floquet multipliers.\

The remaining VSR bifurcations, corresponding to the multiplicities $q=6,7,8,9,10$, are summarized in the bifurcation diagram in Figure~\ref{figure_bif_4}, with numerical data that can be found in Tables~\ref{data_10}--\ref{data_14}.\ Only the red family bifurcating from the 5:6 DRO forms a bridge, terminating at the period-doubling bifurcation of the P3DRO family.\ The continuation of all remaining red and blue families illustrate that the typical continuation behavior of this geometric region is the evolution toward orbits of long-period multi-excursion type, rather than the formation of bridges connecting basic periodic orbit families.\ In all cases of the VSR bifurcations, stable solutions near the corresponding DRO can be identified within the red family.\

\begin{figure}[t!]
	\centering
	
	\begin{subfigure}{\textwidth}
		\centering
		\includegraphics[width=0.32\textwidth]{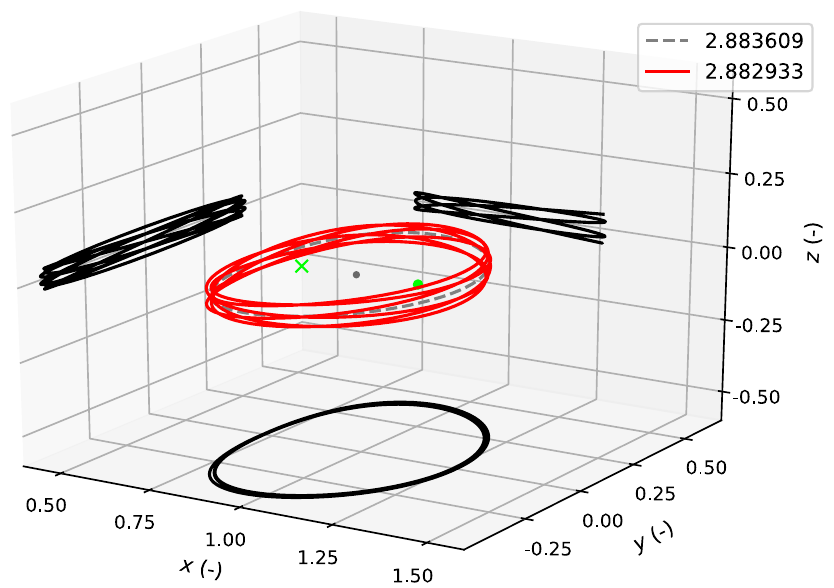}
		\hfill
		\includegraphics[width=0.32\textwidth]{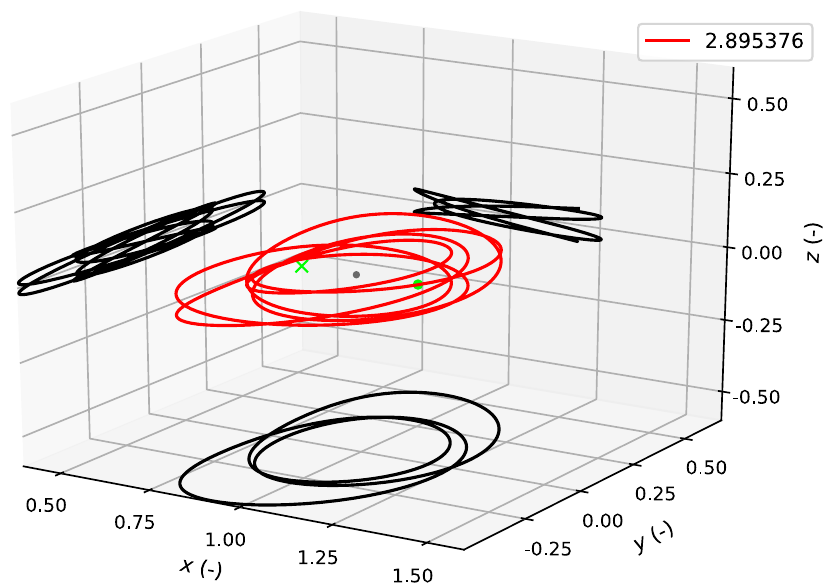}
		\hfill
		\includegraphics[width=0.32\textwidth]{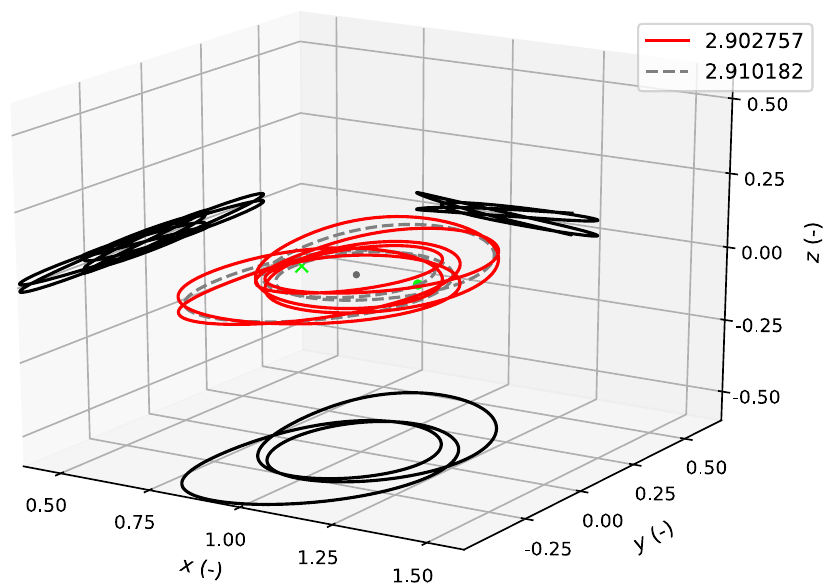}
		\caption{Representative orbits of the red family, continued from left to right.\ They bifurcate from the 5:6 DRO (gray dashed) and terminate at the 1:2 P3DRO (gray dashed).\ The middle orbit is a 1:4 Earth--Moon resonant orbit.}
		\label{plot_7_a}
	\end{subfigure}
	\vspace{0.5em}
	
	\begin{subfigure}{\textwidth}
		\centering
		\includegraphics[width=0.32\textwidth]{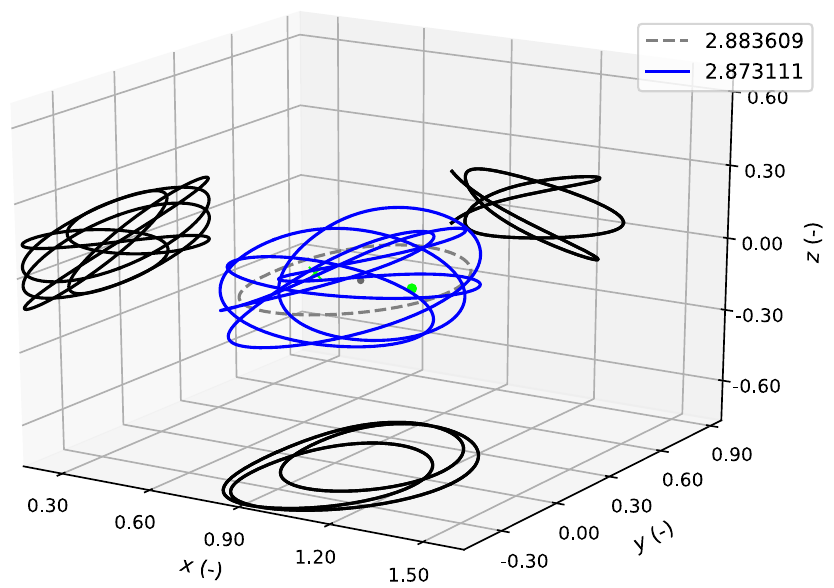}
		\hfill
		\includegraphics[width=0.32\textwidth]{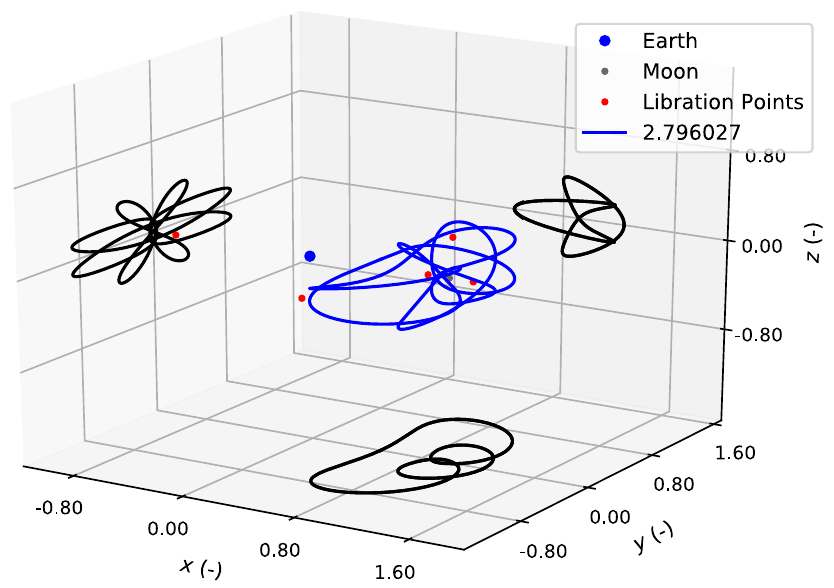}
		\hfill
		\includegraphics[width=0.32\textwidth]{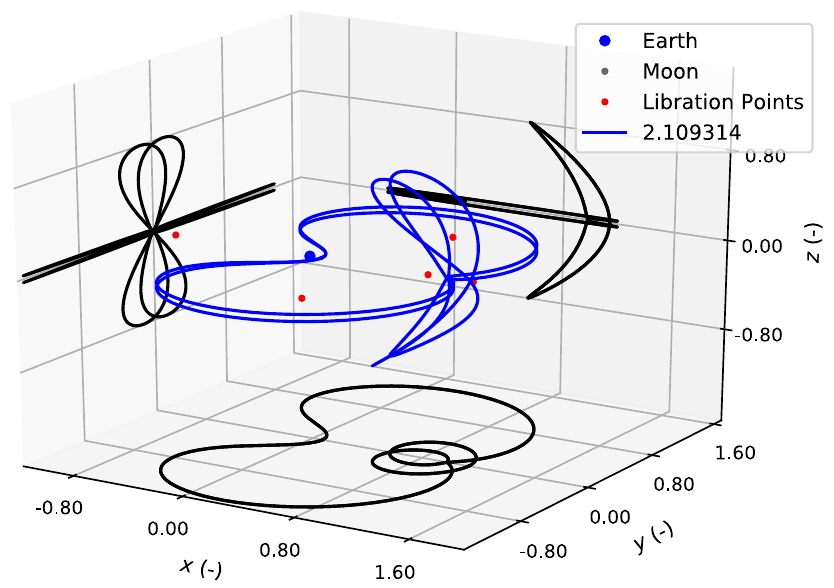}
		\caption{Representative orbits of the blue family, continued from left to right.\ They bifurcate from the 5:6 DRO (gray dashed) and are continued until being of long-period multi-excursion type.}
		\label{plot_7_b}
	\end{subfigure}
	\vspace{0.5em}
	
	\begin{subfigure}{\textwidth}
		\centering
		\includegraphics[width=0.45\textwidth]{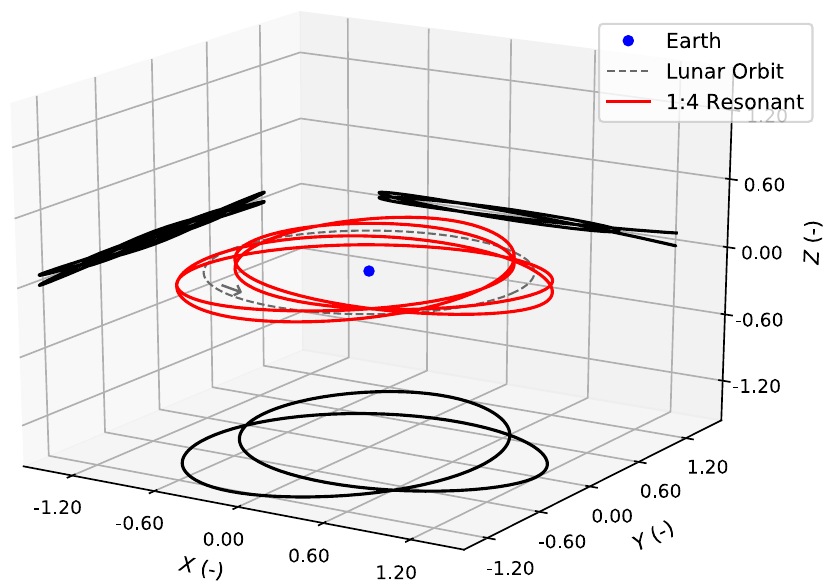}
		\includegraphics[width=0.45\textwidth]{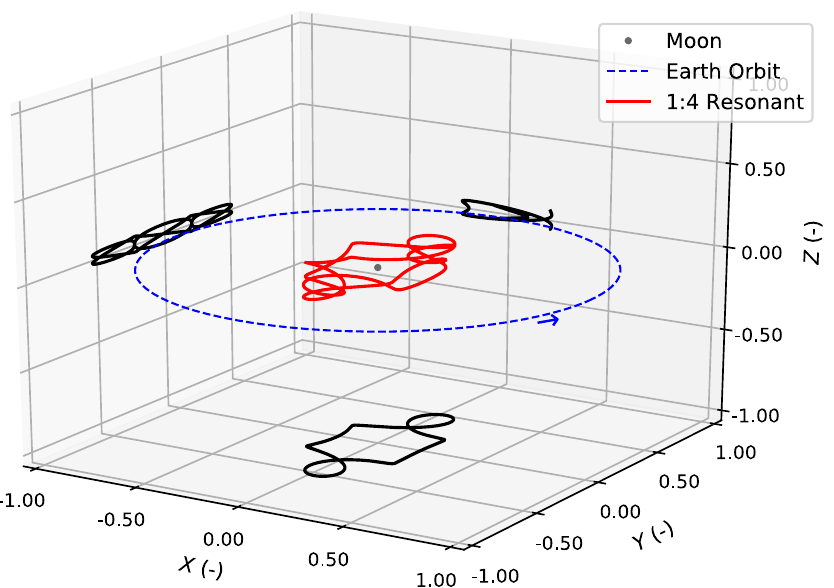}
		\caption{The 1:4 Earth--Moon resonant orbit from (a) is shown in the Earth-centered inertial frame (left) and in the Moon-centered inertial frame (right).}
		\label{plot_7_c}
	\end{subfigure}
	
	\caption{Doubly-symmetric orbits associated to the VSR bifurcations from the 5:6 DRO that moves around the Moon, $L_1$ and $L_2$.\ Within the plots (a) and the left plot in (b):\ \textcolor{gray}{$\bullet$}, \textcolor{green}{$\times$} and \textcolor{green}{$\bullet$} denotes the Moon and the libration points $L_1$ and $L_2$, respectively.}
	\label{figure_plot_7}
\end{figure}

Figure~\ref{figure_plot_7} shows representative orbits of the VSR bifurcations from the 5:6 DRO.\ The orbits in Figure~\ref{plot_7_a} belong to the red family that provides the bridge connecting the DRO with the P3DRO family.\ We notice that a similar bridge was previously identified in the Hill problem \cite{aydin_batkhin}, where an additional bridge between the DRO and P3DRO families was also reported.\ In contrast, the present computations in the Earth--Moon CR3BP reveal only this bridge connecting 5:6 DRO and 1:2 P3DRO.\ Along this bridge, a 1:4 Earth--Moon resonant orbit is encountered, shown in Figure~\ref{plot_7_c}.\ When viewed in the Moon-centered inertial frame, the orbit continuously connects the $L_1$ and $L_2$ libration regions.\ This resonant orbit has CZ index 17 and the Floquet multipliers are elliptic and positive hyperbolic of the form $\{ e^{\pm i \theta}, \lambda, \lambda^{-1} \}$, where $\theta \approx 1.113$ and $\lambda \approx 10.901$.\ The final computed long-period multi-excursion orbit of the blue family $(C \approx 2.109314)$ in Figure~\ref{plot_7_b} has the closest approach of 0.010839 (normalized units) to the Moon.\

\begin{figure}[t!]
	\centering
	
	\begin{subfigure}{0.48\textwidth}
		\centering
		\begin{tabular}{cc}
			\includegraphics[width=.47\linewidth]{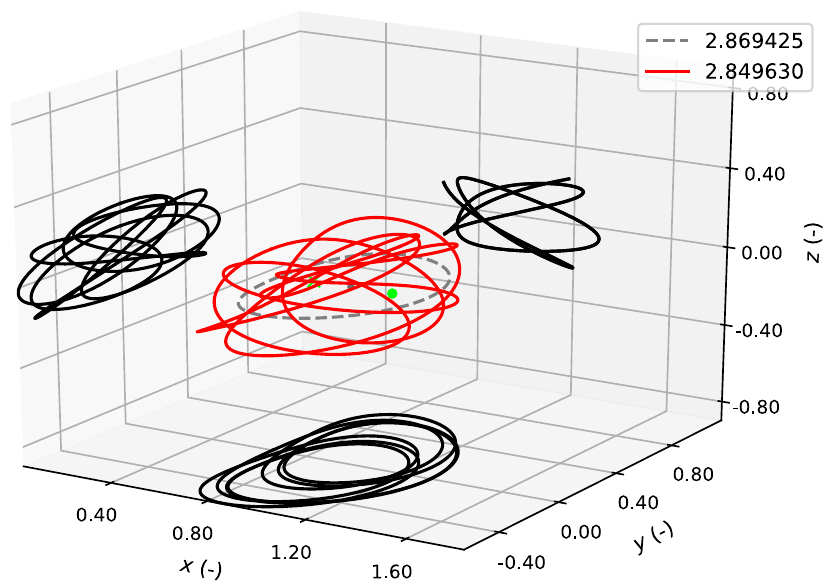} &
			\includegraphics[width=.47\linewidth]{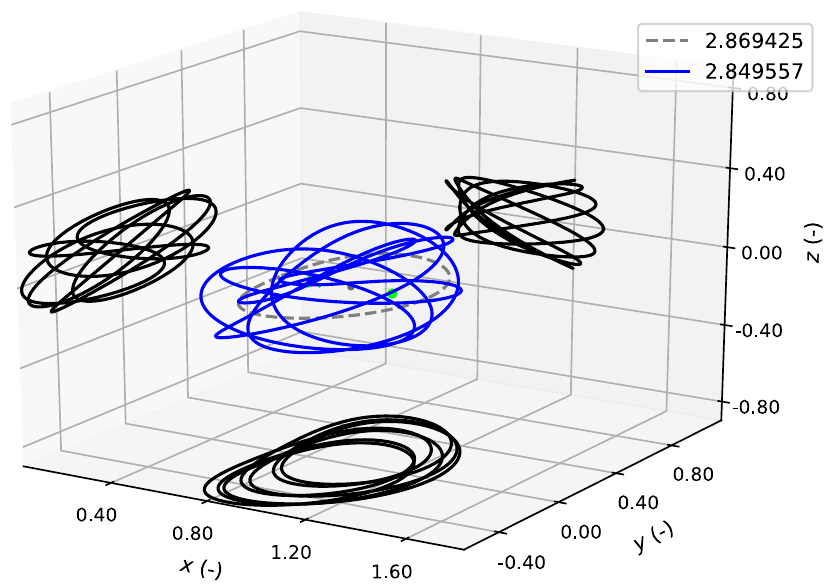} \\
			\includegraphics[width=.47\linewidth]{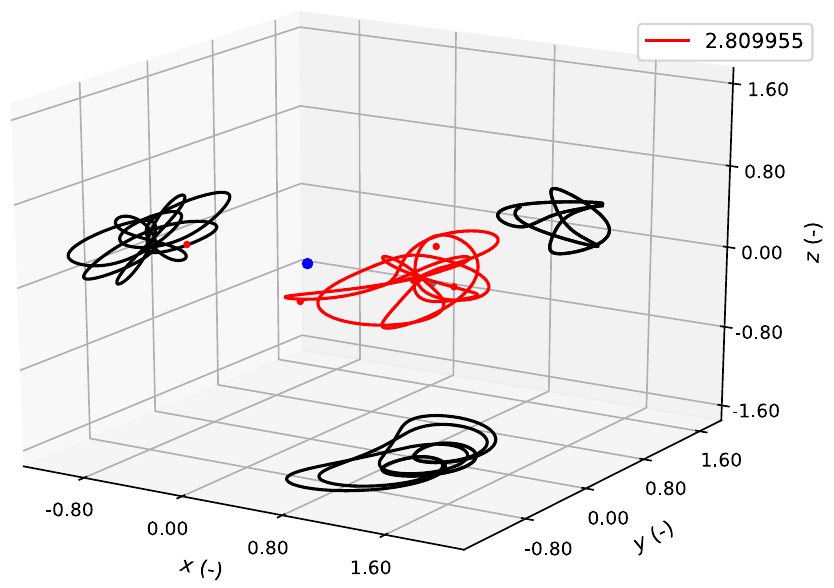} &
			\includegraphics[width=.47\linewidth]{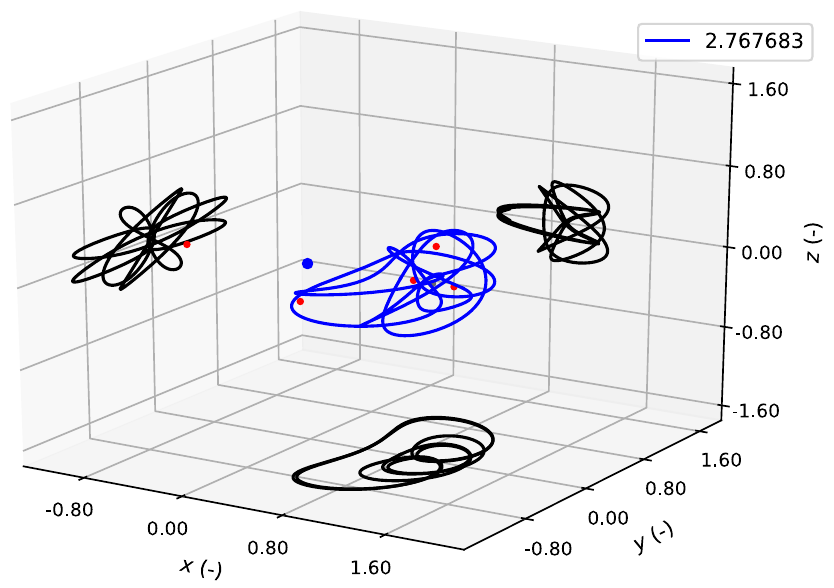} \\
			\includegraphics[width=.47\linewidth]{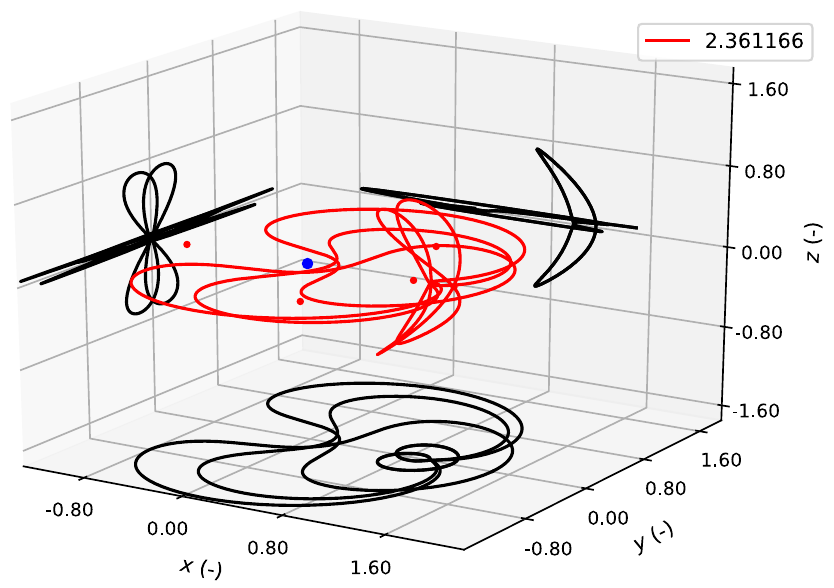} &
			\includegraphics[width=.47\linewidth]{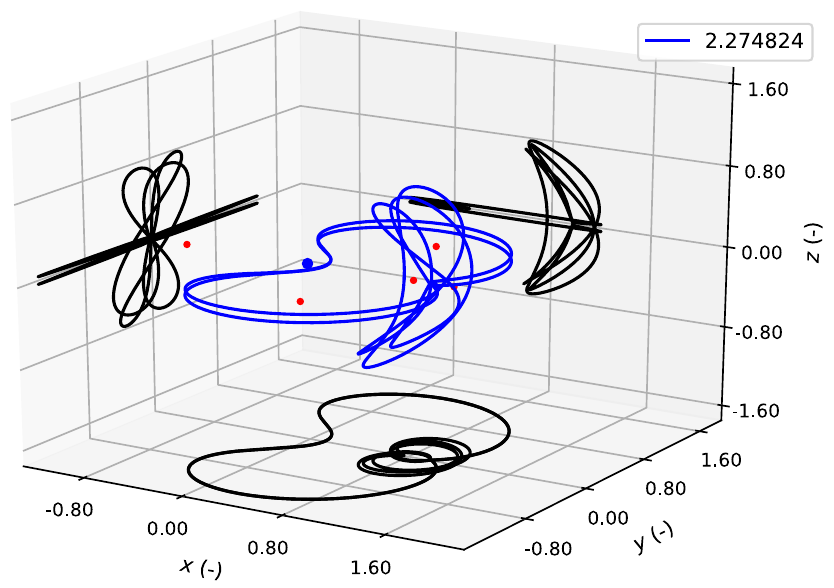}
		\end{tabular}
		\caption{VSR bifurcations from the 6:7 DRO.}
	\end{subfigure}
	\hfill
	\begin{subfigure}{0.48\textwidth}
		\centering
		\begin{tabular}{cc}
			\includegraphics[width=.47\linewidth]{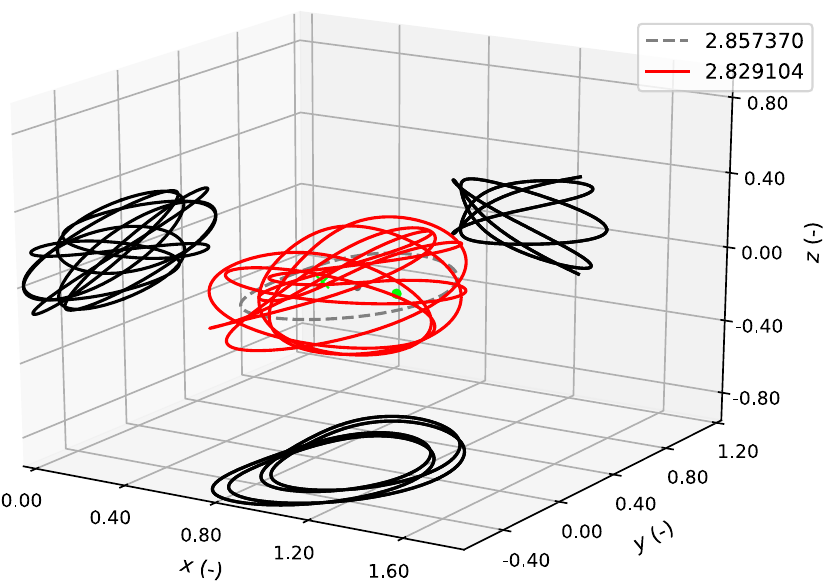} &
			\includegraphics[width=.47\linewidth]{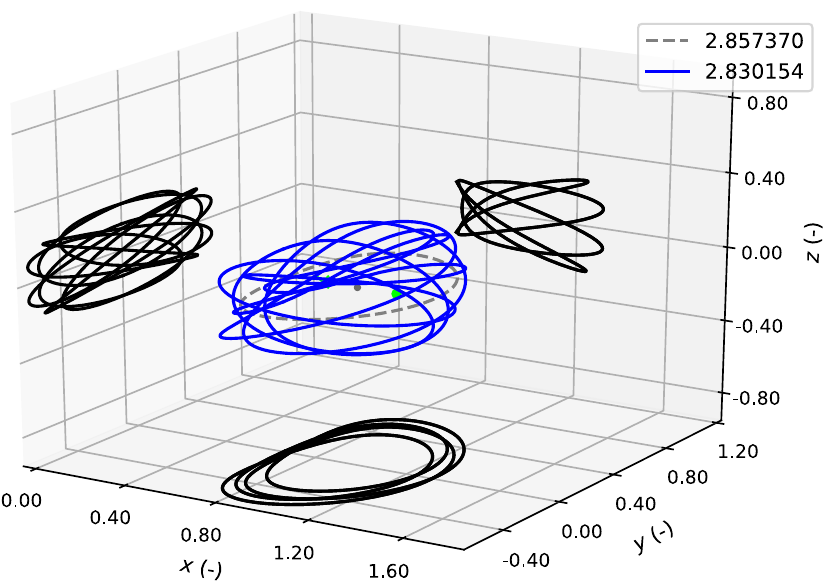} \\
			\includegraphics[width=.47\linewidth]{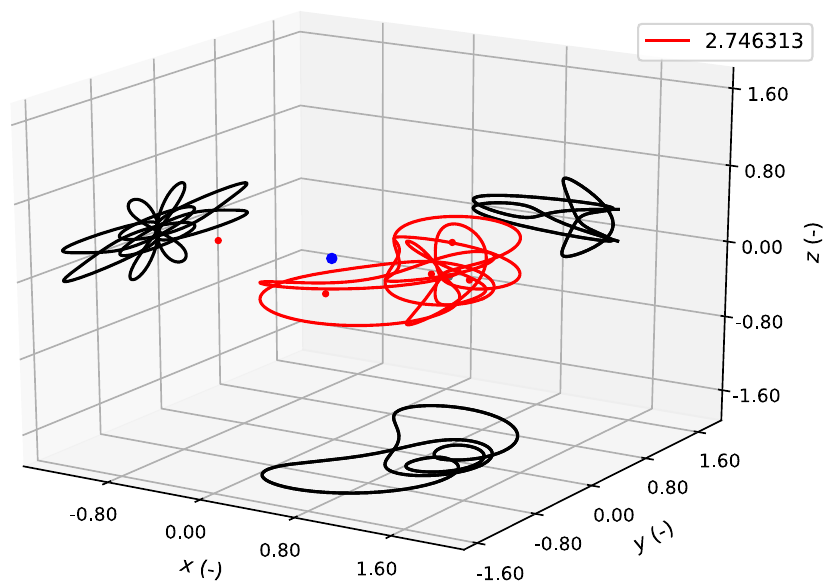} &
			\includegraphics[width=.47\linewidth]{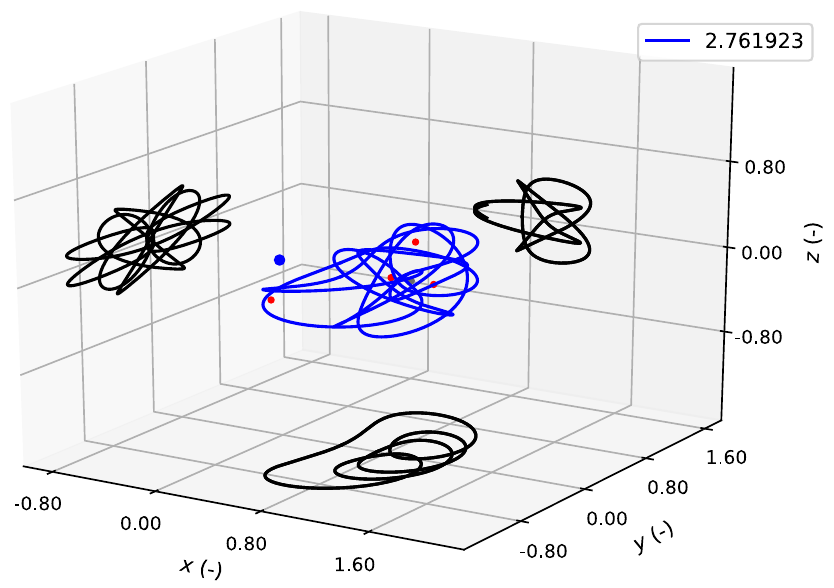} \\
			\includegraphics[width=.47\linewidth]{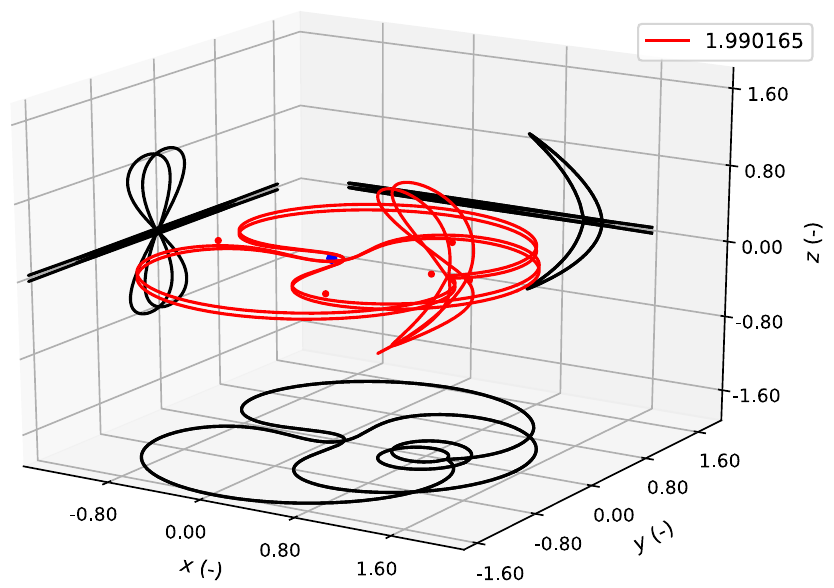} &
			\includegraphics[width=.47\linewidth]{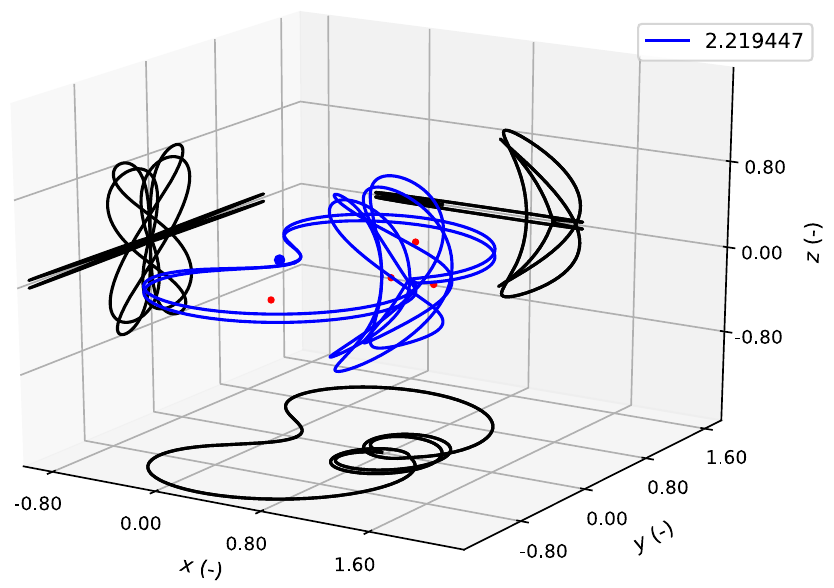}
		\end{tabular}
		\caption{VSR bifurcations from the 7:8 DRO.}
	\end{subfigure}
	
	\vspace{0.8em}
	
	\begin{subfigure}{0.48\textwidth}
		\centering
		\begin{tabular}{cc}
			\includegraphics[width=.47\linewidth]{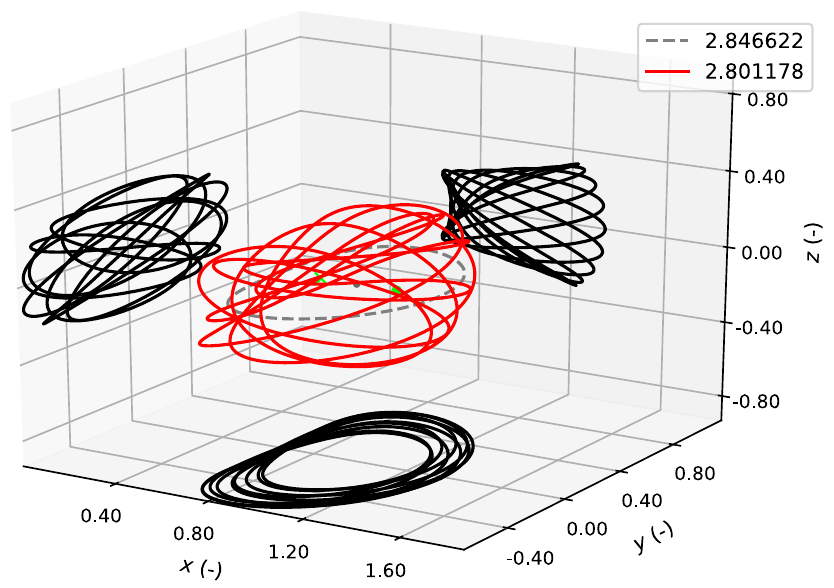} &
			\includegraphics[width=.47\linewidth]{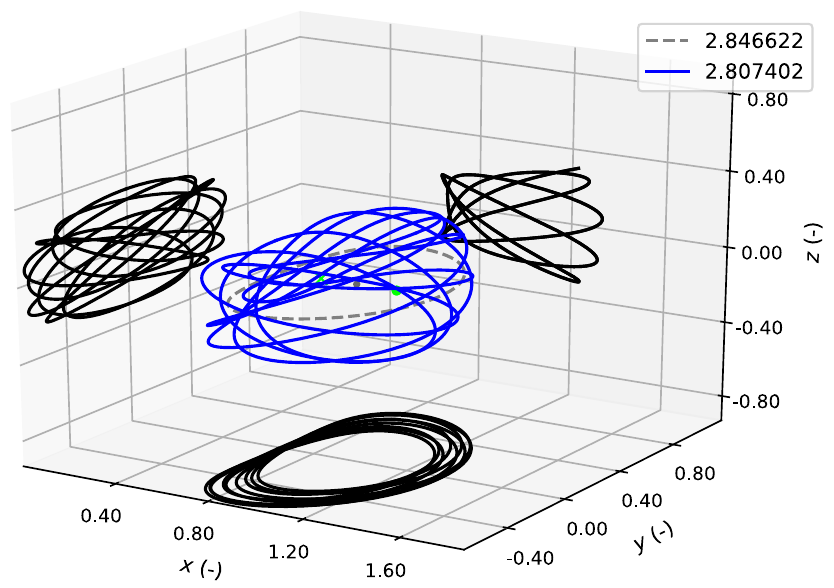} \\
			\includegraphics[width=.47\linewidth]{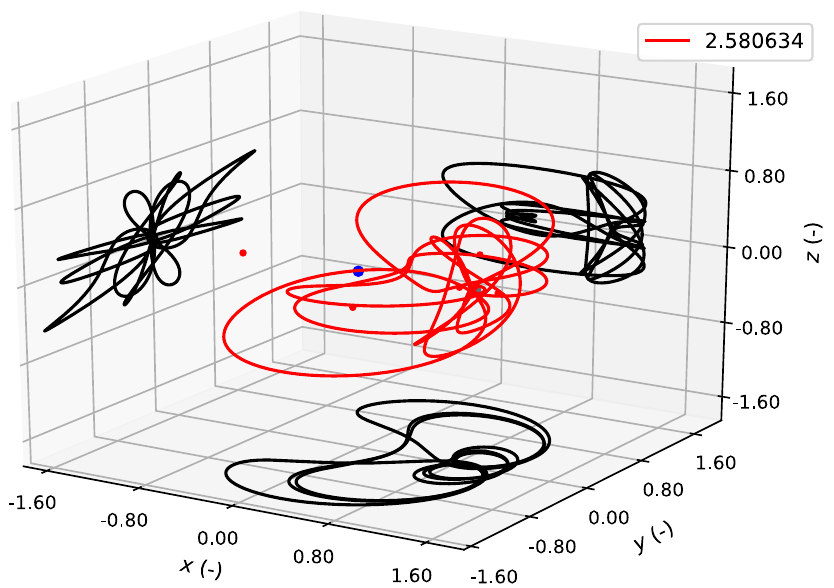} &
			\includegraphics[width=.47\linewidth]{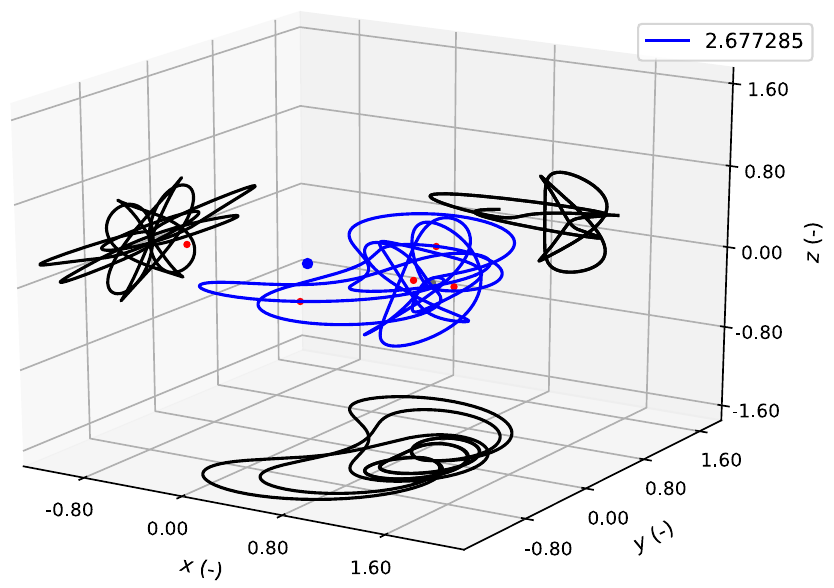} \\
			\includegraphics[width=.47\linewidth]{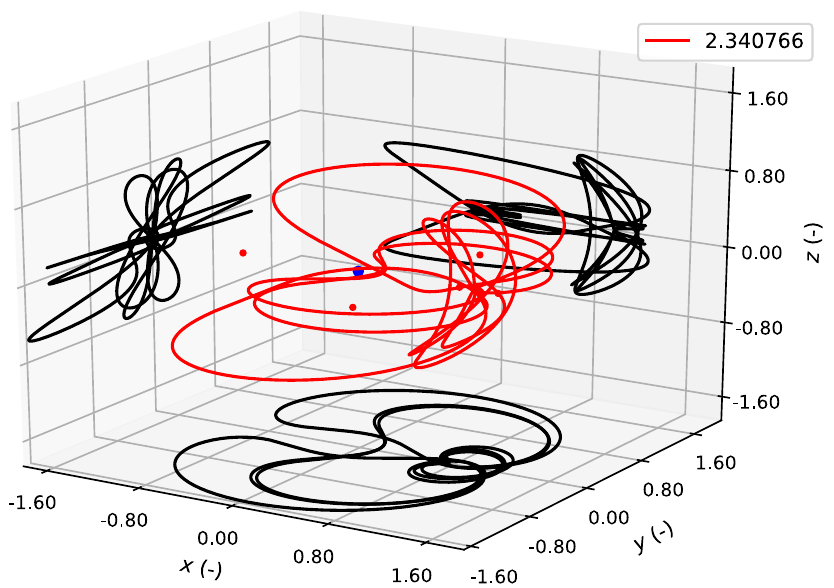} &
			\includegraphics[width=.47\linewidth]{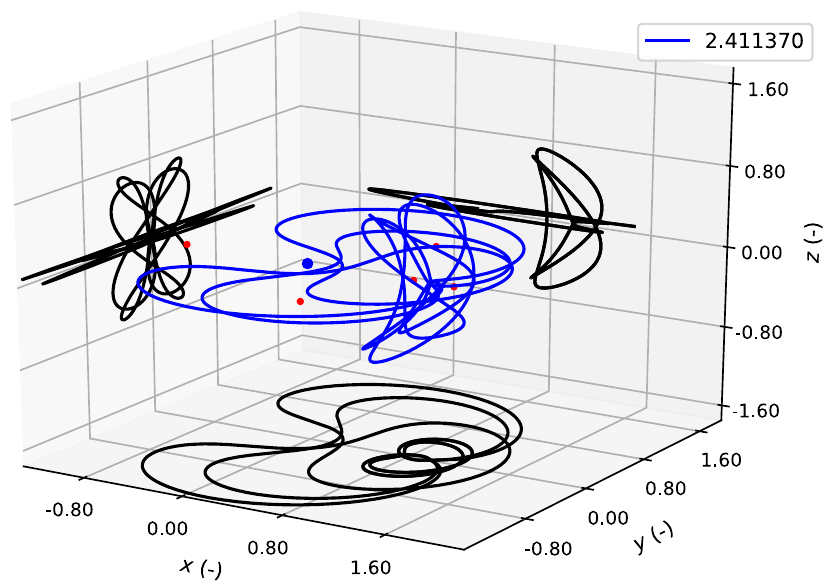}
		\end{tabular}
		\caption{VSR bifurcations from the 8:9 DRO.}
	\end{subfigure}
	\hfill
	\begin{subfigure}{0.48\textwidth}
		\centering
		\begin{tabular}{cc}
			\includegraphics[width=.47\linewidth]{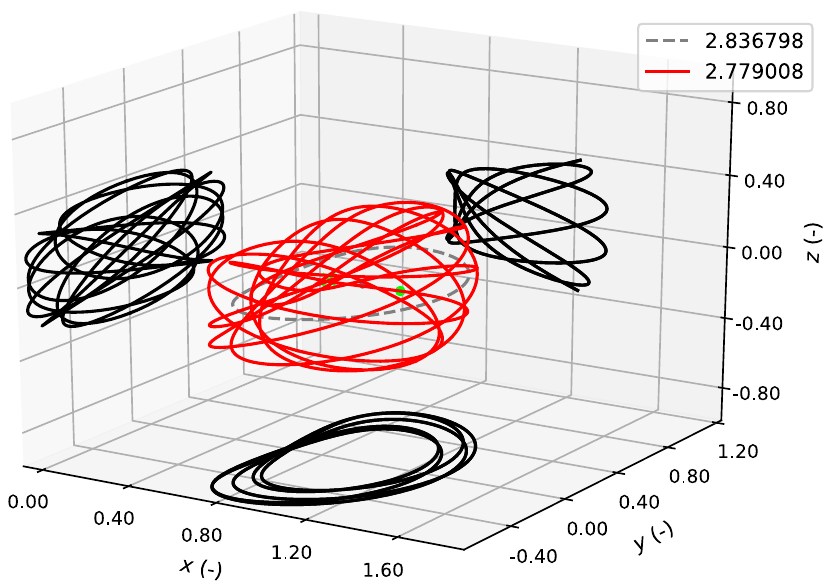} &
			\includegraphics[width=.47\linewidth]{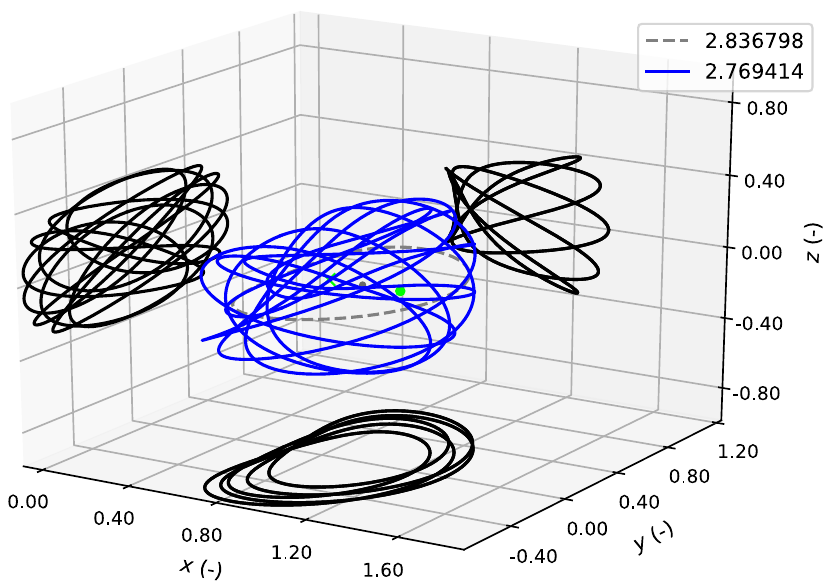} \\
			\includegraphics[width=.47\linewidth]{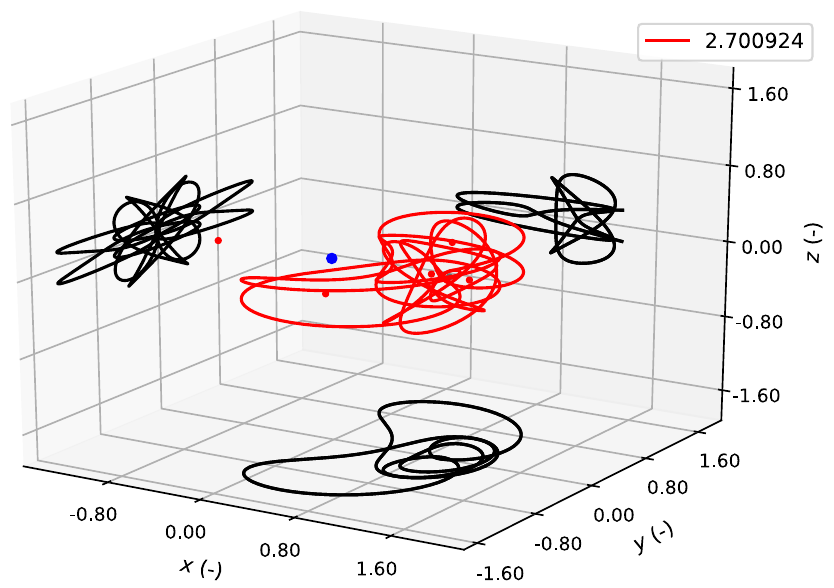} &
			\includegraphics[width=.47\linewidth]{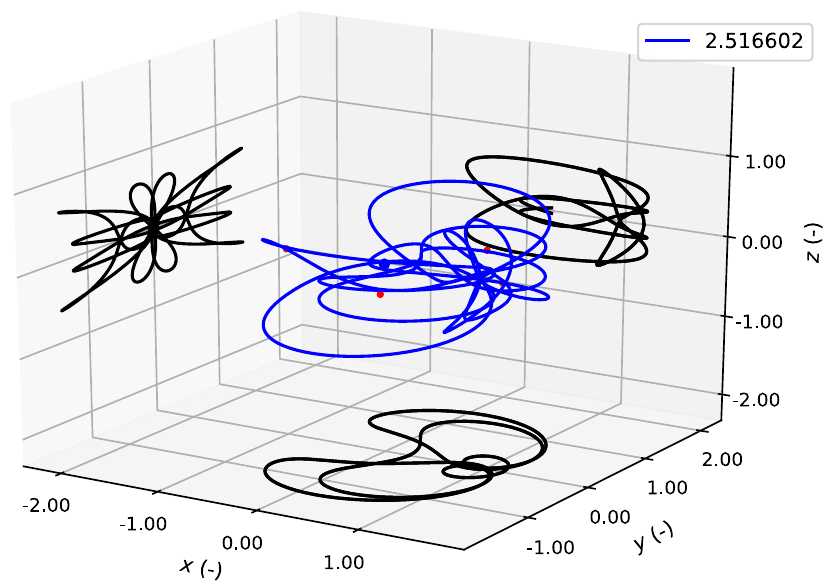} \\
			\includegraphics[width=.47\linewidth]{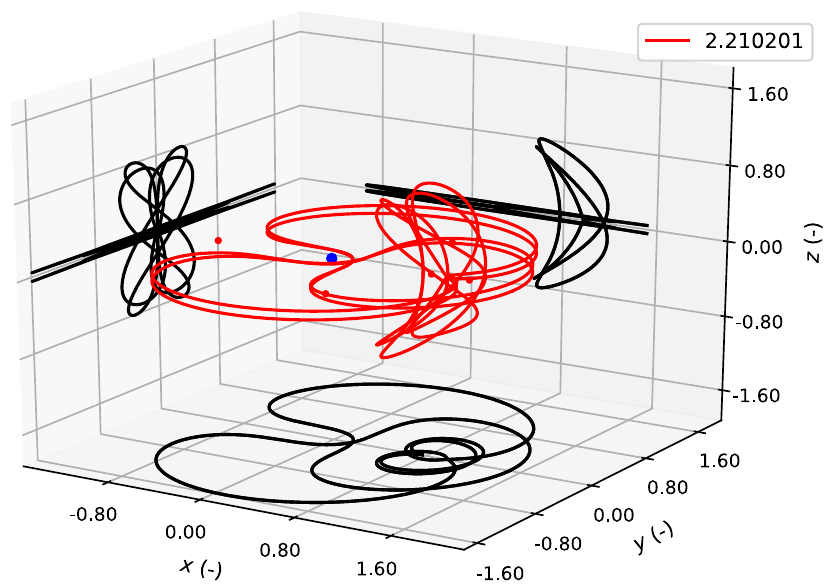} &
			\includegraphics[width=.47\linewidth]{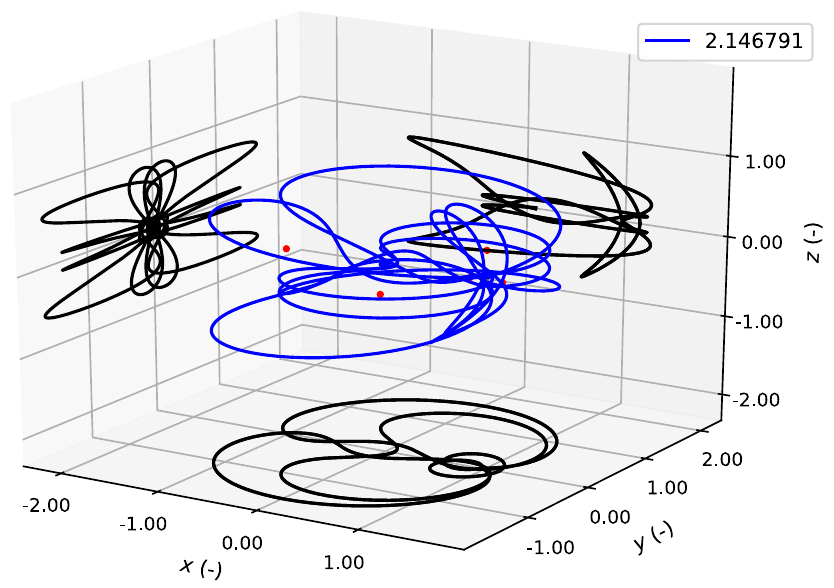}
		\end{tabular}
		\caption{VSR bifurcations from the 9:10 DRO.}
	\end{subfigure}
	
	\caption{Representative orbits of the bifurcated families from the higher-order VSR bifurcations of the DROs that move around the Moon, $L_1$ and $L_2$ $(q=7,8,9,10)$.\ In each case, the left and right column shows orbits of the red and blue family, respectively.\ The families are continued from top to bottom, starting at the corresponding DRO (gray dashed) and are continued until being of long-period multi-excursion type.\ For the odd-period resonances, the $\rho_{xz}$-symmetric and $\rho_x$-symmetric families correspond to the red and blue families in (a), and to the blue and red families in (c), respectively.\ For the even-period resonances, (b) and (d), both families consist of doubly symmetric orbits.\ Within the first row in each of (a)--(d):\ \textcolor{gray}{$\bullet$}, \textcolor{green}{$\times$} and \textcolor{green}{$\bullet$} denotes the Moon and the libration points $L_1$ and $L_2$, respectively.\ Within the remaining plots in (a)--(d):\ \textcolor{blue}{$\bullet$}, \textcolor{gray}{$\bullet$} and \textcolor{red}{$\bullet$} denotes the Earth, the Moon and the libration points, respectively.}
	\label{figure_plot_8}
\end{figure}

\begin{figure}[t!]
	\centering
	
	\begin{subfigure}{0.45\textwidth}
		\centering
		\includegraphics[width=\textwidth]{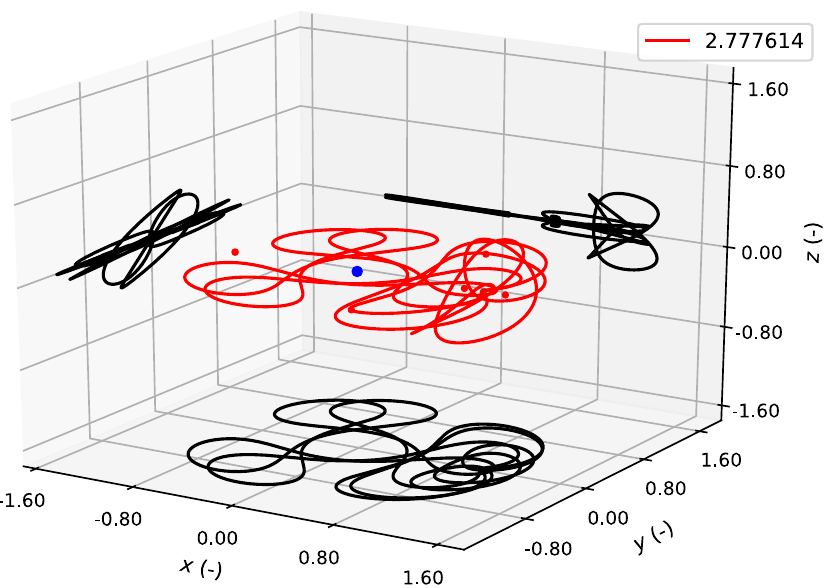}
	\end{subfigure}
	\begin{subfigure}{0.45\textwidth}
		\centering
		\includegraphics[width=\textwidth]{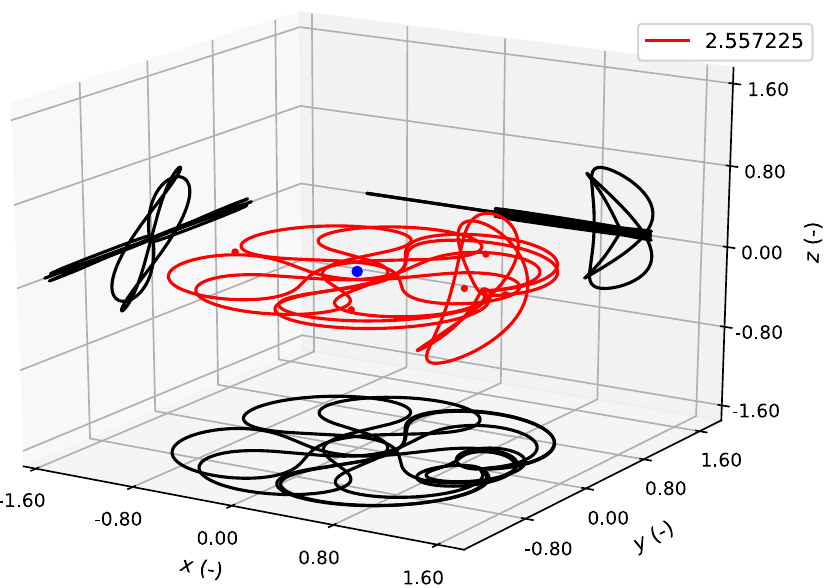}
	\end{subfigure}
	\caption{Spatial Earth--Moon cycling-type orbits that are symmetric with respect to $\rho_x$.\ Within the plots:\ \textcolor{blue}{$\bullet$}, \textcolor{gray}{$\bullet$} and \textcolor{red}{$\bullet$} denotes the Earth, the Moon and the libration points, respectively.}
	\label{figure_plot_9}
\end{figure}

\begin{table}[th!]
	\centering
	\begin{tabular}{c|cc|cc}
		\multirow{2}{*}{$(q-1)$:$q$ DRO}
		& \multicolumn{2}{c|}{Red family}
		& \multicolumn{2}{c}{Blue family} \\
		\cline{2-5}
		& $C$ & $d_{\min}$ & $C$ & $d_{\min}$ \\
		\hline
		6:7  & 2.361166 & 0.013940 & 2.274824 & 0.012765 \\
		7:8  & 1.990165 & 0.009446 & 2.219447 & 0.013174 \\
		8:9  & 2.340766 & 0.016672 & 2.411370 & 0.017469 \\
		9:10 & 2.210201 & 0.012691 & 2.146791 & 0.018325 \\
		\hline
	\end{tabular}
	\caption{Jacobi constants $C$ and minimum lunar distances $d_{\min}$ of the final computed long-period multi-excursion orbits shown in Figure~\ref{figure_plot_8}.\ Distances are given in normalized units.}
	\label{tab:min_lunar_distance}
\end{table}

Representative orbits associated with the VSR bifurcations of multiplicities $q=7,8,9,10$ are presented in Figure~\ref{figure_plot_8}.\ Although the geometric complexity of these orbit families increases along the continuation, the final computed long-period multi-excursion orbits were continued until their minimum lunar distances reached values between 0.009446 and 0.018325 normalized units (see Table~\ref{tab:min_lunar_distance}).\ Furthermore, an additional dynamical feature is achieved along the continuation of the red family generated by the 8:9 VSR bifurcation.\ As the continued family approaches the Moon, the continuation transferred to another family of long-period multi-excursion orbits whose near-ecliptic components indicate the Earth--Moon cycling geometry studied in~\cite{deleo_pontani}.\ We refer to these spatial orbits as \textit{spatial Earth--Moon cycling-type orbits}.\ Representative members of this family are shown in Figure~\ref{figure_plot_9}, and their numerical data are collected in Table~\ref{data_13}.\ Using the CZ index computational algorithm from~\cite{moreno_aydin}, these orbits are found to have CZ index 42.\ Although the continuation could not be traced back to the bifurcation that generates this family of Earth--Moon cycling-type orbits, its occurrence shows the richness of the dynamics in this region of the Earth--Moon CR3BP.\

\subsection{Vertical Pitchfork Bifurcations from large DROs around the Moon, $L_1$, $L_2$, $L_4$ and $L_5$}

\begin{figure}[t!]
	\centering
	
	\begin{subfigure}{0.25\textwidth}
		\centering
		\includegraphics[scale=0.5]{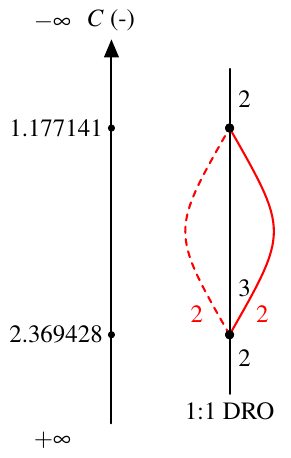}
		\caption{Bifurcation diagram.}
		\label{plot_10_a}
	\end{subfigure}
	\hfill
	\begin{subfigure}{0.35\textwidth}
		\centering
		\includegraphics[width=\textwidth]{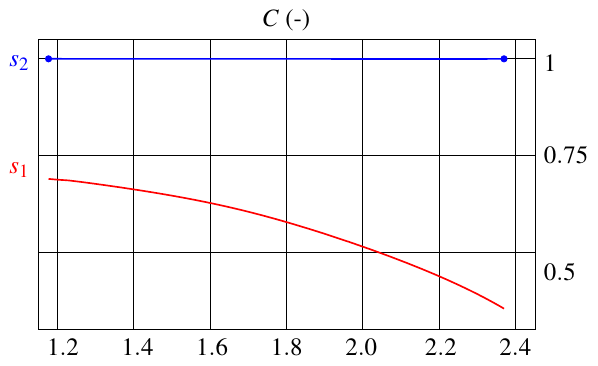}
		\caption{Stability index curves of $s_1$ and $s_2$ of the bifurcated family.}
		\label{plot_10_b}
	\end{subfigure}
	\hfill
	\begin{subfigure}{0.35\textwidth}
		\centering
		\includegraphics[width=\textwidth]{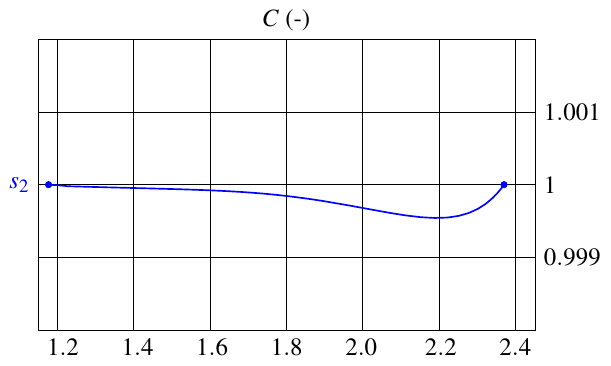}
		\caption{Zoomed view of the stability index curve $s_2$.}
		\label{plot_10_c}
	\end{subfigure}
	\vspace{0.5cm}
	
	\begin{subfigure}{\textwidth}
		\centering
		\includegraphics[width=0.32\textwidth]{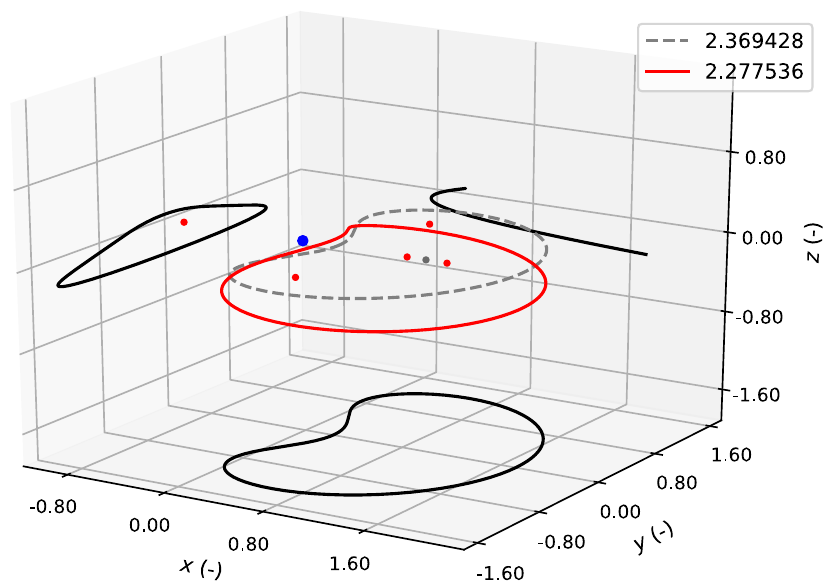}
		\hfill
		\includegraphics[width=0.32\textwidth]{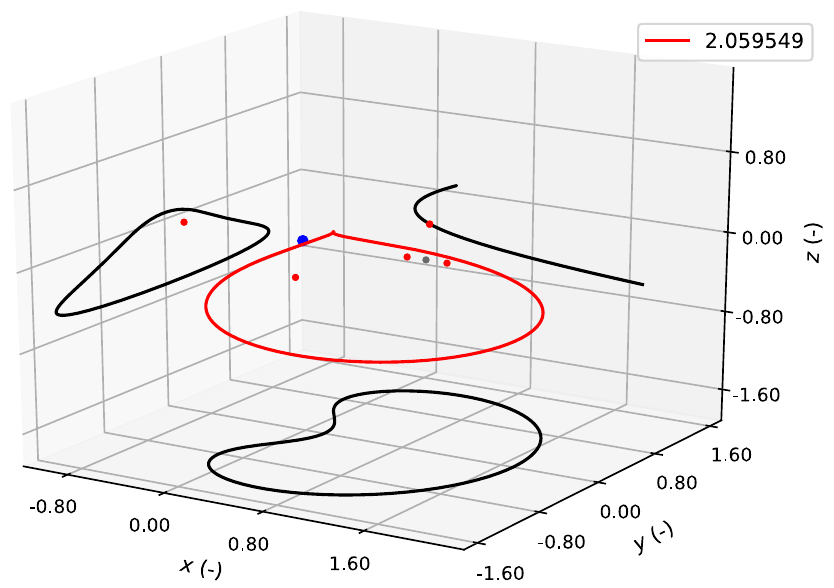}
		\hfill
		\includegraphics[width=0.32\textwidth]{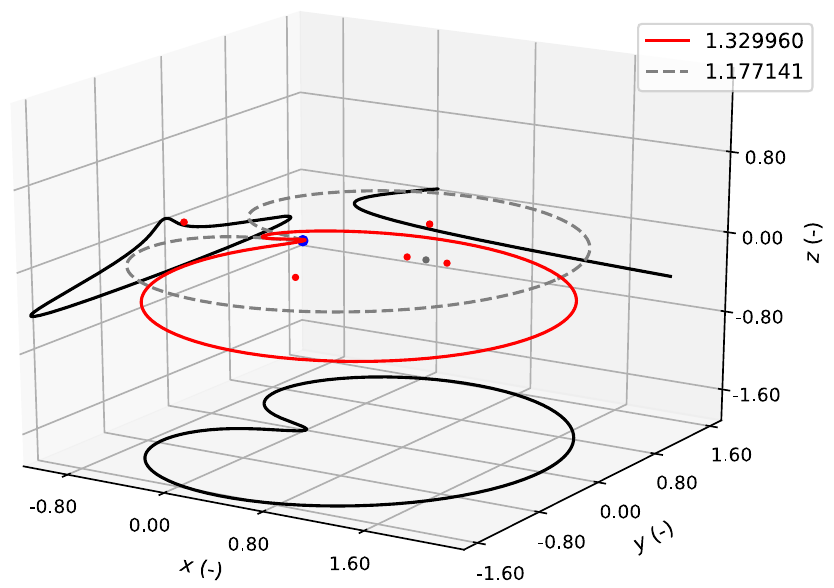}
		\caption{Representative orbits of the bifurcated family.\ The middle orbit is a 1:1 Earth--Moon resonant orbit.\ \textcolor{blue}{$\bullet$}, \textcolor{gray}{$\bullet$} and \textcolor{red}{$\bullet$} denotes the Earth, the Moon and the libration points, respectively.}
		\label{plot_10_d}
	\end{subfigure}
	\vspace{0.1em}
	
	\begin{subfigure}{\textwidth}
		\centering
		\includegraphics[width=0.45\textwidth]{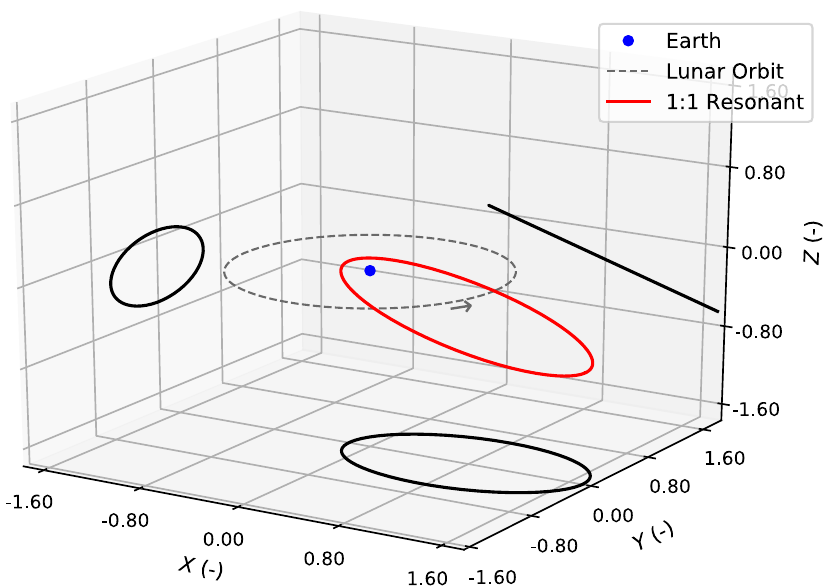}
		\includegraphics[width=0.45\textwidth]{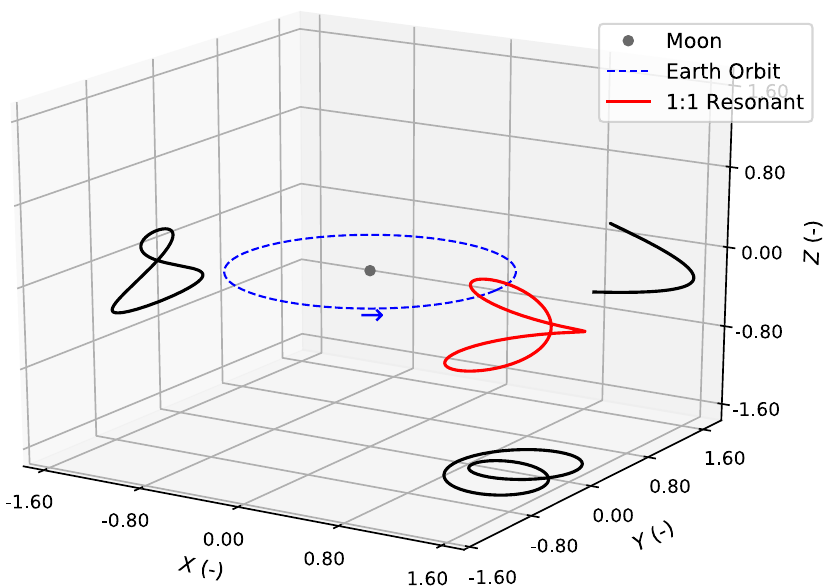}
		\caption{The 1:1 Earth--Moon resonant orbit from (d) is shown in the Earth-centered inertial frame (left) and in the Moon-centered inertial frame (right).}
		\label{plot_10_e}
	\end{subfigure}
	\caption{Bridge family connecting the two vertical pitchfork bifurcations of the large 1:1 DROs.}
\end{figure}

We finally consider the two vertical single-turn bifurcations that occur within the third geometric regime, where the DROs move around the Moon together with the collinear and triangular libration points $L_1$, $L_2$, $L_4$ and $L_5$.\ These bifurcations correspond to vertical pitchfork bifurcations that arise from the monotonic crossing of the vertical stability index through $s_v = 1$ (see Figure~\ref{plot_1_b}).\ The main result is that the two vertical pitchfork bifurcations are connected by a single family of spatial periodic orbits, thereby forming an internal bridge between these two large vertical-critical DROs, as illustrated in the bifurcation diagram in Figure~\ref{plot_10_a}.\ As for the halo orbits, the bifurcated orbit families occur in two $\sigma$-symmetric families, which we also classify into northern and southern families according to whether the orbits extend above or below the ecliptic plane $\{z=0\}$.\ Representative orbits of the southern family are shown in Figure~\ref{plot_10_d}.\ Based on the database provided in Table~\ref{data_15}, the stability indices are shown in Figures~\ref{plot_10_b} and~\ref{plot_10_c}.\ In particular, the stability indices indicate that all spatial periodic orbits are (linearly) stable with constant CZ index 2.\ Furthermore, along the continuation, we identify a 1:1 Earth--Moon resonant orbit that is shown in the rotating frame in the middle of Figure~\ref{plot_10_d}.\ This resonant orbit is plotted in the Earth- and Moon-centered inertial frames in Figure~\ref{plot_10_e}.\ Similar spatial resonant orbits have previously been reported in~\cite{vaquero,oshima}, while the present computations complement them within the global bifurcation structure of the large DROs.\

\subsection{Summary}

\begin{figure}[t!]
	\centering
	\includegraphics[width=1\linewidth]{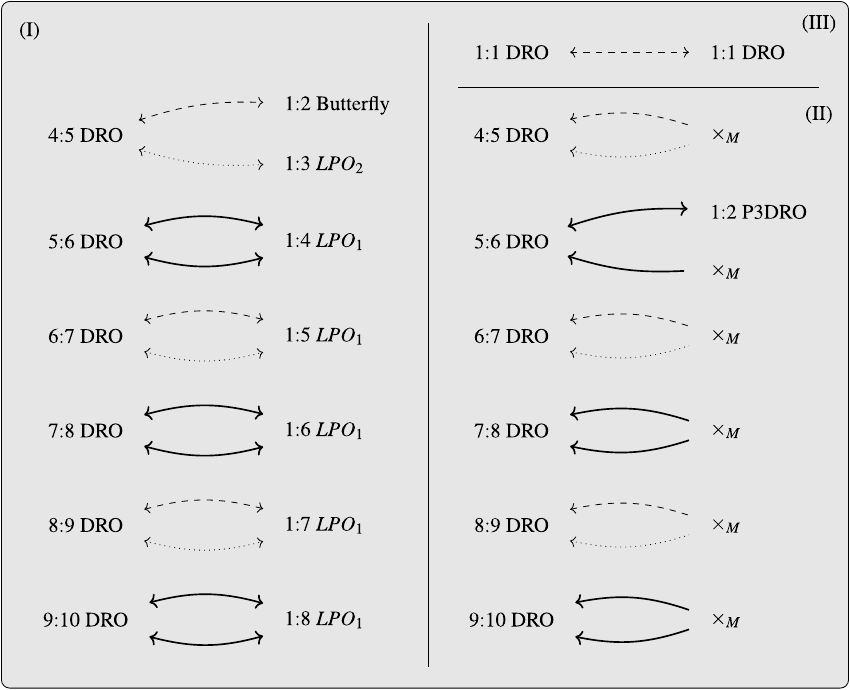}
	\caption{Summary of the bifurcation structures generated by the fourteen vertical self-resonant DROs studied in the present work.\ (I) Bridge families connecting the near-Moon DROs with the Butterfly and planar prograde families $LPO_2$ and $LPO_1$.\ (II) Higher-order VSR bifurcations of the DROs that move around the Moon and the libration points $L_1$ and $L_2$, where only the 5:6 self-resonance connects to the double-covered P3DRO family, while the remaining families evolve toward orbits of long-period multi-excursion type approaching lunar collision (denoted by ``$\times_M$'').\ (III) Bridge of spatial stable orbits connecting the two vertical pitchfork bifurcations of the large DROs.\ Thick arrows indicate families with doubly symmetric orbits, while dashed (or dotted) arrows represent $\rho_{xz}$- (or $\rho_x$-)symmetric orbits.}
	\label{figure_summary}
\end{figure}

The main bifurcation results obtained in the present work are summarized in Figure~\ref{figure_summary}.\ The identified fourteen vertical-critical DROs generate three qualitatively distinct continuation scenarios.\ The first consists of bridges connecting the near-Moon DROs with the Butterfly and basic planar prograde families.\ The bridges connecting the DRO and planar prograde families exhibit a systematic pattern:\ the $q$-covered DRO family is connected to the $(q-2)$-covered planar prograde family, for $q=5,...,10$.\ This relation is consistently reflected by the associated CZ indices.\ The second scenario considers DROs that move around the Moon, $L_1$ and $L_2$, which is characterized by predominantly open continuation families, where only the 6-covered DROs form a bridge to the double-covered P3DRO family, while the remaining families evolve toward highly unstable long-period multi-excursion orbits approaching collision with the Moon.\ The third scenario consists of a stable bridge connecting two vertical pitchfork bifurcations of large DROs that move around the Moon, $L_1$, $L_2$, $L_4$ and $L_5$.\ Together, these three scenarios provide a classification of the bifurcation topology generated by the Earth--Moon DRO family.\

\section{Conclusion}
\label{sec5}

In this work, we developed a computational symplectic framework for the systematic investigation of vertical self-resonant bifurcations of symmetric periodic orbit families in the Circular Restricted Three-Body Problem (CR3BP).\ By combining numerical continuation with the analysis of local bifurcations through Floquet multiplier transitions, Conley--Zehnder indices, Krein signatures and the Euler characteristic of the local Floer homology, we established a methodology for characterizing and interpreting the global bifurcation topology.\ The local symplectic classification underlying this methodology is summarized in Figure~\ref{fig:bifurcations}.\ The framework was applied to the higher-order vertical self-resonant bifurcations of the Earth--Moon Distant Retrograde Orbit (DRO) family.\ The main outcomes are the bifurcation diagrams shown in Figures~\ref{figure_bif_1}, \ref{figure_bif_2}, \ref{plot_6_b}, and~\ref{figure_bif_4}.\ These diagrams go beyond graphical representations of numerical continuations by revealing the global bifurcation network connecting different periodic orbit families.\ Numerical continuation establishes these connections, while symplectic invariants provide complementary tools for their classification, interpretation, and topological consistency checks.\ A synthetic overview of the resulting bifurcation network is presented in Figure~\ref{figure_summary}.\

The present study has two principal limitations.\ First, it is restricted to symmetric periodic orbit families and to the vertical self-resonant bifurcations generated from symmetric planar periodic orbits.\ Therefore, asymmetric periodic orbits, their bifurcations and more general continuation procedures have not been considered.\ Second, all computations have been performed in the physical coordinates of the CR3BP and therefore do not incorporate regularization techniques.\ Such techniques could be applied to continue the highly unstable long-period multi-excursion orbits beyond their close approaches to the Moon.\ Moreover, regularization allows the DRO family itself to be continued through collision with the Earth \cite{broucke}, thereby providing access to additional vertical self-resonant bifurcations beyond those considered in the present work.\ Nevertheless, the computational symplectic framework and the periodic orbit families established in the present study provide a natural foundation for several future research directions.\ One immediate extension is the investigation of asymmetric periodic orbit families.\ A second direction is the application and extension of the proposed methodology to higher-fidelity models that arise through time-dependent perturbations of the CR3BP.\ Representative examples include the bicircular restricted four-body problem and the elliptic restricted three-body problem, where resonant periodic orbits of the CR3BP provide natural starting points for continuation into periodic orbits of the perturbed dynamics.\

\clearpage
\appendix

\section{Predictor--Corrector Equations for Spatial Symmetric Periodic Orbits}
\label{appendix_s}

This appendix summarizes the predictor--corrector equations discussed in Section \ref{sec3} for the four types of spatial symmetric periodic orbits.\ Tables~\ref{tab:spatial_summary_1} and~\ref{tab:spatial_summary_2} present the simple and doubly symmetric cases, respectively.\ For each symmetry type, the initial and final state together with the corresponding predictor--corrector equations are listed.\ Among the equations, the ones whose denominator has the largest absolute value determines the continuation parameter.\ The corrected variational coefficients are defined by
$$ u_{3j} = v_{3j} - v_{2j} \frac{p_{z,cut}}{p_{y,cut} - x_{cut}},\quad u_{4j} = v_{4j} - v_{2j} \frac{\dot{p}_{x,cut}}{p_{y,cut} - x_{cut}},\quad u_{6j} = v_{6j} - v_{2j} \frac{\dot{p}_{z,cut}}{p_{y,cut} - x_{cut}},\quad j=1,3,5,6. $$

\begin{table}[ht!]
	\centering
	\renewcommand{\arraystretch}{1.3}
	\begin{tabular}{lll}
		\hline	
		\textbf{Symmetry type:} $\rho_x$ & \textbf{Initial state:} $(x_0,0,0,0,p_{y0},p_{z0})$ & \textbf{Final state:} $(x_{cut},0,0,0,p_{y,cut},p_{z,cut})$ \\
		\hline
		\textbf{Fixed parameter / Step size} & \textbf{Corrections / Predictors} & \\
		\hline
		& & \vspace{-1em} \\
		$p_{z0}$ & $ \delta x_0 = \dfrac{u_{35} p_{x,cut} - u_{45} z_{cut}}{u_{31}u_{45} - u_{35}u_{41}}$ & $\delta p_{y0} = \dfrac{u_{41} z_{cut} - u_{31} p_{x,cut}}{u_{31}u_{45} - u_{35}u_{41}} $ \vspace{-0.75em} \\
		& &\\
		$\Delta p_{z0}$ &  $ \Delta x_0 = - \dfrac{u_{35} u_{46} - u_{45} u_{36}}{u_{31}u_{45} - u_{35}u_{41}} \Delta p_{z0}$ & $\Delta p_{y0} = - \dfrac{u_{41} u_{36} - u_{31} u_{46}}{u_{31}u_{45} - u_{35}u_{41}} \Delta p_{z0} $ \vspace{0.75em}\\
		\hdashline
		& & \vspace{-0.75em}\\
		$p_{y0}$ & $\delta x_0 = \dfrac{u_{36} p_{x,cut} - u_{46} z_{cut}}{u_{31}u_{46} - u_{36}u_{41}}$ & $\delta p_{z0} = \dfrac{u_{41} z_{cut} - u_{31} p_{x,cut}}{u_{31}u_{46} - u_{36}u_{41}}$ \vspace{-0.75em}\\
		& &\\
		$\Delta p_{y0}$ & $\Delta x_0 = - \dfrac{u_{36} u_{45} - u_{46} u_{35}}{u_{31}u_{46} - u_{36}u_{41}} \Delta p_{y0}$ & $\Delta p_{z0} = - \dfrac{u_{41} u_{35} - u_{31} u_{45}}{u_{31}u_{46} - u_{36}u_{41}} \Delta p_{y0}$ \vspace{0.75em}\\
		\hdashline
		& & \vspace{-0.75em}\\
		$x_0$ & $\delta p_{y0} = \dfrac{u_{36} p_{x,cut} - u_{46} z_{cut}}{u_{35}u_{46} - u_{36}u_{45}}$ & $\delta p_{z0} = \dfrac{u_{45} z_{cut} - u_{35} p_{x,cut}}{u_{35}u_{46} - u_{36}u_{45}}$ \vspace{-0.75em}\\
		& &\\
		$\Delta x_0$ & $\Delta p_{y0} = - \dfrac{u_{36} u_{41} - u_{46} u_{31}}{u_{35}u_{46} - u_{36}u_{45}} \Delta x_0$ & $\Delta p_{z0} = - \dfrac{u_{45} u_{31} - u_{35} u_{41}}{u_{35}u_{46} - u_{36}u_{45}} \Delta x_0$ \vspace{0.75em}\\
		\hline	
		\textbf{Symmetry type:} $\rho_{xz}$ & \textbf{Initial state:} $(x_0,0,z_0,0,p_{y0},0)$ & \textbf{Final state:} $(x_{cut},0,z_{cut},0,p_{y,cut},0)$\\
		\hline
		\textbf{Fixed parameter / Step size} & \textbf{Corrections / Predictors} & \\
		\hline
		& & \vspace{-1em} \\
		$z_0$ & $ \delta x_0 = \dfrac{u_{45} p_{z,cut} - u_{65} p_{x,cut}}{u_{41}u_{65} - u_{61}u_{45}}$ & $\delta p_{y0} = \dfrac{u_{61} p_{x,cut} - u_{41} p_{z,cut}}{u_{41}u_{65} - u_{61}u_{45}} $ \vspace{-0.75em} \\
		& &\\
		$\Delta z_0$ &  $ \Delta x_0 = - \dfrac{u_{45} u_{63} - u_{65} u_{43}}{u_{41}u_{65} - u_{61}u_{45}} \Delta p_{z0}$ & $\Delta p_{y0} = - \dfrac{u_{61} u_{43} - u_{41} u_{63}}{u_{41}u_{65} - u_{61}u_{45}} \Delta p_{z0} $ \vspace{0.75em}\\
		\hdashline
		& & \vspace{-0.75em}\\
		$p_{y0}$ & $\delta x_0 = \dfrac{u_{43} p_{z,cut} - u_{63} p_{x,cut}}{u_{41}u_{63} - u_{61}u_{43}}$ & $\delta z_0 = \dfrac{u_{61} p_{x,cut} - u_{41} p_{z,cut}}{u_{41}u_{63} - u_{61}u_{43}}$ \vspace{-0.75em} \\
		& &\\
		$\Delta p_{y0}$ & $\Delta x_0 = - \dfrac{u_{43} u_{65} - u_{63} u_{45}}{u_{41}u_{63} - u_{61}u_{43}} \Delta p_{y0}$ & $\Delta z_0 = - \dfrac{u_{61} u_{45} - u_{41} u_{65}}{u_{41}u_{63} - u_{61}u_{43}} \Delta p_{y0}$ \vspace{0.75em}\\
		\hdashline
		& & \vspace{-0.75em}\\
		$x_0$ & $\delta z_0 = \dfrac{u_{45} p_{z,cut} - u_{65} p_{x,cut}}{u_{43}u_{65} - u_{63}u_{45}}$ & $\delta p_{y0} = \dfrac{u_{63} p_{x,cut} - u_{43} p_{z,cut}}{u_{43}u_{65} - u_{63}u_{45}}$ \vspace{-0.75em} \\
		& &\\
		$\Delta x_0$ & $\Delta z_0 = - \dfrac{u_{45} u_{61} - u_{65} u_{41}}{u_{43}u_{65} - u_{63}u_{45}} \Delta x_0$ & $\Delta p_{y0} = - \dfrac{u_{63} u_{41} - u_{43} u_{61}}{u_{43}u_{65} - u_{63}u_{45}} \Delta x_0$
	\end{tabular}
	\caption{Summary of the predictor--corrector equations for spatial simple symmetric periodic orbits.}
	\label{tab:spatial_summary_1}
\end{table}

\begin{table}[ht!]
	\centering
	\renewcommand{\arraystretch}{1.3}
	\begin{tabular}{lll}
		\hline	
		\textbf{Symmetry type:} $\rho_x$-$\rho_{xz}$ & \textbf{Initial state:} $(x_0,0,0,0,p_{y0},p_{z0})$ & \textbf{Final state:} $(x_{cut},0,z_{cut},0,p_{y,cut},0)$\\
		\hline
		\textbf{Fixed parameter / Step size} & \textbf{Corrections / Predictors} & \\
		\hline
		& & \vspace{-1em} \\
		$p_{z0}$ & $ \delta x_0 = \dfrac{u_{45} p_{z,cut} - u_{65} p_{x,cut}}{u_{41}u_{65} - u_{61}u_{45}}$ & $\delta p_{y0} = \dfrac{u_{61} p_{x,cut} - u_{41} p_{z,cut}}{u_{41}u_{65} - u_{61}u_{45}} $ \vspace{-0.75em} \\
		& &\\
		$\Delta p_{z0}$ &  $ \Delta x_0 = - \dfrac{u_{45} u_{66} - u_{65} u_{46}}{u_{41}u_{65} - u_{61}u_{45}} \Delta p_{z0}$ & $\Delta p_{y0} = - \dfrac{u_{61} u_{46} - u_{41} u_{66}}{u_{41}u_{65} - u_{61}u_{45}} \Delta p_{z0} $ \vspace{0.75em}\\
		\hdashline
		& & \vspace{-0.75em}\\
		$p_{y0}$ & $\delta x_0 = \dfrac{u_{46} p_{z,cut} - u_{66} p_{x,cut}}{u_{41}u_{66} - u_{61}u_{46}}$ & $\delta p_{z0} = \dfrac{u_{61} p_{x,cut} - u_{41} p_{z,cut}}{u_{41}u_{66} - u_{61}u_{46}}$ \vspace{-0.75em}\\
		& &\\
		$\Delta p_{y0}$ & $\Delta x_0 = - \dfrac{u_{46} u_{65} - u_{66} u_{45}}{u_{41}u_{66} - u_{61}u_{46}} \Delta p_{y0}$ & $\Delta p_{z0} = - \dfrac{u_{61} u_{45} - u_{41} u_{65}}{u_{41}u_{66} - u_{61}u_{46}} \Delta p_{y0}$ \vspace{0.75em}\\
		\hdashline
		& & \vspace{-0.75em}\\
		$x_0$ & $\delta p_{y0} = \dfrac{u_{46} p_{z,cut} - u_{66} p_{x,cut}}{u_{45}u_{66} - u_{65}u_{46}}$ & $\delta p_{z0} = \dfrac{u_{65} p_{x,cut} - u_{45} p_{z,cut}}{u_{45}u_{66} - u_{65}u_{46}}$ \vspace{-0.75em}\\
		& &\\
		$\Delta x_0$ & $\Delta p_{y0} = - \dfrac{u_{46} u_{61} - u_{66} u_{41}}{u_{45}u_{66} - u_{65}u_{46}} \Delta x_0$ & $\Delta p_{z0} = - \dfrac{u_{65} u_{41} - u_{45} u_{61}}{u_{45}u_{66} - u_{65}u_{46}} \Delta x_0$ \vspace{0.75em}\\
		\hline	
		\textbf{Symmetry type:} $\rho_{xz}$-$\rho_x$ & \textbf{Initial state:} $(x_0,0,z_0,0,p_{y0},0)$ & \textbf{Final state:} $(x_{cut},0,0,0,p_{y,cut},p_{z,cut})$\\
		\hline
		\textbf{Fixed parameter / Step size} & \textbf{Corrections / Predictors} & \\
		\hline
		& & \vspace{-1em} \\
		$z_0$ & $ \delta x_0 = \dfrac{u_{35} p_{x,cut} - u_{45} z_{cut}}{u_{31}u_{45} - u_{41}u_{35}}$ & $\delta p_{y0} = \dfrac{u_{41} z_{cut} - u_{31} p_{x,cut}}{u_{31}u_{45} - u_{41}u_{35}} $ \vspace{-0.75em} \\
		& &\\
		$\Delta z_0$ &  $ \Delta x_0 = - \dfrac{u_{35} u_{43} - u_{45} u_{33}}{u_{31}u_{45} - u_{41}u_{35}} \Delta p_{z0}$ & $\Delta p_{y0} = - \dfrac{u_{41} u_{33} - u_{31} u_{43}}{u_{31}u_{45} - u_{41}u_{35}} \Delta p_{z0} $ \vspace{0.75em}\\
		\hdashline
		& & \vspace{-0.75em}\\
		$p_{y0}$ & $\delta x_0 = \dfrac{u_{33} p_{x,cut} - u_{43} z_{cut}}{u_{31}u_{43} - u_{41}u_{33}}$ & $\delta z_0 = \dfrac{u_{41} z_{cut} - u_{31} p_{x,cut}}{u_{31}u_{43} - u_{41}u_{33}}$ \vspace{-0.75em} \\
		& &\\
		$\Delta p_{y0}$ & $\Delta x_0 = - \dfrac{u_{33} u_{45} - u_{43} u_{35}}{u_{31}u_{43} - u_{41}u_{33}} \Delta p_{y0}$ & $\Delta z_0 = - \dfrac{u_{41} u_{35} - u_{31} u_{45}}{u_{31}u_{43} - u_{41}u_{33}} \Delta p_{y0}$ \vspace{0.75em}\\
		\hdashline
		& & \vspace{-0.75em}\\
		$x_0$ & $\delta z_0 = \dfrac{u_{35} p_{x,cut} - u_{45} z_{cut}}{u_{33}u_{45} - u_{43}u_{35}}$ & $\delta p_{y0} = \dfrac{u_{43} z_{cut} - u_{33} p_{x,cut}}{u_{33}u_{45} - u_{43}u_{35}}$ \vspace{-0.75em} \\
		& &\\
		$\Delta x_0$ & $\Delta z_0 = - \dfrac{u_{35} u_{41} - u_{45} u_{31}}{u_{33}u_{45} - u_{43}u_{35}} \Delta x_0$ & $\Delta p_{y0} = - \dfrac{u_{43} u_{31} - u_{33} u_{41}}{u_{33}u_{45} - u_{43}u_{35}} \Delta x_0$
	\end{tabular}
	\caption{Summary of the predictor--corrector equations for spatial doubly symmetric periodic orbits.}
	\label{tab:spatial_summary_2}
\end{table}

\newpage

\section{Tables of Numerical Data}
\label{appendix}

The numerical database contains Jacobi constants, initial data, periods, Floquet multipliers and Conley--Zehnder indices of selected orbits associated to the symmetric orbit families we have studied.\ For a pair of hyperbolic Floquet multipliers we denote by $\lambda$ the one with $|\lambda| > 1$.\ For a pair of elliptic Floquet multipliers we denote by~$\theta$ the rotation angle modulo $2 \pi$ characterized by Krein signature in round brackets.\ In the complex instability case, we note the two Floquet multipliers that are outside the unit circle.\ Horizontal dashed lines indicate birth--death bifurcations.\

\begin{table}[ht]\fontsize{10}{10}\selectfont \centering
	\caption{Data for planar symmetric orbits: DRO (first block), P3DRO (second block) and $L_2$ planar Lyapunov orbits (third block).}
	\begin{tabular}{c|c|c|c|c|c}
		$C$ & $x(0)$ & $p_y(0)$ & $T/2$ & Floquet multipliers & $\mu_{CZ}$ / $\mu_{CZ}^p$ / $\mu_{CZ}^s$\\
		\hline
		3.115178 & 1.044501300 & 0.520438208 & 0.346550 & $(-)$ $\theta_p = 5.572$, $(-)$ $\theta_s \approx 2 \pi$ $9/10$ & 2 / 1 / 1\\
		3.094976 & 1.049724451 & 0.539462235 & 0.390865 & $(-)$ $\theta_p = 5.483$, $(-)$ $\theta_s \approx 2 \pi$ $8/9$ & 2 / 1 / 1\\
		3.073375 & 1.056389717 & 0.560205596 & 0.448793 & $(-)$ $\theta_p = 5.369$, $(-)$ $\theta_s \approx 2 \pi$ $7/8$ & 2 / 1 / 1\\		
		3.049858 & 1.065314651 & 0.583155461 & 0.528417 & $(-)$ $\theta_p = 5.216$, $(-)$ $\theta_s \approx 2 \pi$ $6/7$ & 2 / 1 / 1\\		
		3.023442 & 1.078210256 & 0.609095797 & 0.646741 & $(-)$ $\theta_p = 5.000$, $(-)$ $\theta_s \approx 2 \pi$ $5/6$ & 2 / 1 / 1\\
		2.991353 & 1.100023113 & 0.639781288 & 0.852969 & $(-)$ $\theta_p = 4.659$, $(-)$ $\theta_s \approx 2 \pi$ $4/5$ & 2 / 1 / 1\\		
		2.954927 & 1.138211689 & 0.669474505 & 1.220987 & $(-)$ $\theta_p \approx 2 \pi$ $2/3$, $(-)$ $\theta_s = 4.819$ & 2 / 1 / 1\\
		2.902460 & 1.227777496 & 0.680344137 & 2.011714 & $(-)$ $\theta_p = 3.860$, $(-)$ $\theta_s \approx 2 \pi$ $4/5$ & 2 / 1 / 1\\
		2.883609 & 1.263819475 & 0.672130625 & 2.263090 & $(-)$ $\theta_p = 3.950$, $(-)$ $\theta_s \approx 2 \pi$ $5/6$ & 2 / 1 / 1\\
		2.869425 & 1.289516882 & 0.663559669 & 2.413348 & $(-)$ $\theta_p = 4.047$, $(-)$ $\theta_s \approx 2 \pi$ $6/7$ & 2 / 1 / 1\\		
		2.857370 & 1.310084964 & 0.655415149 & 2.516817 & $(-)$ $\theta_p = 4.132$, $(-)$ $\theta_s \approx 2 \pi$ $7/8$ & 2 / 1 / 1\\		
		2.849071 & 1.323558296 & 0.649558155 & 2.577028 & $(-)$ $\theta_p \approx 2 \pi$ $2/3$, $(-)$ $\theta_s = 5.566$ & 2 / 1 / 1\\
		2.846622 & 1.327431094 & 0.647806037 & 2.593281 & $(-)$ $\theta_p = 4.204$, $(-)$ $\theta_s \approx 2 \pi$ $8/9$ & 2 / 1 / 1\\
		2.836798 & 1.342511949 & 0.640716906 & 2.652398 & $(-)$ $\theta_p = 4.266$, $(-)$ $\theta_s \approx 2 \pi$ $9/10$ & 2 / 1 / 1\\
		2.369428 & 1.725740160 & 0.392126651 & 3.120727 & $(-)$ $\theta_p = 5.075$, $\lambda_s \approx 1.000$ & $2 \shortto 3$ / 1 / $1 \shortto 2$\\		
		2.180860 & 1.809333107 & 0.323710140 & 3.137033 & $(-)$ $\theta_p = 5.1656$, $\lambda_s = 1.020$ & 3 / 1 / 2\\
		1.177141 & 2.007218176 & 0.044057442 & 3.154815 & $(-)$ $\theta_p = 5.468$, $\lambda_s \approx 1.000$ & $3 \shortto 2$ / 1 / $2 \shortto 1$\\
		1.149891 & 2.008414118 & 0.037261006 & 3.154932 & $(-)$ $\theta_p = 5.376$, $(-)$ $\theta_s = 6.273$ & 2 / 1 / 1\\
		\hline
		2.954927 & 1.138211689 & 0.669474505 & 3.662961 & $\lambda_p \approx 1.000$, $(+)$ $\theta_s = 1.893$ & $10 \shortto 8$ / $5 \shortto 3$ / 5\\
		2.932161 & 1.140077439 & 0.648174645 & 4.676482 & $\lambda_p = 3.663$, $(+)$ $\theta_s = 2.121$ & 9 / 4 / 5\\
		2.910182 & 1.157730765 & 0.645828225 & 5.696314 & $\lambda_p = 5.874$, $\lambda_s \approx -1.000$ & 9 / 4 / 5\\
		2.897814 & 1.177602262 & 0.650699246 & 6.209380 & $\lambda_p = 5.166$, $\lambda_s = -1.448$ & 9 / 4 / 5\\
		2.879404 & 1.222889259 & 0.657934786 & 6.896250 & $\lambda_p = 3.011$, $\lambda_s \approx -1.000$ & 9 / 4 / 5\\
		2.849071 & 1.323558296 & 0.649558155 & 7.731085 & $\lambda_p \approx 1.000$, $(-)$ $\theta_s = 4.132$ & $8 \shortto 10$ / $3 \shortto 5$ / 5\\
		\hline
		3.171668 & 1.150675461 & 1.177467785 & 1.687145 & $\lambda_p = 1446.7$, $(-)$ $\theta_s = 6.028$ & 3 / 2 / 1\\		
		3.161303 & 1.130689163 & 1.257533668 & 1.697551 & $\lambda_p = 1318.6$, $(-)$ $\theta_s = 6.107$ & 3 / 2 / 1\\
		3.152146 & 1.120395188 & 1.296473224 & 1.707786 & $\lambda_p = 1212.1$, $\lambda_s \approx 1.000$ & $3 \shortto 4$ / 2 / $1 \shortto 2$\\
		3.141768 & 1.110698400 & 1.332709180 & 1.720718 &  $\lambda_p = 1099.1$, $\lambda_s = 1.209$ & 4 / 2 / 2
	\end{tabular}
	\label{data_1}
\end{table}

\begin{table}[ht]\fontsize{10}{10}\selectfont \centering
	\caption{Data for $L_2$ halo (first block) and Butterfly orbits (second block).}
	\begin{tabular}{c|c|c|c|c|c|c}
		$C$ & $x(0)$ & $z(0)$ & $p_y(0)$ & $T/2$ & Floquet multipliers & $\mu_{CZ}$\\
		\hline
		3.152146 & 1.18092135 & 0 & 1.025033880 & 1.707786 & $(+)$ $\lambda_p = 1212$, $\lambda_s \approx 1.000$ & $3 \shortto 4$\\
		3.150734 & 1.180601447 & $-$0.017977855 & 1.022876897 & 1.706461 & $\lambda = 1182$, $(-)$ $\theta = 6.182$ & 3\\
		3.098956 & 1.162963891 & $-$0.116997836 & 0.958135828 & 1.644892 & $\lambda = 421$, $(-)$ $\theta = 5.253$ & 3\\
		3.023395 & 1.109360769 & $-$0.194640671 & 0.888202524 & 1.381807 & $\lambda_1 = 22.68$, $\lambda_2 \approx -1.000$ & 3\\
		3.022991 & 1.108729409 & $-$0.195006649 & 0.887870465 & 1.377603 & $\lambda_1 = 21.79$, $\lambda_2 \approx -1.000$ & 3\\
		3.015176 & 1.082874298 & $-$0.202346868 & 0.881915950 & 1.191093 & $\lambda = 1.005$, $(+)$ $\theta = 2.313$ & 3\\
		\hdashline
		3.017142 & 1.068774054 & $-$0.201385883 & 0.884397328 & 1.084549 & $(+)$ $\theta = 1.658$, $\lambda \approx -1.000$ & 4\\
		3.058029 & 1.011828416 & $-$0.173936086 & 0.931882270 & 0.687173 & $(+)$ $\theta = 0.694$, $\lambda \approx -1.000$ & 4\\
		3.086044 & 0.997378723 & $-$0.153359896 & 0.955359901 & 0.563300 & $(+)$ $\theta_1 = 0.58$, $(-)$ $\theta_2 = 4.381$ & 4\\
		\hline
		3.058029 & 1.011828416 & $-$0.173936086 & 0.931882270 & 1.374347 & $(+)$ $\theta = 1.388$, $\lambda \approx 1.000$ & $\overline{7} \shortto 8$\\
		3.059628 & 1.029984409 & $-$0.174265251 & 0.948298600 & 1.384137 & $(+)$ $\theta = 1.395$, $\lambda = 2.304$ & 7\\
		3.071122 & 1.054982912 & $-$0.169939935 & 0.997846809 & 1.468887 & $(+)$ $\theta = 1.507$, $\lambda = 10.45$ & 7\\
		3.087915 & 1.070107927 & $-$0.158061513 & 1.070763080 & 1.698928 & $(+)$ $\theta = 2.338$, $\lambda = 24.58$ & 7\\
		3.090308 & 1.071395741 & $-$0.155301425 & 1.089723687 & 1.781031 & $\lambda_1 \approx -1.000$, $\lambda_2 = 18.31$ & 7\\
		3.090980 & 1.071437726 & $-$0.154012931 & 1.100790430 & 1.833690 & $\lambda_1 \approx -1.000$, $\lambda_2 = 9.691$ & 7\\
		3.091076 & 1.071220268 & $-$0.153458738 & 1.106609374 & 1.862914 & $(-)$ $\theta = 4.634$, $\lambda = 1.008$ & 7\\
		\hdashline
		3.091076 & 1.071208091 & $-$0.153438822 & 1.106838745 & 1.864086 & $(-)$ $\theta_1 = 5.11$, $(+)$ $\theta_2 = 1.157$ & 8\\
		3.091068 & 1.071125129 & $-$0.153317154 & 1.108280164 & 1.871513 & $-0.016 \pm 2.973 i$ & 8
	\end{tabular}
	\label{data_2}
\end{table}

\begin{table}[ht]\fontsize{10}{10}\selectfont \centering
	\caption{Data for two families branching off from 4:5 DRO (near-Moon).\ Orbits in the first block terminate at 1:2 Butterfly, and orbits in the second block terminate at 1:3 $LPO_2$.}
	\begin{tabular}{c|c|c|c|c|c|c}
		$C$ & $x(0)$ & $z(0)$ & $p_y(0)$ & $T/2$ & Floquet multipliers & $\mu_{CZ}$\\
		\hline
		2.991353 & 0.878519575 & 0 & 1.348629522 & 4.264849 & $(-)$ $\theta_p = 4.450$, $\lambda_s \approx 1.000$ & $16 \shortto 14$\\
		3.002894 & 0.896648068 & $-$0.053362325 & 1.345389316 & 3.992994 & $(-)$ $\theta_1 = 5.86$, $(-)$ $\theta_2 = 6.235$ & 14\\
		3.002971 & 0.895841591 & $-$0.054388749 & 1.341941001 & 3.994214 & $(-)$ $\theta_1 = 5.84$, $(-)$ $\theta_2 = 6.277$ & 14\\
		\hdashline
		2.996087 & 0.880472299 & $-$0.078194918 & 1.289204629 & 4.362041 & $\lambda_1 = -1.332$, $\lambda_2 = 1.241$ & 15\\
		2.995443 & 0.882171647 & $-$0.082565128 & 1.286947698 & 4.416919 & $(+)$ $\theta = 2.902$, $\lambda = 1.001$ & 15\\
		\hdashline
		2.995746 & 0.884878061 & $-$0.085226385 & 1.287043690 & 4.429567 & $\lambda \approx -1.000$, $(-)$ $\theta = 6.155$ & 14\\
		2.996091 & 0.886515013 & $-$0.086437307 & 1.287393349 & 4.429444 & $\lambda_1 = -1.512$, $\lambda_2 \approx 1.000$ & $14 \shortto 15$\\
		3.048061 & 0.993558332 & $-$0.121148669 & 1.293398384 & 3.997389 & $\lambda_1 = -1070$, $\lambda_2 = 3.245$ & 15\\
		3.087852 & 1.058735562 & $-$0.147757768 & 1.153583170 & 3.699461 & $\lambda_1 = -469.8$, $\lambda_2 = 9.689$ & 15\\
		3.090995 & 1.071426404 & $-$0.153964159 & 1.101266567 & 3.672072 & $(-)$ $\theta = 5.926$, $\lambda = 84.72$ & 15\\
		3.090980 & 1.071437726 & $-$0.154012931 & 1.100790430 & 3.667381 & $\lambda_1 \approx 1.000$, $\lambda_2 = 93.92$ & $15 \shortto \overline{14}$\\
		\hline
		$C$ & $x(0)$ & $p_y(0)$ & $p_z(0)$ & $T/2$ & Floquet multipliers & $\mu_{CZ}$\\
		\hline
		2.991353 & 0.878519575 & 1.348629522 & 0 & 4.264849 & $(-)$ $\theta_p = 4.450$, $\lambda_s \approx 1.000$ & $16 \shortto 14$\\
		3.007065 & 0.881942812 & 1.302869520 & 0.181199664 & 3.908588 & $(-)$ $\theta = 6.234$, $\lambda = 2.591$ & 15\\
		3.019162 & 0.884779146 & 1.262850112 & 0.243555725 & 3.707202 & $(-)$ $\theta = 6.199$, $\lambda = 6.931$ & 15\\
		3.079720 & 0.899626154 & 1.006585647 & 0.400725437 & $\approx \pi$ & $(-)$ $\theta = 5.998$, $\lambda = 29.695$ & 15\\
		3.095867 & 0.903746519 & 0.921968551 & 0.408319742 & 3.053279 & $(-)$ $\theta = 5.941$, $\lambda = 20.87$ & 15\\
		3.098307 & 0.904407952 & 0.908532728 & 0.408082217 & 3.041127 & $(-)$ $\theta = 5.926$, $\lambda = 18.56$ & 15\\
		3.143397 & 0.929623322 & 0.630875479 & 0.383310362 & 2.764689 & $\lambda_1 = -1.746$, $\lambda_2 = 1.446$ & 15\\
		3.170595 & 0.948969209 & 0.348763232 & 0.225021757 & 2.575393 & $(-)$ $\theta = 3.374$, $\lambda \approx 1.000$ & $15 \shortto 14$\\
		3.173441 & 0.950879545 & 0.304817848 & 0.150028963 & 2.558345 & $(-)$ $\theta_1 = 3.45$, $(-)$ $\theta_2 = 6.257$ & 14\\
		3.175326 & 0.952158810 & 0.272962994 & 0 & 2.548222 & $(-)$ $\theta = 3.511$, $\lambda \approx 1.000$ & $14 \shortto 16$		
	\end{tabular}
	\label{data_3}
\end{table}

\begin{table}[ht]\fontsize{10}{10}\selectfont \centering
	\caption{Data for two bridge families between 5:6 DRO and 1:4 $LPO_1$ (near-Moon).}
	\begin{tabular}{c|c|c|c|c|c|c}
		$C$ & $x(0)$ & $p_y(0)$ & $p_z(0)$ & $T/4$ & Floquet multipliers & $\mu_{CZ}$\\
		\hline
		3.023442 & 0.899133015 & 1.375632885 & 0 & 1.940225 & $(-)$ $\theta_p = 4.869$, $\lambda_s \approx 1.000$ & $20 \shortto 18$\\		
		3.034573 & 0.900344280 & 1.331771785 & 0.181345924 & 1.857216 & $(-)$ $\theta_1 = 6.25$, $(-)$ $\theta_2 = 6.256$ & 18\\
		3.034575 & 0.900344466 & 1.331765822 & 0.181357790 & 1.857209 & $(-)$ $\theta = 6.269$, $\lambda \approx 1.000$ & $18 \shortto 19$\\
		3.047691 & 0.902150503 & 1.277700806 & 0.263998359 & 1.777938 & $(-)$ $\theta = 6.247$, $\lambda = 4.139$ & 19\\
		3.101297 & 0.910626162 & 1.019508055 & 0.414640995 & $\approx \pi/2$ & $(-)$ $\theta = 6.191$, $\lambda = 12.967$ & 19\\
		3.114242 & 0.912621749 & 0.947011277 & 0.420742443 & 1.537298 & $(-)$ $\theta = 6.194$, $\lambda = 12.1$ & 19\\
		3.119264 & 0.913378191 & 0.917768495 & 0.419623715 & 1.525418 & $(-)$ $\theta = 6.196$, $\lambda = 11.15$ & 19\\
		3.144116 & 0.916949739 & 0.763674405 & 0.377538810 & 1.474071 & $(-)$ $\theta = 6.219$, $\lambda \approx 1.000$ & $19 \shortto 18$\\
		3.145214 & 0.917099822 & 0.756507781 & 0.373899193 & 1.472054 & $(-)$ $\theta_1 = 5.58$, $(-)$ $\theta_2 = 6.220$ & 18\\
		3.158334 & 0.918804331 & 0.668625884 & 0.312412181 & 1.449272 & $\lambda \approx -1.000$, $(-)$ $\theta = 6.242$ & 18\\
		3.167111 & 0.917382291 & 0.618290197 & 0.235554772 & 1.434263 & $(+)$ $\theta_1 = 0.10$, $(-)$ $\theta_2 = 6.231$ & 18\\
		3.167112 & 0.917359130 & 0.618389366 & 0.235489441 & 1.434248 & $1.016 \pm 0.073 i$ & 18\\
		3.167115 & 0.917290425 & 0.618685332 & 0.235298783 & 1.434182 & $\lambda_1 = 1.068$, $\lambda_2 = 1.084$ & 18\\
		3.167224 & 0.915145908 & 0.628320945 & 0.230163537 & 1.431550 & $\lambda_1 \approx 1.000$, $\lambda_2 = 1.840$ & $18 \shortto 19$\\
		3.167737 & 0.912518901 & 0.639317633 & 0.222893621 & 1.425631 & $(+)$ $\theta = 0.041$, $\lambda = 2.048$ & 19\\
		3.190463 & 0.894141140 & 0.674092060 & 0 & 1.294847 & $(+)$ $\theta_p = 1.966$, $\lambda_s \approx 1.000$ & $18 \shortto 20$\\
		\hline
		$C$ & $x(0)$ & $z(0)$ & $p_y(0)$ & $T/4$ & Floquet multipliers & $\mu_{CZ}$\\
		\hline
		3.023442 & 0.899133015 & 0 & 1.375632885 & 1.940225 & $(-)$ $\theta_p = 4.869$, $\lambda_s \approx 1.000$ & $20 \shortto 18$\\
		3.028474 & 0.906550073 & 0.031885062 & 1.376625552 & 1.900585 & $(-)$ $\theta = 5.254$, $\lambda = 1.005$ & 19\\
		3.034577 & 0.915212927 & 0.045774140 & 1.377852624 & 1.857177 & $\lambda_1 \approx 1.000$, $\lambda_2 = 1.014$ & $19 \shortto 20$\\
		3.046948 & 0.931786866 & 0.061759827 & 1.380447603 & 1.781959 & $\lambda_1 = 1.035$, $\lambda_2 = 3.983$ & 20\\
		3.088828 & 0.980395160 & 0.077744583 & 1.390339443 & 1.607633 & $\lambda_1 = 1.088$, $\lambda_2 = 12.05$ & 20\\
		3.101293 & 0.993151949 & 0.076822939 & 1.393522500 & $\approx \pi/2$ & $\lambda_1 = 1.095$, $\lambda_2 = 12.975$ & 20\\
		3.120127 & 1.011279925 & 0.071859659 & 1.398446850 & 1.523452 & $\lambda_1 = 1.090$, $\lambda_2 = 10.95$ & 20\\
		3.144127 & 1.032683355 & 0.058922865 & 1.404889667 & 1.474108 & $\lambda_1 \approx 1.000$, $\lambda_2 = 1.066$ & $20 \shortto 19$\\
		3.149081 & 1.036889398 & 0.055072879 & 1.406270421 & 1.465137 & $(-)$ $\theta = 4.708$, $\lambda = 1.058$ & 19\\
		3.165672 & 1.050078082 & 0.036792367 & 1.413503943 & 1.437536 & $(+)$ $\theta = 1.465$, $\lambda = 1.026$ & 19\\
		3.167241 & 1.046251099 & 0.032312257 & 1.442989430 & 1.431689 & $\lambda_1 \approx 1.000$, $\lambda_2 = 1.838$ & $19 \shortto 18$\\
		3.168000 & 1.042136490 & 0.030128933 & 1.469637231 & 1.423395 & $\lambda_1 = 1.057$, $\lambda_2 = 2.021$ & 18\\
		3.170174 & 1.037764703 & 0.026695877 & 1.503985372 & 1.434248 & $1.375 \pm 0.021 i$ & 18\\
		3.172853 & 1.034735074 & 0.023441079 & 1.533012954 & 1.389323 & $(+)$ $\theta_1 = 0.40$, $(-)$ $\theta_2 = 5.917$ & 18\\
		3.189044 & 1.026441526 & 0.005296439 & 1.651571367 & 1.301618 & $(+)$ $\theta_1 = 1.89$, $(-)$ $\theta_2 = 6.272$ & 18\\
		3.190463 & 1.026012085 & 0 & 1.660229960 & 1.294847 & $(+)$ $\theta_p = 1.966$, $\lambda_s \approx 1.000$ & $18 \shortto 20$
	\end{tabular}
\end{table}

\begin{table}[ht]\fontsize{10}{10}\selectfont \centering
	\caption{Data for two bridge families between 6:7 DRO and 1:5 $LPO_1$ (near-Moon).}
	\begin{tabular}{c|c|c|c|c|c|c}
		$C$ & $x(0)$ & $z(0)$ & $p_y(0)$ & $T/2$ & Floquet multipliers & $\mu_{CZ}$\\
		\hline
		3.049858 & 0.911471544 & 0 & 1.399616868 & 3.698925 & $(-)$ $\theta_p = 5.099$, $\lambda_s \approx 1.000$ & $24 \shortto 22$\\
		3.052483 & 0.914889051 & 0.020999984 & 1.400048039 & 3.668199 & $(-)$ $\theta_1 = 5.26$, $(-)$ $\theta_2 = 6.281$ & 22\\
		3.060770 & 0.924859536 & 0.041059758 & 1.398566177 & 3.578701 & $\lambda \approx 1.000$, $(-)$ $\theta = 6.249$ & $22 \shortto 23$\\
		3.070627 & 0.937177980 & 0.052999061 & 1.403150606 & 3.484861 & $\lambda = 2.852$, $(-)$ $\theta = 6.282$ & 23\\
		3.113445 & 0.982861218 & 0.069176997 & 1.411681391 & 3.186548 & $\lambda = 7.536$, $(-)$ $\theta = 6.282$ & 23\\
		3.121999 & 0.991075037 & 0.068674242 & 1.413650867 & $\approx \pi$ & $\lambda = 7.825$, $(-)$ $\theta = 6.264$ & 23\\
		3.166989 & 1.029149398 & 0.049012836 & 1.437389565 & 2.953076 & $\lambda \approx 1.000$, $(-)$ $\theta = 6.247$ & $23 \shortto 22$\\
		3.171543 & 1.033725205 & 0.046351626 & 1.429910005 & 2.937816 & $(-)$ $\theta_1 = 5.14$, $(-)$ $\theta_2 = 6.282$ & 22\\
		3.188070 & 1.044642765 & 0.025853312 & 1.446774428 & 2.882935 & $\lambda \approx -1.000$, $(-)$ $\theta = 6.281$ & 22\\
		3.195982 & 1.045240884 & 0 & 1.479461010 & 2.846460 & $(+)$ $\theta_p = 2.004$, $\lambda_s \approx 1.000$ & $22 \shortto 24$\\
		\hline
		$C$ & $x(0)$ & $p_y(0)$ & $p_z(0)$ & $T/2$ & Floquet multipliers & $\mu_{CZ}$\\
		\hline
		3.049858 & 0.911471544 & 1.399616868 & 0 & 3.698925 & $(-)$ $\theta_p = 5.099$, $\lambda_s \approx 1.000$ & $24 \shortto 22$\\
		3.053302 & 0.911683970 & 1.383277181 & 0.114999807 & 3.658839 & $(-)$ $\theta = 5.315$, $\lambda = 1.002$ & 23\\
		3.060770 & 0.912245765 & 1.347584233 & 0.201828801 & 3.578587 & $\lambda_1 \approx 1.000$, $\lambda_2 = 1.004$ & $23 \shortto 24$\\
		3.072506 & 0.913324244 & 1.290149484 & 0.283995426 & 3.468362 & $\lambda_1 = 3.116$, $\lambda_2 = 1.003$ & 24\\
		3.121992 & 0.918861573 & 1.020286235 & 0.431870436 & $\approx \pi$ & $\lambda_1 = 7.825$, $\lambda_2 = 1.002$ & 24\\
		3.132402 & 0.920031680 & 0.956792265 & 0.436585582 & 3.091485 & $\lambda_1 = 7.633$, $\lambda_2 = 1.001$ & 24\\
		3.166987 & 0.923347883 & 0.729313280 & 0.369804972 & 2.953462 & $\lambda_1 \approx 1.000$, $\lambda_2 = 1.052$ & $24 \shortto 23$\\
		3.171561 & 0.923627589 & 0.697900158 & 0.346995193 & 2.937754 & $(-)$ $\theta = 5.140$, $\lambda = 1.002$ & 23\\
		3.188022 & 0.923231608 & 0.590498121 & 0.205008489 & 2.883107 & $\lambda_1 \approx -1.000$, $\lambda_2 = 1.002$ & 23\\
		3.195982 & 0.919311357 & 0.565865253 & 0 & 2.846460 & $(+)$ $\theta_p = 2.004$, $\lambda_s \approx 1.000$ & $22 \shortto 24$
	\end{tabular}
\end{table}

\begin{table}[ht]\fontsize{10}{10}\selectfont \centering
	\caption{Data for two bridge families between 7:8 DRO and 1:6 $LPO_1$ (near-Moon).}
	\begin{tabular}{c|c|c|c|c|c|c}
		$C$ & $x(0)$ & $p_y(0)$ & $p_z(0)$ & $T/4$ & Floquet multipliers & $\mu_{CZ}$\\
		\hline
		3.073375 & 1.056389717 & 0.560205596 & 0 & 1.795175 & $(-)$ $\theta_p = 5.254$, $\lambda_s \approx 1.000$ & $28 \shortto 26$\\
		3.079033 & 1.055556436 & 0.587591186 & 0.159991073 & 1.769417 & $(-)$ $\theta_1 = 5.58$, $(-)$ $\theta_2 = 6.280$ & 26\\
		3.084180 & 1.054878497 & 0.613676691 & 0.219450559 & 1.747486 & $\lambda \approx 1.000$, $(-)$ $\theta = 6.282$ & $26 \shortto 27$\\
		3.100250 & 1.053147179 & 0.701259102 & 0.334510388 & 1.687008 & $\lambda = 3.121$, $(-)$ $\theta = 6.280$ & 27\\
		3.141516 & 1.050280887 & 0.955919827 & 0.449017575 & $\approx \pi/2$ & $\lambda = 5.610$, $(-)$ $\theta = 6.258$ & 27\\
		3.147898 & 1.049943469 & 0.998051328 & 0.451047664 & 1.556407 & $\lambda = 5.586$, $(-)$ $\theta = 6.280$ & 27\\
		3.186774 & 1.047769984 & 1.268435529 & 0.363347332 & 1.482531 & $\lambda \approx 1.000$, $(-)$ $\theta = 6.281$ & $27 \shortto 26$\\
		3.189115 & 1.047590585 & 1.285693075 & 0.350482502 & 1.478687 & $(-)$ $\theta_1 = 5.63$, $(-)$ $\theta_2 = 6.281$ & 26\\
		3.211617 & 1.044349431 & 1.467797343 & 0 & 1.443191 & $(-)$ $\theta_p = 3.526$, $\lambda_s \approx 1.000$ & $26 \shortto 28$\\
		\hline
		3.073375 & 0.920081572 & 1.421306845 & 0 & 1.795175 & $(-)$ $\theta_p = 5.2543$, $\lambda_s \approx 1.000$ & $28 \shortto 26$\\
		3.078944 & 0.920350629 & 1.390852484 & 0.159998135 & 1.769809 & $(-)$ $\theta = 5.580$, $\lambda = 1.005$ & 27\\
		3.084180 & 0.920654083 & 1.362051851 & 0.219999985 & 1.747486 & $\lambda_1 \approx 1.000$, $\lambda_2 = 1.001$ & $27 \shortto 28$\\
		3.141516 & 0.925368140 & 1.018989093 & 0.449981646 & $\approx \pi / 2$ & $\lambda_1 = 5.610$, $\lambda_2 = 1.002$ & 28\\
		3.150489 & 0.926150204 & 0.959390132 & 0.453645807 & 1.550771 & $\lambda_1 = 2.856$, $\lambda_2 = 1.001$ & 28\\
		3.186774 & 0.928762794 & 0.702555006 & 0.368603425 & 1.482532 & $\lambda_1 \approx 1.000$, $\lambda_2 = 1.004$ & $28 \shortto 27$\\
		3.190925 & 0.928932673 & 0.672181743 & 0.343641775 & 1.475742 & $(-)$ $\theta = 5.402$, $\lambda = 1.005$ & 27\\
		3.211617 & 0.928295714 & 0.529367971 & 0 & 1.443191 & $(-)$ $\theta_p = 3.526$, $\lambda_s \approx 1.000$ & $26 \shortto 28$	
	\end{tabular}
\end{table}

\begin{table}[ht]\fontsize{10}{10}\selectfont \centering
	\caption{Data for two bridge families between 8:9 DRO and 1:7 $LPO_1$ (near-Moon).}
	\begin{tabular}{c|c|c|c|c|c|c}
		$C$ & $x(0)$ & $z(0)$ & $p_y(0)$ & $T/2$ & Floquet multipliers \& Krein sign & $\mu_{CZ}$\\
		\hline
		3.094976 & 0.926549812 & 0 & 1.441161202 & 3.517787 & $(-)$ $\theta_p = 5.367$, $\lambda_s \approx 1.000$ & $32 \shortto 30$\\
		3.100402 & 0.932400272 & 0.024998649 & 1.441916135 & 3.476305 & $(-)$ $\theta_1 = 5.65$, $(-)$ $\theta_2 = 6.281$ & 30\\
		3.105699 & 0.937989216 & 0.034079913 & 1.442662720 & 3.438154 & $\lambda \approx 1.000$, $(-)$ $\theta = 6.275$ & $30 \shortto 31$\\
		3.118206 & 0.950735785 & 0.046488131 & 1.444525766 & 3.356045 & $\lambda = 2.471$, $(-)$ $\theta = 6.281$ & 31\\
		3.154661 & 0.984826386 & 0.057309824 & 1.450904297 & 3.165059 & $\lambda = 4.363$, $(-)$ $\theta = 6.280$ & 31\\
		3.160029 & 0.989509789 & 0.057096114 & 1.452009918 & $\approx \pi$ & $\lambda = 4.427$, $(-)$ $\theta = 6.264$ & 31\\
		3.201442 & 1.023012199 & 0.042291688 & 1.463547284 & 2.988422 & $\lambda = 1.873$, $(-)$ $\theta = 6.280$ & 31\\
		3.204832 & 1.025535246 & 0.039783995 & 1.464873745 & 2.977633 & $\lambda \approx 1.000$, $(-)$ $\theta = 6.281$ & $31 \shortto 30$\\	
		3.222223 & 1.037675902 & 0.019346841 & 1.474619836 & 2.925028 & $(-)$ $\theta_1 = 4.53$, $(-)$ $\theta_2 = 6.280$ & 30\\
		3.227134 & 1.040760999 & 0 & 1.478776434 & 2.910670 & $(-)$ $\theta_p = 4.199$, $\lambda_s \approx 1.000$ & $30 \shortto 32$\\
		\hline
		$C$ & $x(0)$ & $p_y(0)$ & $p_z(0)$ & $T/2$ & Floquet multipliers & $\mu_{CZ}$\\
		\hline
		3.094976 & 0.926549812 & 1.441161202 & 0 & 3.517787 & $(-)$ $\theta_p = 5.367$, $\lambda_s \approx 1.000$ & $32 \shortto 30$\\
		3.097622 & 0.926648397 & 1.425028209 & 0.119990453 & 3.497236 & $(-)$ $\theta = 5.496$, $\lambda = 1.004$ & 31\\
		3.105700 & 0.927007344 & 1.375404733 & 0.236317951 & 3.438148 & $\lambda_1 \approx 1.000$, $\lambda_2 = 1.004$ & $31 \shortto 32$\\
		3.115030 & 0.927512665 & 1.317309713 & 0.314304136 & 3.375918 & $\lambda_1 = 2.207$, $\lambda_2 = 1.001$ & 32\\
		3.160029 & 0.930582323 & 1.017307059 & 0.467618885 & $\approx \pi$ & $\lambda_1 = 4.427$, $\lambda_2 = 1.002$ & 32\\
		3.173532 & 0.931523601 & 0.919491564 & 0.468292399 & 3.086662 & $\lambda_1 = 4.247$, $\lambda_2 = 1.003$ & 32\\
		3.204835 & 0.933269070 & 0.680265555 & 0.370571435 & 2.977628 & $\lambda_1 \approx 1.000$, $\lambda_2 = 1.003$ & $32 \shortto 31$\\
		3.207261 & 0.933351199 & 0.661322880 & 0.354336216 & 2.970032 & $(-)$ $\theta = 5.721$, $\lambda = 1.002$ & 31\\
		3.227134 & 0.933280794 & 0.510047445 & 0 & 2.910670 & $(-)$ $\theta_p = 4.199$, $\lambda_s \approx 1.000$ & $30 \shortto 32$
	\end{tabular}
\end{table}

\begin{table}[ht]\fontsize{10}{10}\selectfont \centering
	\caption{Data for two bridge families between 9:10 DRO and 1:8 $LPO_1$ (near-Moon).}
	\begin{tabular}{c|c|c|c|c|c|c}
		$C$ & $x(0)$ & $p_y(0)$ & $p_z(0)$ & $T/4$ & Floquet multipliers & $\mu_{CZ}$\\
		\hline
		3.115178 & 0.931641145 & 1.459521304 & 0 & 1.732752 & $(-)$ $\theta_p = 5.456$, $\lambda_s \approx 1.000$ & $36 \shortto 34$\\
		3.119600 & 0.931785355 & 1.429853511 & 0.164999730 & 1.718059 & $(-)$ $\theta_1 = 5.66$, $(-)$ $\theta_2 = 6.278$ & 34\\
		3.125813 & 0.932021804 & 1.387856472 & 0.251103312 & 1.698494 & $\lambda \approx 1.000$, $(-)$ $\theta = 6.279$ & $34 \shortto 35$\\
		3.134644 & 0.932417837 & 1.327505946 & 0.329977216 & 1.672591 &  $\lambda = 2.008$, $(-)$ $\theta = 6.280$ & 35\\
		3.177698 & 0.934850251 & 1.015645302 & 0.484456560 & $\approx \pi / 2$ &  $\lambda = 3.706$, $(-)$ $\theta = 6.255$ & 35\\
		3.184219 & 0.935233477 & 0.965324108 & 0.486840713 & 1.558104 &  $\lambda = 3.686$, $(-)$ $\theta = 6.278$ & 35\\
		3.221771 & 0.937044737 & 0.660999965 & 0.374237744 & 1.494744 & $\lambda \approx 1.000$, $(-)$ $\theta = 6.280$ & $35 \shortto 34$\\
		3.226928 & 0.937179216 & 0.618084670 & 0.332420467 & 1.487085 & $(-)$ $\theta_1 = 5.549$, $(-)$ $\theta_2 = 6.278$ & 34\\
		3.242208 & 0.937193800 & 0.493146677 & 0 & 1.465384 & $(-)$ $\theta_p = 4.598$, $\lambda_s \approx 1.000$ & $34 \shortto 36$\\
		\hline
		3.115178 & 1.044501300 & 0.520438208 & 0 & 1.732752 & $(-)$ $\theta_p = 5.456$, $\lambda_s \approx 1.000$ & $36 \shortto 34$\\
		3.119115 & 1.044135813 & 0.544784372 & 0.154994331 & 1.719639 & $(-)$ $\theta = 5.636$, $\lambda = 1.004$ & 35\\
		3.125813 & 1.043576961 & 0.587679840 & 0.250976313 & 1.698494 & $\lambda_1 \approx 1.000$, $\lambda_2 = 1.002$ & $35 \shortto 36$\\
		3.134647 & 1.042938743 & 0.646548567 &0.331247212 & 1.672590 & $\lambda_1 = 2.008$, $\lambda_2 = 1.003$ & 36\\
		3.177700 & 1.040820897 & 0.959736089 & 0.483960150 & $\approx \pi / 2$ & $\lambda_1 = 3.706$, $\lambda_2 = 1.003$ & 36\\
		3.182188 & 1.040650606 & 0.994257224 & 0.485364677 & 1.561999 & $\lambda_1 = 3.702$, $\lambda_2 = 1.003$ & 36\\
		3.221771 & 1.038998140 & 1.312006530 & 0.371484335 & 1.494744 & $\lambda_1 \approx 1.000$, $\lambda_2 = 1.004$ & $36 \shortto 35$\\
		3.225941 & 1.038755178 & 1.347343747 & 0.339609109 & 1.488536 & $(-)$ $\theta = 5.629$, $\lambda = 1.006$ & 35\\
		3.242208 & 1.037429436 & 1.492054156 & 0 & 1.465384 & $(-)$ $\theta_p = 4.598$, $\lambda_s \approx 1.000$ & $34 \shortto 36$
	\end{tabular}
	\label{data_8}
\end{table}

\begin{table}[ht]\fontsize{10}{10}\selectfont \centering
	\caption{Data for two families branching off from 4:5 DRO (orbiting around Moon, $L_1$ and $L_2$).}
	\begin{tabular}{c|c|c|c|c|c|c}
		$C$ & $x(0)$ & $z(0)$ & $p_y(0)$ & $T/2$ & Floquet multipliers & $\mu_{CZ}$\\
		\hline
		2.902460 & 1.227777496 & 0 & 0.680344137 & 10.058570 & $(+)$ $\theta_p = 0.450$, $\lambda_s \approx 1.000$ & $14 \shortto 16$\\
		2.887008 & 1.233558431 & 0.152193081 & 0.692214949 & 10.657543 & $(+)$ $\theta_1 = 2.26$, $(-)$ $\theta_2 = 5.923$ & 14\\
		2.884441 & 1.230229886 & 0.167351206 & 0.693979032 & 10.721459 & $\lambda \approx -1.000$, $(-)$ $\theta = 5.845$ & 14\\
		2.862993 & 1.190939276 & 0.263095604 & 0.708498892 & 11.111291 & $\lambda \approx -1.000$, $(-)$ $\theta = 5.446$ & 14\\
		2.861486 & 1.181404253 & 0.269469374 & 0.710071406 & 11.127041 & $(-)$ $\theta_1 = 4.63$, $(-)$ $\theta_2 = 6.261$ & 14\\
		\hdashline
		2.901545 & 1.095678971 & 0.223394478 & 0.714679318 & 9.550657 & $\lambda_1 = -84.6$, $\lambda_2 = 1.046$ & 15\\
		\hdashline
		2.901501 & 1.094189255 & 0.225212177 & 0.715703033 & 9.548664 & $\lambda = -89.2$, $(-)$ $\theta = 4.924$ & 14\\
		2.881400 & 1.058643839 & 0.297521045 & 0.747946484 & 10.253985 & $-26.31 \pm 240.4i$ & 14\\
		2.818602 & 0.995556906 & 0.417278841 & 0.790148572 & 11.336385 & $\lambda_1 = 411$, $\lambda_2 = 4953$ & 14\\
		2.723713 & 0.923752612 & 0.531345808 & 0.819855209 & 11.865775 & $-645 \pm 4056i$ & 14\\
		1.101449 & $-$0.004691495 & 1.045382961 & 0.892698358 & $\approx 4 \pi$ & $\lambda_1 = - 37869$, $\lambda_2 = -3490$ & 14\\
		0.065386 & $-$0.573372922 & 0.891976988 & 0.893472446 & 12.612233 & $3709 \pm 1534i$ & 14\\
		$-$0.855996 & $-$1.223280692 & 0.151305754 & 0.772234834 & 12.839582 & $\lambda_1 = -446.3$, $\lambda_2 = -244.2$ & 14\\
		$-$0.856079 & $-$1.227939641 & 0.146627723 & 0.769109347 & 12.845075 & $\lambda_1 = -547.6$, $(+)$ $\theta = 1.772$ & 14\\
		\hdashline
		0.458403 & $-$2.101511973 & 0.017577852 & 0.114651646 & 13.886584 & $\lambda_1 = -6 \cdot 10^8$, $\lambda_2 = 9 \cdot 10^6$ & 13\\
		\hline
		$C$ & $x(0)$ & $p_y(0)$ & $p_z(0)$ & $T/2$ & Floquet multipliers & $\mu_{CZ}$\\
		\hline 2.902460 & 0.760777742 & 1.343426379 & 0 & 10.058570 & $(+)$ $\theta_p = 0.450$, $\lambda_s \approx 1.000$ & $14 \shortto 16$\\
		2.892579 & 0.736806021 & 1.345615450 & 0.120001160 & 10.499085 & $(+)$ $\theta = 1.538$, $\lambda = 1.192$ & 15\\
		2.848123 & 0.720232263 & 1.306281093 & 0.340002427 & 11.296762 & $(-)$ $\theta = 4.964$, $\lambda = 4.577$ & 15\\
		2.844857 & 0.720113576 & 1.302739933 & 0.350846667 & 11.330178 & $\lambda_1 \approx 1.000$, $\lambda_2 = 8.781$ & $15 \shortto 16$\\
		2.410131 & 0.767272579 & 0.886008194 & 0.899714988 & 12.317352 & $\lambda_1 = 48.3$, $\lambda_2 = 22960$ & 16\\
		0.979348 & 0.838255406 & $-$0.062576131 & 1.182183565 & $\approx 4 \pi$ & $\lambda_1 = 17.484$, $\lambda_2 = 32796$ & 16\\
		0.738941 & 0.845321222 & $-$0.205437105 & 1.160248944 & 12.578450 & $\lambda_1 \approx 1.000$, $\lambda_2 = 27816$ & $16 \shortto 15$\\
		0.689354 & 0.846678404 & $-$0.234596400 & 1.153688664 & 12.580728 & $(-)$ $\theta = 4.738$, $\lambda = 26706$ & 15\\
		0.313072 & 0.851736238 & $-$0.452509302 & 1.084850417 & 12.597104 & $(+)$ $\theta = 0.020$, $\lambda = 18485$ & 15\\
		\hdashline
		0.329663 & 0.847426320 & $-$0.442533666 & 1.093631801 & 12.596693 & $\lambda_1 = 1.024$, $\lambda_2 = 19121$ & 14\\
		\hdashline
		0.329500 & 0.847341164 & $-$0.442616705 & 1.093698924 & 12.596703 & $(+)$ $\theta = 0.003$, $\lambda = 19121$ & 15\\
		$-$0.017901 & 0.858651209 & $-$0.640548951 & 0.982117910 & 12.610735 & $\lambda_1 = -2.19$, $\lambda_2 = 9074$ & 15\\
		$-$0.401101 & 0.877229234 & $-$0.851419993 & 0.789676614 & 12.633233 & $(-)$ $\theta = 5.015$, $\lambda = 3341$ & 15\\
		$-$0.468907 & 0.881262982 & $-$0.887685217 & 0.745537965 & 12.639787 & $\lambda_1 \approx 1.000$, $\lambda_2 = 2721$ & $15 \shortto 16$\\
		$-$0.766344 & 1.076639926 & $-$0.871821841 & 0.466195898 & 12.720661 & $\lambda_1 = 12.60$, $\lambda_2 = 1084$ & 16\\
		$-$0.812528 & 1.073911941 & $-$0.896353559 & 0.430638157 & 12.767665 & $\lambda_1 \approx 1.000$, $\lambda_2 = 994$ & $16 \shortto 15$\\
		$-$0.814362 & 1.073834989 & $-$0.897264065 & 0.429201865 & 12.770353 & $(-)$ $\theta = 4.970$, $\lambda = 990$ & 15\\
		$-$0.859711 & 1.076672891 & $-$0.910050532 & 0.400067514 & 12.957931 & $(+)$ $\theta = 1.028$, $\lambda = 1828$ & 15\\
		$-$0.859717 & 1.076732689 & $-$0.909931196 & 0.400176022 & 12.959709 & $\lambda_1 \approx 1.000$, $(+)$ $\lambda_2 = 1968$ & $15 \shortto 14$\\
		$-$0.859717 & 1.076768576 & $-$0.909858646 & 0.400243043 & 12.960769 & $\lambda_1 = 1.052$, $\lambda_2 = 2057$ & 14\\
		\hdashline
		$-$0.859718 & 1.076769254 & $-$0.909857269 & 0.400244322 & 12.960784 & $(-)$ $\theta = 6.265$, $\lambda = 2058$ & 13\\
		$-$0.660349 & 1.092923707 & $-$0.790543466 & 0.571455651 & 13.658223 & $\lambda_1 = - 728.6$, $\lambda_2 = 10^9$ & 13\\
		$-$0.659855 & 1.092935872 & $-$0.790300865 & 0.571765931 & 13.659229 & $\lambda_1 = 6814$, $\lambda_2 = 10^9$ & 12\\
		$-$0.112798 & 1.099380161 & $-$0.533645246 & 0.806266047 & 14.494830 & $\lambda_1 = 10^8$, $\lambda_2 = 10^{12}$ & 12
	\end{tabular}
	\label{data_9}
\end{table}

\begin{table}[ht]\fontsize{10}{10}\selectfont \centering
	\caption{Data for two families bifurcating from 5:6 DRO (orbiting around Moon, $L_1$ and $L_2$).\ Orbits in the first block terminate at period-doubling bifurcation of P3DRO.}
	\begin{tabular}{c|c|c|c|c|c|c}
		$C$ & $x(0)$ & $z(0)$ & $p_y(0)$ & $T/4$ & Floquet multipliers & $\mu_{CZ}$\\
		\hline
		2.883609 & 1.263819475 & 0 & 0.672130625 & 6.789272 & $(-)$ $\theta_p = 4.856$, $\lambda_s \approx 1.000$ & $16 \shortto 18$\\
		2.882933 & 1.257702369 & 0.049942693 & 0.672811605 & 6.789052 & $(-)$ $\theta_1 = 4.90$, $(-)$ $\theta_2 = 6.256$ & 16\\
		2.882472 & 1.247425196 & 0.073942909 & 0.673403647 & 6.775078 & $(-)$ $\theta_1 = 4.98$, $(-)$ $\theta_2 = 6.281$ & 16\\
		\hdashline
		2.883547 & 1.232074452 & 0.084371028 & 0.672801120 & 6.715695 & $(-)$ $\theta = 5.156$, $\lambda = 1.607$ & 17\\
		2.895376 & 1.189275043 & 0.057269604 & 0.662331503 & $\approx 2 \pi$ & $(-)$ $\theta = 5.170$, $\lambda = 10.901$ & 17\\
		2.902757 & 1.171985055 & 0.037701797 & 0.654546148 & 5.995918 & $(-)$ $\theta = 5.316$, $\lambda = 21.14$ & 17\\
		2.910182 & 1.157730765 & 0 & 0.645828225 & 5.696314 & $\lambda_s \approx 1.000$, $(+)$ $\lambda_p = 34.51$ & $17 \shortto \overline{18}$\\
		\hline
		2.883609 & 0.727149837 & 0 & 1.368047580 & 6.789272 & $(-)$ $\theta_p = 4.856$, $\lambda_s \approx 1.000$ & $16 \shortto 18$\\
		2.873111 & 0.699533893 & 0.109999797 & 1.361725284 & 6.740019 & $(-)$ $\theta = 5.066$, $\lambda = 17.11$ & 17\\
		2.866265 & 0.679099637 & 0.120629138 & 1.374673105 & 6.762837 & $\lambda_1 \approx 1.000$, $\lambda_2 = 72.44$ & $17 \shortto 18$\\
		2.855942 & 0.650249375 & 0.124401702 & 1.402450738 & 6.817202 & $\lambda_1 = 20.29$, $\lambda_2 = 244.61$ & 18\\
		2.796027 & 0.499643883 & 0.070087950 & 1.651813954 & 6.791945 & $\lambda_1 = 5408$, $\lambda_2 = 2.4 \cdot 10^6$ & 18\\
		2.109314 & 0.071005416 & 0.001741355 & 4.726602250 & 6.538984 & $\lambda_1 = 2 \cdot 10^8$, $\lambda_2 = 4 \cdot 10^{11}$ & 18
	\end{tabular}
	\label{data_10}
\end{table}

\begin{table}[ht]\fontsize{10}{10}\selectfont \centering
	\caption{Data for two families branching off from 6:7 DRO (orbiting around Moon, $L_1$ and $L_2$).}
	\begin{tabular}{c|c|c|c|c|c|c}
		$C$ & $x(0)$ & $z(0)$ & $p_y(0)$ & $T/2$ & Floquet multipliers & $\mu_{CZ}$\\
		\hline 
		2.869425 & 0.702822759 & 0 & 1.390656122 & 16.893438 & $(-)$ $\theta_p = 3.200$, $\lambda_s \approx 1.000$ & $20 \shortto 22$\\
		2.862220 & 0.697060926 & 0.080006015 & 1.384419080 & 16.985100 & $\lambda = -1.334$, $(-)$ $\theta = 6.279$ & 20\\
		2.849630 & 0.660230761 & 0.124421284 & 1.398418133 & 16.884013 & $(-)$ $\theta_1 = 3.86$, $(-)$ $\theta_2 = 6.016$ & 20\\
		2.849433 & 0.657360714 & 0.124999395 & 1.400523748 & 16.844260 & $(-)$ $\theta = 5.181$, $\lambda \approx 1.000$ & $20 \shortto 21$\\
		2.846566 & 0.639085109 & 0.126377153 & 1.417271966 & 16.682028 & $\lambda_1 = -2.783$, $\lambda_2 = 91.96$ & 21\\
		2.809955 & 0.528773143 & 0.084606106 & 1.591908752 & 16.647334 & $\lambda_1 = -3 \cdot 10^5$, $\lambda_2 = 1481$ & 21\\
		2.514648 & 0.235911021 & 0.010562262 & 2.587610595 & 16.611511 & $\lambda_1 = -7 \cdot 10^7$, $\lambda_2 = 3 \cdot 10^7$ & 21\\
		2.361166 & 0.157623404 & 0.005160084 & 3.211220633 & 16.522951 & $\lambda_1 = -3 \cdot 10^8$, $\lambda_2 = 5 \cdot 10^7$ & 21\\
		\hline 
		$C$ & $x(0)$ & $p_y(0)$ & $p_z(0)$ & $T/2$ & Floquet multipliers & $\mu_{CZ}$\\
		\hline 
		2.869425 & 0.702822759 & 1.390656122 & 0 & 16.893438 & $(-)$ $\theta_p = 3.201$, $\lambda_s \approx 1.000$ & $20 \shortto 22$\\
		2.866721 & 0.704051045 & 1.387261296 & 0.080004782 & 16.933426 & $(-)$ $\theta = 3.285$, $\lambda = 1.001$ & 21\\
		2.863451 & 0.705896110 & 1.382741920 & 0.120002196 & 16.972988 & $\lambda_1 = -1.185$, $\lambda_2 = 1.001$ & 21\\
		2.849557 & 0.719806033 & 1.352920440 & 0.239997803 & 16.958692 & $\lambda_1 = -2.493$, $\lambda_2 = 1.429$ & 21\\
		2.847256 & 0.725441559 & 1.343991881 & 0.259998147 & 16.823461 & $(-)$ $\theta = 5.512$, $\lambda = 7.613$ & 21\\
		2.844979 & 0.725800280 & 1.342583035 & 0.267218510 & 16.816939 & $\lambda_1 \approx 1.000$, $\lambda_2 = 26.87$ & $21 \shortto 22$\\
		2.840253 & 0.724688508 & 1.341281574 & 0.280001683 & 16.863633 & $\lambda_1 = 5.366$, $\lambda_2 = 78.33$ & 22\\
		2.767683 & 0.734231651 & 1.243401595 & 0.504994627 & 16.680694 & $\lambda_1 = 4921$, $\lambda_2 = 10^7$ & 22\\
		2.274824 & 0.762747528 & 0.817171133 & 0.980608953 & 16.332169 & $\lambda_1 = 10^8$, $\lambda_2 = 9 \cdot 10^{11}$ & 22
	\end{tabular}
\end{table}

\begin{table}[ht]\fontsize{10}{10}\selectfont \centering
	\caption{Data for two families bifurcating from 7:8 DRO (orbiting around Moon, $L_1$ and $L_2$).}
	\begin{tabular}{c|c|c|c|c|c|c}
		$C$ & $x(0)$ & $z(0)$ & $p_y(0)$ & $T/4$ & Floquet multipliers & $\mu_{CZ}$\\
		\hline
		2.857370 & 1.310084964 & 0 & 0.655415149 & 10.067271 & $(+)$ $\theta_p = 1.644$, $\lambda_s \approx 1.000$ & $24 \shortto 26$\\
		2.852068 & 1.294922701 & 0.120014506 & 0.659555090 & 10.106681 & $(+)$ $\theta_1 = 2.01$, $(-)$ $\theta_2 = 6.262$ & 24\\
		2.829104 & 1.215575674 & 0.267261947 & 0.676919616 & 10.189391 & $(-)$ $\theta_1 = 3.30$, $(-)$ $\theta_2 = 5.656$ & 24\\
		2.827624 & 1.198415709 & 0.274859748 & 0.678501959 & 10.146970 & $(-)$ $\theta_1 = 4.01$, $(-)$ $\theta_2 = 6.278$ & 24\\
		\hdashline
		2.827633 & 1.197187433 & 0.274974112 & 0.678539106 & 10.141769 & $(-)$ $\theta = 4.122$, $\lambda = 1.321$ & 25\\
		2.828279 & 1.188583148 & 0.273546343 & 0.678298865 & 10.092039 & $\lambda_1 \approx 1.000$, $\lambda_2 = 2.655$ & $25 \shortto 26$ \\
		2.829595 & 1.181732661 & 0.268559925 & 0.676972645 & 10.027259 & $\lambda_1 = 4.1024$, $\lambda_2 = 14.19$ & 26\\
		2.830651 & 1.176537867 & 0.259889939 & 0.673549569 & 9.943485 & $\lambda_1 = 1.025$, $\lambda_2 = 72.22$ & 26\\
		\hdashline
		2.830651 & 1.176537462 & 0.259888773 & 0.673548980 & 9.943477 & $(+)$ $\theta = 0.014$, $\lambda = 72.23$ & 27\\
		2.829564 & 1.175590817 & 0.253083398 & 0.669164274 & 9.892097 & $\lambda_1 \approx 1.000$, $\lambda_2 = 184.80$ & $27 \shortto 26$\\
		2.746334 & 1.267835403 & 0.170067449 & 0.585651807 & 9.874530 & $\lambda_1 = 17990$, $\lambda_2 = 10^8$ & 26\\
		1.990165 & 1.554905352 & 0.030118072 & 0.239267041 & 9.737470 & $\lambda_1 = 4 \cdot 10^8$, $\lambda_2 = 3 \cdot 10^{12}$ & 26\\
		\hline
		2.857370 & 0.683131870 & 0 & 1.411572193 & 10.067271 & $(+)$ $\theta_p = 1.644$, $\lambda_s \approx 1.000$ & $24 \shortto 26$\\
		2.838238 & 0.672431060 & 0.120004882 & 1.399465730 & 10.188101 & $(+)$ $\theta = 2.745$, $\lambda = 1.360$ & 25\\
		2.830154 & 0.664243725 & 0.140003843 & 1.399281536 & 10.224162 & $(+)$ $\theta = 3.101$, $\lambda = 1.903$ & 25\\
		2.812588 & 0.635120127 & 0.164768531 & 1.416924068 & 10.276771 & $(-)$ $\theta = 4.449$, $\lambda = 5.160$ & 25\\
		2.808903 & 0.625240609 & 0.166011387 & 1.426826471 & 10.279462 & $(-)$ $\theta = 5.621$, $\lambda = 5.916$ & 25\\
		2.808463 & 0.623854653 & 0.165989329 & 1.428344083 & 10.279289 & $\lambda_1 \approx 1.000$, $\lambda_2 = 5.963$ & $25 \shortto 26$\\
		2.793010 & 0.520792111 & 0.101888191 & 1.602562239 & 10.025204 & $\lambda_1 = 48.08$, $\lambda_2 = 1.3 \cdot 10^5$ & 26\\
		2.761923 & 0.456420375 & 0.064647879 & 1.748169825 & 9.911641 & $\lambda_1 = 1393$, $\lambda_2 = 1.4 \cdot 10^7$ & 26\\
		2.219447 & 0.103854094 & 0.003014433 & 3.956587053 & 9.705477 & $\lambda_1 = 2 \cdot 10^8$, $\lambda_2 = 4 \cdot 10^{12}$ & 26
	\end{tabular}
\end{table}

\begin{table}[ht]\fontsize{10}{10}\selectfont \centering
	\caption{Data for two families branching off from 8:9 DRO (orbiting around Moon, $L_1$ and $L_2$).\ Orbits in the third block are the Earth--Moon cycling-type orbits.}
	\begin{tabular}{c|c|c|c|c|c|c}
		$C$ & $x(0)$ & $p_y(0)$ & $p_z(0)$ & $T/2$ & Floquet multipliers & $\mu_{CZ}$\\
		\hline
		2.846622 & 0.666380578 & 1.431077272 & 0 & 23.339534 & $(+)$ $\theta_p = 0.144$, $\lambda_s \approx 1.000$ & $28 \shortto 30$\\
		2.801178 & 0.671893540 & 1.386510394 & 0.323663230 & 23.806531 & $\lambda = -1.000$, $(-)$ $\theta = 6.171$ & 28\\
		2.793465 & 0.675852863 & 1.374152362 & 0.353861512 & 23.845881 & $\lambda = -1.000$, $(-)$ $\theta = 6.169$ & 28\\
		2.788525 & 0.680327232 & 1.363034468 & 0.374617259 & 23.857801 & $(-)$ $\theta = 4.131$, $\lambda \approx 1.000$ & $28 \shortto 29$\\
		2.782971 & 0.692611174 & 1.337791260 & 0.407614474 & 23.821899 & $\lambda_1 = 1.000$, $\lambda_2 = 11.45$ & $29 \shortto 30$\\
		2.580634 & 0.740471948 & 1.069092000 & 0.763547926 & 22.919475 & $\lambda_1 = 5 \cdot 10^5$, $\lambda_2 = 2 \cdot 10^{12}$ & 30\\
		2.340766 & 0.757791871 & 0.869044682 & 0.944872825 & 22.719169 & $\lambda_1 = 4 \cdot 10^7$, $\lambda_2 = 2 \cdot 10^{13}$ & 30\\
		\hline 
		$C$ & $x(0)$ & $z(0)$ & $p_y(0)$ & $T/2$ & Floquet multipliers & $\mu_{CZ}$\\
		\hline
		2.846622 & 0.666380578 & 0 & 1.431077272 & 23.339534 & $(+)$ $\theta_p = 0.144$, $\lambda_s \approx 1.000$ & $28 \shortto 30$\\
		2.807402 & 0.644109263 & 0.159999229 & 1.416397967 & 23.764939 & $(+)$ $\theta = 2.527$, $\lambda = 1.087$ & 29\\
		2.783571 & 0.592342315 & 0.177156145 & 1.464006738 & 23.791412 & $(-)$ $\theta = 6.280$, $\lambda = 9.027$ & 29\\
		\hdashline
		2.786069 & 0.565806639 & 0.155676376 & 1.508592349 & 23.608409 & $7.718 \pm 15.53 i$ & 30\\
		2.787253 & 0.555448775 & 0.145489969 & 1.528232795 & 23.522133 & $\lambda_1 = -222.8$, $\lambda_2 = -2.971$ & 30\\
		2.787985 & 0.540861166 & 0.131006107 & 1.557274076 & 23.401380 & $\lambda = -3002$, $(+)$ $\theta = 0.002$ & 30\\
		\hdashline
		2.677285 & 0.357494019 & 0.031221157 & 2.035243544 & 22.879807 & $\lambda_1 = -2 \cdot 10^9$, $\lambda_2 = 6 \cdot 10^4$ & 29\\
		2.411372 & 0.180822175 & 0.007082813 & 2.988400837 & 22.785698 & $\lambda_1 = - 6 \cdot 10^8$, $\lambda_2 = 2 \cdot 10^7$ & 29\\
		\hline
		$C$ & $x(0)$ & $p_y(0)$ & $p_z(0)$ & $T/2$ & Floquet multipliers & $\mu_{CZ}$\\
		\hline
		2.777614 & 0.739115485 & 1.244757378 & 0.490447794 & 25.840747 & $\lambda_1 = -8 \cdot 10^9$, $\lambda_2 = -2726$ & 42\\
		2.557225 & 0.743486599 & 1.046058329 & 0.786173087 & 26.235830 & $\lambda_1 = - 10^{11}$, $\lambda_2 = -2 \cdot 10^6$ & 42
	\end{tabular}
	\label{data_13}
\end{table}

\begin{table}[ht]\fontsize{10}{10}\selectfont \centering
	\caption{Data for two families bifurcating from 9:10 DRO (orbiting around Moon, $L_1$ and $L_2$).}
	\begin{tabular}{c|c|c|c|c|c|c}
		$C$ & $x(0)$ & $z(0)$ & $p_y(0)$ & $T/4$ & Floquet multipliers & $\mu_{CZ}$\\
		\hline
		2.836798 & 1.342511949 & 0 & 0.640716906 & 13.261994 & $(-)$ $\theta_p = 4.967$, $\lambda_s \approx 1.000$ & $30 \shortto 32$\\
		2.824894 & 1.316076201 & 0.180002085 & 0.648178048 & 13.341285 & $(-)$ $\theta_1 = 5.751$, $(-)$ $\theta_2 = 6.24$ & 30\\
		2.817515 & 1.300316889 & 0.226458596 & 0.652470658 & 13.383209 & $\lambda \approx 1.000$, $(-)$ $\theta = 6.214$ & $30 \shortto 31$\\
		2.814967 & 1.294990954 & 0.240002773 & 0.653895150 & 13.396596 & $\lambda = 1.090$, $(-)$ $\theta = 6.191$ & 31\\
		2.814076 & 1.293142550 & 0.244512930 & 0.654386327 & 13.401150 & $\lambda \approx 1.000$, $(-)$ $\theta = 6.182$ & $31 \shortto 32$\\
		2.779008 & 1.225053965 & 0.370003935 & 0.670977397 & 13.539836 & $(+)$ $\theta_1 = 2.09$, $(-)$ $\theta_2 = 5.612$ & 32\\
		2.754819 & 1.177986449 & 0.421574106 & 0.677810805 & 13.582852 & $(-)$ $\theta = 5.212$, $\lambda \approx 1.000$ & $32 \shortto 33$\\
		2.753677 & 1.173247833 & 0.420731482 & 0.676869346 & 13.570586 & $(-)$ $\theta = 6.266$, $\lambda = 5.509$ & 33\\
		\hdashline
		2.768218 & 1.184087633 & 0.358462832 & 0.659294702 & 13.375902 & $\lambda_1 \approx 1.000$, $\lambda_2 = 152.8$ & $34 \shortto 33$\\
		2.772804 & 1.200832787 & 0.313572688 & 0.643817915 & 13.240235 & $(+)$ $\theta = 0.037$, $\lambda = 3 \cdot 10^4$ & 33\\
		\hdashline
		2.700924 & 1.304161893 & 0.167035003 & 0.557934058 & 13.031790 & $\lambda_1 = 10^4$, $\lambda_2 = 7 \cdot 10^8$ & 32\\
		2.210201 & 1.513826723 & 0.042797313 & 0.321464201 & 12.950987 & $\lambda_1 = 2 \cdot 10^8$, $\lambda_2 = 5 \cdot 10^{12}$ & 32\\
		\hline
		2.836798 & 0.651717320 & 0 & 1.449382820 & 13.261994 & $(-)$ $\theta_p = 4.967$, $\lambda_s \approx 1.000$ & $30 \shortto 32$\\
		2.823515 & 0.648140022 & 0.095000994 & 1.440051366 & 13.349526 & $(-)$ $\theta = 5.838$, $\lambda = 1.023$ & 31\\
		2.817782 & 0.646197466 & 0.113015599 & 1.436793082 & 13.381878 & $\lambda_1 \approx 1.000$, $\lambda_2 = 1.063$ & $31 \shortto 32$\\
		2.814485 & 0.644971948 & 0.122001100 & 1.435127303 & 13.399256 & $\lambda_1 \approx 1.095$, $\lambda_2 = 1.112$ & 32\\
		2.813427 & 0.644562140 & 0.124718629 & 1.434625033 & 13.404650 & $\lambda_1 \approx 1.000$, $\lambda_2 = 1.108$ & $32 \shortto 33$\\
		2.769414 & 0.620174249 & 0.199999406 & 1.427818754 & 13.573348 & $(+)$ $\theta = 2.499$, $\lambda = 2.405$ & 33\\
		2.741116 & 0.590155496 & 0.222942399 & 1.448470008 & 13.633918 & $\lambda_1 \approx 1.000$, $\lambda_2 = 3.772$ & $33 \shortto 34$\\
		2.736929 & 0.573367410 & 0.217860166 & 1.471230910 & 13.623612 & $\lambda_1 = 1.006$, $\lambda_2 = 23.87$ & 34\\
		\hdashline
		2.736929 & 0.573361172 & 0.217856177 & 1.471240864 & 13.623611 & $(+)$ $\theta = 0.004$, $\lambda = 23.87$ & 35\\
		2.753888 & 0.509699151 & 0.145406197 & 1.605841466 & 13.355487 & $(+)$ $\theta = 0.005$, $\lambda = 2 \cdot 10^5$ & 35\\
		\hdashline
		2.516602 & 0.246648578 & 0.034540110 & 2.514651352 & 13.068666 & $\lambda_1 = 5 \cdot 10^5$, $\lambda = 10^{15}$ & 34\\
		2.146791 & 0.084725665 & 0.007598956 & 4.352964672 & 12.871614 & $\lambda_1 = 9 \cdot 10^7$, $\lambda = 2 \cdot 10^{17}$ & 34
	\end{tabular}
	\label{data_14}
\end{table}

\begin{table}[ht]\fontsize{10}{10}\selectfont \centering
	\caption{Data for internal bridge family between two 1:1 DROs (orbiting around Moon and $L_i$, $i=1,2,4,5$).}
	\begin{tabular}{c|c|c|c|c|c|c}
		$C$ & $x(0)$ & $z(0)$ & $p_y(0)$ & $T/2$ & Floquet multipliers & $\mu_{CZ}$\\
		\hline
		2.369447 & 1.725740160 & 0 & 0.392126651 & 3.120727 & $(-)$ $\theta_p = 5.075$, $\lambda_s \approx 1.000$ & $2 \shortto 3$\\
		2.277536 & 1.721002750 & $-$0.350096729 & 0.367635240 & 3.129116 & $(-)$ $\theta_1 = 5.127$, $(-)$ $\theta_2 = 6.255$ & 2\\
		2.059431 & 1.699217701 & $-$0.660877861 & 0.310057983 & $\approx \pi$ & $(-)$ $\theta_1 = 5.228$, $(-)$ $\theta_2 = 6.256$ & 2\\
		1.737024 & 1.718190254 & $-$0.848770010 & 0.214206167 & 3.151141 & $(-)$ $\theta_1 = 5.350$, $(-)$ $\theta_2 = 6.267$ & 2\\
		1.329960 & 1.921964874 & $-$0.534165579 & 0.085878901 & 3.154626 & $(-)$ $\theta_1 = 5.450$, $(-)$ $\theta_2 = 6.274$ & 2\\
		1.177141 & 2.007218176 & 0 & 0.044057442 & 3.154815 & $(-)$ $\theta_p = 5.468$, $\lambda_s \approx 1.000$ & $3 \shortto 2$
	\end{tabular}
	\label{data_15}
\end{table}

\clearpage

\noindent
\textbf{Acknowledgement.}\ The author acknowledges support by the Deutsche Forschungsgemeinschaft (DFG, German Research Foundation), Project-ID 541062288.\ He also wishes to thank Alexander Batkhin and Otto van Koert for useful discussions on general aspects related to computations of periodic solutions.\ He is also thankful to Chankyu Joung for helpful discussions.\ He is also grateful to the anonymous referees for valuable comments and remarks that have greatly improved the scientific presentation, quality and readability of the manuscript.\

\end{document}